\documentclass[leqno,a4paper,12pt]{article}

\newtheorem{Thm}{\indent Theorem}[subsection]
\newtheorem{MainThm}[Thm]{\indent Main Theorem}
\newtheorem{Prop}[Thm]{\indent Proposition}
\newtheorem{Lem}[Thm]{\indent Lemma}

\newtheorem{Cor}[Thm]{\indent Corollary}

\newtheorem{Def}[Thm]{\indent Definition}

\newtheorem{Rem}[Thm]{\indent Remark}
\newtheorem{Ex}[Thm]{\indent Example}

\newtheorem{Conj}[Thm]{\indent Conjecture}
\newtheorem{OP}[Thm]{\indent Open Problem}

\def\qed{{\hskip0pt\unskip\unskip\nobreak\hfil\penalty50
          \hskip1em\hbox{}\nobreak\hfil
          {\bf q.e.d.}%
          \parfillskip=0pt\finalhyphendemerits=0
          \par}\medskip}

\newenvironment{Proof}
               {{\it Proof.}\quad}
               {\qed}

\newenvironment{Proofof}[1]
               {{\it Proof of #1.}\quad}
               {\qed}

\newcommand{\Prime}{\kern3\fontdimen1\font$'$\kern-7\fontdimen1\font}

\long\def\forget#1{}

\long\def\beginSIDEREMARK#1\endSIDEREMARK
    {{\par\bigskip\advance\leftskip by 2cm
                  \advance\rightskip by -2cm\noindent
      {\bf Our own side remark:} #1
      \par\bigskip\noindent}}

\long\def\beginFORGET#1\endFORGET{#1}
\long\def\beginFORGET#1\endFORGET{}

\def\?{\ ???\ \immediate\write16{}%
\immediate\write16{Warning: There was still a question mark . . . }%
\immediate\write16{}}

\usepackage{exscale}
\usepackage{amsmath} 
\usepackage{amscd}
\usepackage{amssymb}
\usepackage{cdextra} 

\newcommand{\BC}{{\mathbb{C}}}

\newcommand{\BQ}{{\mathbb{Q}}}
\newcommand{\BR}{{\mathbb{R}}}

\newcommand{\BZ}{{\mathbb{Z}}}

\newcommand{\Fa}{{\mathfrak{a}}}
\newcommand{\Fb}{{\mathfrak{b}}}

\newcommand{\FF}{{\mathfrak{F}}}
\newcommand{\FG}{{\mathfrak{G}}}

\newcommand{\CC}{{\cal C}}

\newcommand{\CF}{{\cal F}}

\newcommand{\CH}{{\cal H}}

\newcommand{\CO}{{\cal O}}

\forget{
\usepackage{calligra}
\newfont{\callignormal}{callig15 scaled 720}
\newfont{\calligscript}{callig15 scaled 500}

\let\SUB_ 
\let\SUPER^ 
\let\PRIME'

\def\MAKEIT#1#2#3#4#5#6#7#8#9{
\expandafter\edef\csname tildeC#1\endcsname%
  {\noexpand\mathchoice%
   {\mbox{\noexpand\makebox[0pt][l]{\noexpand\hskip#8
         $\noexpand\widetilde{\noexpand\phantom{t}}%
         $\noexpand\hss}}}
   {\mbox{\noexpand\makebox[0pt][l]{\noexpand\hskip#8
         $\noexpand\widetilde{\noexpand\phantom{t}}$\noexpand\hss}}}
   {\mbox{\noexpand\makebox[0pt][l]{\noexpand\hskip#9
  $\noexpand\scriptstyle\noexpand\widetilde{\noexpand\phantom{t}}%
         $\noexpand\hss}}}
   {\mbox{\noexpand\makebox[0pt][l]{\noexpand\hskip#9
  $\noexpand\scriptstyle\noexpand\widetilde{\noexpand\phantom{t}}%
         $\noexpand\hss}}}
   \csname C#1\endcsname}
\expandafter\edef\csname C#1\endcsname%
  {\noexpand\futurelet\noexpand\next\csname C#1GO\endcsname}
\expandafter\edef\csname C#1GO\endcsname%
  {\noexpand\ifx\noexpand\next\SUB
   \noexpand\let\noexpand\next\csname C#1b\endcsname
   \noexpand\else\noexpand\let\noexpand\next\csname C#1DO\endcsname
   \noexpand\fi\noexpand\next}
\expandafter\edef\csname C#1b\endcsname_##1%
  {\noexpand\def\noexpand\BOT{##1}
   \noexpand\futurelet\noexpand\next\csname C#1bGO\endcsname}
\expandafter\edef\csname C#1bGO\endcsname%
  {\noexpand\ifx\noexpand\next\noexpand\SUPER
   \noexpand\let\noexpand\next\csname C#1buDO\endcsname
   \noexpand\else\noexpand\ifx\noexpand\next\noexpand\PRIME
   \noexpand\let\noexpand\next\csname C#1bpDO\endcsname
   \noexpand\else\noexpand\let\noexpand\next\csname C#1bDO\endcsname
   \noexpand\fi\noexpand\fi\noexpand\next}
\expandafter\edef\csname C#1buDO\endcsname^##1%
  {\csname C#1DO\endcsname%
   \csname C#1kern\endcsname_{\noexpand\BOT}%
 ^{\csname C#1backern\endcsname##1}}
\expandafter\edef\csname C#1bpDO\endcsname'%
  {\csname C#1DO\endcsname%
   \csname C#1kern\endcsname_{\noexpand\BOT}%
 ^{\csname C#1backern\endcsname\prime}}
\expandafter\edef\csname C#1bDO\endcsname%
  {\csname C#1DO\endcsname%
   \csname C#1kern\endcsname_{\noexpand\BOT}}
\expandafter\edef\csname C#1DO\endcsname%
 {\noexpand\mathchoice{\mbox{\kern#2\callignormal#1\kern#3}}
                      {\mbox{\kern#2\callignormal#1\kern#3}}
                      {\mbox{\kern#4\calligscript#1\kern#5}}
                      {\mbox{\kern#4\calligscript#1\kern#5}}}
\expandafter\edef\csname C#1kern\endcsname%
 {\noexpand\mathchoice{\kern-#6}{\kern-#6}{\kern-#7}{\kern-#7}}
\expandafter\edef\csname C#1backern\endcsname%
 {\noexpand\mathchoice{\kern#6}{\kern#6}{\kern#6}{\kern#7}}
}

\MAKEIT{A}{-0.5pt}{5.7pt}{-0.5pt}{3.4pt}{3.5pt}{2.0pt}{9.0pt}{5.5pt}
\MAKEIT{B}{-1.0pt}{4.0pt}{-1.0pt}{2.1pt}{2.0pt}{1.0pt}{7.0pt}{4.0pt}
\MAKEIT{C}{-2.0pt}{4.3pt}{-2.0pt}{2.7pt}{3.0pt}{1.5pt}{6.0pt}{3.0pt}
\MAKEIT{D}{-1.8pt}{3.0pt}{+0.5pt}{1.9pt}{1.5pt}{0.2pt}{5.0pt}{5.5pt}
\MAKEIT{E}{-2.2pt}{4.0pt}{-2.1pt}{2.3pt}{2.5pt}{1.5pt}{5.0pt}{2.5pt}
\MAKEIT{F}{-0.4pt}{6.5pt}{-0.4pt}{4.4pt}{4.5pt}{3.0pt}{7.0pt}{5.0pt}
\MAKEIT{G}{-2.0pt}{4.0pt}{-2.0pt}{2.2pt}{2.5pt}{1.0pt}{7.0pt}{4.0pt}
\MAKEIT{H}{-1.0pt}{6.1pt}{-1.0pt}{4.0pt}{4.5pt}{3.0pt}{7.0pt}{5.0pt}
\MAKEIT{I}{-0.5pt}{5.9pt}{-0.5pt}{3.9pt}{4.0pt}{2.5pt}{6.5pt}{4.5pt}
\MAKEIT{J}{+1.5pt}{6.2pt}{+1.0pt}{4.0pt}{4.0pt}{2.5pt}{6.0pt}{4.0pt}
\MAKEIT{K}{-0.8pt}{6.5pt}{-0.8pt}{4.1pt}{4.5pt}{3.0pt}{8.0pt}{5.5pt}
\MAKEIT{L}{-0.5pt}{5.2pt}{-0.8pt}{3.2pt}{3.0pt}{1.5pt}{7.0pt}{4.0pt}
\MAKEIT{M}{-0.3pt}{4.8pt}{-0.5pt}{2.8pt}{3.0pt}{1.5pt}{8.5pt}{5.5pt}
\MAKEIT{N}{-0.5pt}{6.4pt}{-0.9pt}{4.0pt}{5.0pt}{3.0pt}{9.0pt}{6.0pt}
\MAKEIT{O}{-0.0pt}{4.2pt}{-1.0pt}{2.9pt}{2.5pt}{1.5pt}{6.0pt}{3.0pt}
\MAKEIT{P}{-1.0pt}{4.0pt}{-1.1pt}{2.7pt}{2.0pt}{1.0pt}{6.0pt}{3.0pt}
\MAKEIT{Q}{-0.0pt}{4.2pt}{-1.0pt}{2.7pt}{2.5pt}{1.0pt}{6.0pt}{3.0pt}
\MAKEIT{R}{-1.2pt}{3.5pt}{-1.4pt}{1.7pt}{1.5pt}{0.5pt}{6.0pt}{3.0pt}
\MAKEIT{S}{-1.0pt}{5.2pt}{-1.1pt}{3.1pt}{4.0pt}{2.0pt}{7.0pt}{4.0pt}
\MAKEIT{T}{-0.5pt}{7.0pt}{-0.7pt}{4.7pt}{5.0pt}{3.5pt}{7.0pt}{4.0pt}
\MAKEIT{U}{-2.0pt}{4.5pt}{-2.2pt}{2.5pt}{2.5pt}{1.0pt}{6.0pt}{3.0pt}
\MAKEIT{V}{-2.0pt}{6.6pt}{-2.2pt}{4.2pt}{5.0pt}{3.5pt}{7.0pt}{4.0pt}
\MAKEIT{W}{-2.0pt}{6.5pt}{-2.3pt}{4.0pt}{5.0pt}{3.5pt}{8.0pt}{5.0pt}
\MAKEIT{X}{-0.2pt}{6.3pt}{-0.5pt}{4.0pt}{4.0pt}{2.5pt}{7.0pt}{4.0pt}
\MAKEIT{Y}{-2.0pt}{4.5pt}{-1.9pt}{2.5pt}{2.5pt}{1.0pt}{5.0pt}{2.5pt}
\MAKEIT{Z}{-1.0pt}{4.3pt}{-1.1pt}{2.4pt}{3.0pt}{1.5pt}{6.0pt}{3.0pt}
}

\newcommand{\Spec}{\mathop{{\bf Spec}}\nolimits}

\newcommand{\SL}{{\rm SL}}

\newcommand{\Gr}{{\rm Gr}}
\newcommand{\image}{\mathop{{\rm Im}}\nolimits}
\newcommand{\imm}{\mathop{{\rm im}}\nolimits}
\newcommand{\kerr}{\mathop{{\rm ker}}\nolimits}
\newcommand{\cokerr}{\mathop{{\rm coker}}\nolimits}

\newcommand{\Sym}{\mathop{\rm Sym}\nolimits}
\newcommand{\Lie}{\mathop{\rm Lie}\nolimits}
\newcommand{\End}{\mathop{\rm End}\nolimits}
\newcommand{\Hom}{\mathop{\rm Hom}\nolimits}

\newcommand{\Ext}{\mathop{\rm Ext}\nolimits}

\newcommand{\UHom}{\mathop{\underline{\rm Hom}}\nolimits}

\newcommand{\Res}{\mathop{\rm Res}\nolimits}

\newcommand{\codim}{{\rm codim}}

\newcommand{\tor}{{\rm tors}}

\newcommand{\uGh}{\underline{\Gh}}
\newcommand{\bE}{\bar{E}}
\newcommand{\bF}{\bar{F}}
\newcommand{\bH}{\bar{H}}
\newcommand{\tE}{\tilde{E}}
\newcommand{\tH}{\tilde{H}}

\newcommand{\pr}{\mathop{\rm pr}\nolimits}

\def\tei{\, | \,}
\def\halb{\frac{1}{2}}

\def\phi{\varphi}
\def\epsilon{\varepsilon}

\def\id{{\rm id}}

\newcommand{\geteilt}[2]{
   \[\begin{array}{p{6cm}|p{6cm}}
   \begin{minipage}{6cm}#1\end{minipage}&
   \begin{minipage}{6cm}#2\end{minipage}
\end{array}\]
}                   

\newbox\mybox
\def\arrover#1{\mathrel{
       \setbox\mybox=\hbox spread 1.4em{\hfil$\scriptstyle#1$\hfil}
       \vbox{\offinterlineskip\copy\mybox
             \hbox to\wd\mybox{\rightarrowfill}}}}
\def\larrover#1{\mathrel{
       \setbox\mybox=\hbox spread 1.4em{\hfil$\scriptstyle#1$\hfil}
       \vbox{\offinterlineskip\copy\mybox
             \hbox to\wd\mybox{\leftarrowfill}}}}

\def\ontoover#1{\mathrel{
       \setbox\mybox=\hbox spread 1.4em{\hfil$\scriptstyle#1$\hfil}
       \vbox{\offinterlineskip\copy\mybox
             \hbox to\wd\mybox{\rightarrowfill\hskip-2.8mm
                               $\rightarrow$}}}}
\def\leftontoover#1{\mathrel{
       \setbox\mybox=\hbox spread 1.4em{\hfil$\scriptstyle#1$\hfil}
       \vbox{\offinterlineskip\copy\mybox
             \hbox to\wd\mybox{$\leftarrow$\hskip-2.8mm
                               \leftarrowfill}}}}
\def\longto{\longrightarrow}

\def\longonto{\ontoover{\ }}
\def\isoto{\arrover{\sim}}

\def\longinto{\lhook\joinrel\longrightarrow}
\def\longleftinto{\longleftarrow\joinrel\rhook}

\usepackage[curve,matrix,arrow,cmtip]{xy}
\NoComputerModernTips

\def\myxymessage{\def\messagetext
   {Here an xy-pic diagram was omitted to speed up compilation . . . }
   \immediate\write16{\messagetext}
   \hbox{\bf \messagetext}}
\def\filxymatrix#1{\myxymessage}
\def\filxyarray#1{\myxymessage}

\newdir^{ (}{{}*!/-3pt/\dir^{(}}    
\newdir_{ (}{{}*!/-3pt/\dir_{(}}    
\newdir^{ )}{{}*!/+3pt/\dir^{)}}    
\newdir_{ )}{{}*!/+3pt/\dir_{)}}    

\def\rscript#1{\hbox to 0pt{$\scriptstyle#1$\hss}}

\newcommand{\abs}{{\rm abs}}

\newcommand{\can}{{\rm can}}

\newcommand{\dR}{{\rm dR}}

\newcommand{\symm}{{\rm sym}}

\newcommand{\topp}{\rm{top}}

\renewcommand{\Bbb}{\mathbb}  

\newcommand{\C}{{\Bbb{C}}}  

\newcommand{\E}{\Bbb{E}}

\newcommand{\Na}{{\Bbb{N}}} 
\newcommand{\Q}{{\Bbb{Q}}}  
\newcommand{\R}{{\Bbb{R}}}  
\newcommand{\V}{\Bbb{V}}

\newcommand{\Z}{{\Bbb{Z}}}  

\newcommand{\Ch}{{\cal C}}
\newcommand{\Dh}{{\cal D}}
\newcommand{\Eh}{{\cal E}}
\newcommand{\Fh}{{\cal F}}
\newcommand{\Gh}{{\cal G}}
\newcommand{\Hh}{{\cal H}}
\newcommand{\Mh}{{\cal M}}
\newcommand{\Oh}{{\cal O}} 
\newcommand{\Rh}{{\cal R}}
\newcommand{\Rr}{{\cal R}}

\newcommand{\tEh}{\widetilde{\Eh}}

\newcommand{\tpi}{\widetilde{\pi}}

\newcommand{\fH}{\mathfrak{H}}

\newcommand{\fS}{\mathfrak{S}}
\newcommand{\eA}{\mathfrak{U}}

\newcommand{\bS}{{\bf S}}
\newcommand{\bT}{{\bf T}}

\newcommand{\DD}{\displaystyle}

\newcommand{\stern}{(\ast)}

\newcommand{\verk}{\circ}

\newcommand{\pfeil}{\longrightarrow}

\newcommand{\ld}{\DD \lim_{\DD\longrightarrow}}
\newcommand{\li}[1]{\DD \lim_{\stackrel{\DD\longleftarrow}{#1}}}
\newcommand{\ls}[1]{\DD \lim_{\stackrel{\DD\longrightarrow}{#1}}}
\newcommand{\silo}{\isoto}

\newcommand{\Bl}{\mathop{\rm Bl}\nolimits}
\newcommand{\DXQ}{D^b_{(\bS,L),m} (X, \Q_\ell)}
\newcommand{\ev}{{\rm ev}}

\newcommand{\Gen}{{\cal G}en}

\newcommand{\Eis}{{\cal E}is}
\newcommand{\Et}{\mathop{\bf Et}\nolimits}

\newcommand{\Log}{{\cal L}og}
\newcommand{\MHM}{\mathop{\bf MHM}\nolimits}

\newcommand{\MH}{\mathop{\bf MHM}\nolimits}

\newcommand{\MHS}{\mathop{\bf MHS}_F\nolimits}
\newcommand{\MHSR}{\mathop{\bf MHS}_{\R}\nolimits}
\newcommand{\MM}{\mathop{\bf MM}\nolimits}
\newcommand{\MMs}{\MM^s\nolimits}
\newcommand{\Ob}{{\rm Ob}}
\newcommand{\Perv}{\mathop{\bf Perv}\nolimits}
\newcommand{\pol}{pol}

\newcommand{\reg}{{\rm reg}}
\newcommand{\res}{{\rm res}}  

\newcommand{\rtop}{\topp}
\newcommand{\sgn}{{\rm sgn}}        
\newcommand{\Sgn}{{\rm Sgn}}

\newcommand{\Sh}{\mathop{\bf Sh}\nolimits}
\newcommand{\SHs}{\Sh^s}

\newcommand{\sing}{{\rm sing}}
\newcommand{\spec}{\Spec}   
\newcommand{\tr}{\operatorname{N}}
\newcommand{\uchi}{\underline{\chi}}
\newcommand{\uHom}{\UHom}

\newcommand{\Var}{\mathop{\bf Var}\nolimits}
\newcommand{\Xt}{X_{\rtop}}
\newcommand{\SXt}{\Sh \Xt}
\newcommand{\SsXt}{\SHs \! \Xt}
\newcommand{\SB}{\Sh B}
\newcommand{\SX}{\Sh X}
\newcommand{\SsB}{\Sh^s \! B}

\newcommand{\SsX}{\Sh^s \! X}
\newcommand{\SsY}{\Sh^s Y}

\newcommand{\SsWB}{\Sh^{s,W} \! B}
\newcommand{\SsWX}{\Sh^{s,W} \! X}
\newcommand{\USs}{U\!\Sh^s}

\newcommand{\USsX}{U\!\SsX}
\newcommand{\pUSsE}{\pi{\rm -}\USs \Eh}
\newcommand{\pUSstE}{\pi{\rm -}\USs \tEh}
\newcommand{\pUSsX}{\pi{\rm -}\USsX}
\newcommand{\pUSsXt}{\pi{\rm -}\USs\!\Xt}
\newcommand{\proUSsE}{pro{\rm -}\pUSsE}
\newcommand{\proUSsX}{pro{\rm -}\pUSsX}
\newcommand{\USsW}{U\!\Sh^{s,W}}
\newcommand{\pUSsW}{\pi{\rm -}\USsW}

\newcommand{\pUSsWX}{\pUSsW \! X}

\newcommand{\DSX}{D^b \SX}

\renewcommand{\bar}{\overline}

\let\oldbullet\bullet
\def\bullet{{\mathchoice{\oldbullet}%
                        {\oldbullet}%
                        {\scriptscriptstyle\oldbullet}%
                        {\oldbullet}}}

\newcommand{\argdot}{{\;\bullet\;}}
\newcommand{\argstern}{{\;\ast\;}}
\newcommand{\punkt}{{\bullet}}
\newcommand{\cpunkt}{{\bullet}}

\begin{document}  

%

\hfuzz=3pt
\overfullrule=10pt                   


\setlength{\abovedisplayskip}{6.0pt plus 3.0pt}
\setlength{\belowdisplayskip}{6.0pt plus 3.0pt}     
\setlength{\abovedisplayshortskip}{6.0pt plus 3.0pt}
\setlength{\belowdisplayshortskip}{6.0pt plus 3.0pt}

\setlength{\baselineskip}{13.0pt}
\setlength{\lineskip}{0.0pt}
\setlength{\lineskiplimit}{0.0pt}

%
%

\title{On the Eisenstein symbol
\forget{
\footnotemark
\footnotetext{To appear in J.\ of Alg.\ Geom.}
}
}
\author{\footnotesize by\\ \\
\mbox{\hskip-2cm
\begin{minipage}{7cm} \begin{center} \begin{tabular}{c}
J\"org Wildeshaus\\[5pt]
\footnotesize Institut Galil\'ee\\[-3pt]
\footnotesize Universit\'e Paris 13\\[-3pt]
\footnotesize Avenue Jean-Baptiste Cl\'ement\\[-3pt]
\footnotesize F-93430 Villetaneuse\\[-3pt]
\footnotesize France\\
{\footnotesize \tt wildesh@zeus.math.univ-paris13.fr}
\end{tabular} \end{center} \end{minipage}
\hskip-2cm} 
\\[2cm]
}

\date{June 5, 2000} 

\maketitle

\begin{abstract}
The main purpose of this paper is the geometric construction, and the
analysis of the formalism of elliptic Bloch groups.
In the setting of absolute cohomology, we obtain a generalization of
Beilinson's Eisenstein symbol to 
divisors of an elliptic curve $\Eh$,
whose support is not necessarily torsion. 
For motivic cohomology, such a generalization
is obtained in low degrees. Our main result shows that the Eisenstein
symbol can be defined in all degrees if the groups
\[
H^i_\Mh(\Eh^m , \BQ(m))_{- , \ldots , -} 
\]
vanish in a certain range of indices. Consequently, 
the weak version of Zagier's conjecture for elliptic curves
is implied by the elliptic analogue of the Beilinson--Soul\'e
vanishing conjecture.\\

\end{abstract}

\vfill


\eject
\quad
\newpage

%
%

\setcounter{section}{-1}
\section{Introduction} \label{Intro}

According to the {\it Beilinson conjectures}, there is an intimate connection
between the $K$-theory, or {\it motivic cohomology}, of a smooth proper variety 
$X$ over a number field, and special values of the $L$-function of $X$. 
These values are expected
to compare two $\Q$-structures on the highest exterior power of 
a certain finite dimensional $\R$-vector space, namely the 
{\it absolute Hodge cohomology} of $X$. The first comes about via an exact 
sequence involving 
algebraic de Rham, and Betti cohomology of $X$. According to the conjecture,
the image of $K$-theory under the so-called {\it regulator} provides 
absolute Hodge 
cohomology with a second $\Q$-structure.\\

In any of those higher dimensional cases, when weak versions of the conjectures
are proved (see \cite{N}, Section 8 for a survey), one constructs {\it special
elements} in motivic cohomology and shows then that the regulators of these 
elements generate a $\Q$-structure, whose 
covolume has the expected relation to the special
value of the $L$-function. However, one does not even know if the motivic
cohomology is finite dimensional, or whether it is generated by the special
elements.\\

The only case when the full conjectures are known concerns varieties of 
dimension zero (\cite{Bo}, \cite{R}). These are basically the 
only varieties for which the ranks of all $K$-groups 
have been determined, and the regulator is known to be injective.\\

The present work focuses on
the construction of special elements in the motivic cohomology of symmetric
powers of elliptic curves, and their images under the regulators. 
More precisely, let
\[
\pi: \Eh \longrightarrow B
\]
be an elliptic curve over a regular base $B$. Denote by $\widetilde{\Eh}$
the complement of the zero section, and by $\Eh^{(n)}$ the kernel of the
summation map
\[
\Eh^{n+1} \pfeil \Eh \; .
\] 
Note that the symmetric group $\fS_{n+1}$ acts on $\Eh^{(n)}$.
For technical reasons (to be explained
later in this introduction), we shall
suppose that $\Eh$ satisfies the
{\it disjointness property} $(DP)$: any two unequal
sections in $\Eh(B)$ have disjoint support.\\ 

The \emph{first} objective
of this work is the \emph{geometric} construction of the
following data, which we refer to as the \emph{formalism of elliptic
Bloch groups} (see Section~\ref{II}):
\begin{itemize}
\item[(a)] The \emph{elliptic Bloch groups} $\Bl_{k,\Mh} = \Bl_{k,\Mh} (\Eh)$,
$k \ge 1$, which are actually $\BQ$-vector spaces.
\item[(b)] The \emph{elliptic symbols}
\[
\{ \argdot \}_k : \widetilde{\Eh} (B) \longrightarrow
\Bl_{k,\Mh} \; .
\]
\item[(c)] The \emph{differentials}
\[
d_k : \Bl_{k,\Mh} \longrightarrow \Bl_{k-1,\Mh}
\otimes_{\Q} \Bl_{1,\Mh} \quad (:= 0 \; \mbox{for} \; k=1)    
\]
mapping $\{ s \}_k$ to $\{ s \}_{k-1} \otimes \{ s \}_1$.
\item[(d)] The \emph{restrictions}
\[
\varrho_k: H^{k-1}_{\Mh} ( \Eh^{(k-2)} , \BQ(k-1) )^\sgn \pfeil \kerr (d_k)
\]
for $k \ge 2$ (the superscript $\sgn$ refers to the action of $\fS_{k-1}$), and
\[
\varrho_1: H^2_{\Mh} ( \Eh , \BQ(1))_- = \Eh(B) \otimes_\BZ \BQ \pfeil 
\kerr (d_1) = \Bl_{1,\Mh} \; .
\]
\end{itemize}

Our \emph{second} aim is a \emph{sheaf theoretical} interpretation
of the analogue of this formalism in absolute cohomology.
In particular (Theorem~\ref{9bE}), this interpretation will allow to describe 
the images of the elements of $\kerr (d_k)$
under the regulator to absolute Hodge cohomology 
in terms of {\it Eisenstein--Kronecker
series}.\\

The \emph{third} objective is a detailed analysis of the maps
$\varrho_k$. Our Main Theorem~\ref{9bI}
implies that the $\varrho_k$ are isomorphisms
as soon as the elliptic analogue of the Beilinson--Soul\'e vanishing conjecture
(Conjecture~\ref{9bL}) holds. As a consequence
(Theorem~\ref{9cE}), we shall see that
Conjecture~\ref{9bL} 
implies Parts~1 and 2 of the \emph{weak version of
Zagier's conjecture for elliptic curves} (\cite{W5}, Conj.~1.6.(B)).\\

For small $k$, one can say more: First, we get a proof of 
the conjecture
for $k=2$ (see Subsection~\ref{10}), thus constructing elements in
\[
H^1_{\Mh} (B, \BQ(1)) = \Oh^*(B) \otimes_{\Z} \Q
\]
from certain formal linear combinations of elements of the Mordell--Weil group
of $\Eh$. The present construction is of \emph{geometric} 
nature. It will be shown to yield the same result as the \emph{sheaf 
theoretical} construction of \cite{W6}.\\

Second, in the case $k=3$, we get an unconditional proof of 
Parts~1 and 2 of the weak version of
Zagier's conjecture for elliptic curves (see Subsection~\ref{11}). 
Let us remark that in the case $B = \spec K$, the
main result of \cite{GL} is considerably stronger than ours.
Indeed, the authors show Zagier's conjecture for 
($k=3$ and) an elliptic curve $\Eh$ over an
{\it arbitrary} field $K$. Translated into the language employed above,
the difference between Zagier's conjecture proved in \cite{GL},
and its weak version considered in the present work,
consists in \emph{surjectivity}:
The Bloch group $B_3(\Eh)$ of loc.~cit.\ is
generated by the elliptic symbols if $K$ is algebraically closed.\\

In the case of $B = \spec K$, \cite{G} contains 
the construction of certain
complexes $B(\Eh,k)^{\bullet}$ for arbitrary $k \ge 2$. The groups 
occurring in $B(\Eh,k)^{\bullet}$
are given, at least up to torsion, by taking the $K$-rational divisors on
$\Eh$, and dividing out certain convolution products of divisors of functions
on $\Eh$.
Goncharov conjectures (loc.\ cit., Conjecture 9.5.b)) that there are canonical
isomorphisms
$$H^i(B(\Eh,k)^{\bullet} \otimes_{\Z} \Q) \silo 
H^{k+i-2}_{\Mh} (\Eh^{(k-2)} , \BQ(k-1))^{\sgn} $$
for any $i$.
(The main result of \cite{GL} concerns the cohomology of the complex
$B(\Eh,3)^{\bullet}$.)
It would naturally be desirable to compare the two approaches (see
Remark~\ref{9bO}).\\

We also get an alternative description of the \emph{Eisenstein symbol 
on torsion} (\cite{B5}, \cite{De1}). More precisely, we show, in 
Subsections~\ref{8} and \ref{12}, that the restriction of the formalism
of elliptic Bloch groups to torsion yields a construction, which
at least in absolute cohomology is identical, up to scaling, 
to the Eisenstein symbol of 
loc.\ cit. It should be
possible to link the two motivic constructions as well.\\

In analogy to the torsion case, let us refer to the inverse of $\varrho_k$
as the \emph{Eisenstein symbol}, whenever that inverse exists.\\

In order to sketch some of the background, 
let us continue by discussing the classical {\it Zagier conjecture} (\cite{Z}).
It gives a description, in terms of a formalism of classical Bloch groups,
of the $K$-theory of a field.
In the case of number fields,
there are two proofs of its {\it weak version}, which is about the
existence of the formalism, but does not say anything about its
image in $K$-theory:
\cite{Jeu} gives the construction of the 
$(\Bl_{k,\Mh},\{ \argdot \}_k,d_k,\varrho_k)$, in the
$K$-theory of \emph{any} field. As in \cite{Es} or \cite{De1},
quite some
technical efforts have to be made in order to identify the image of the 
elements in $\kerr(d_k)$ under
the regulator to absolute Hodge cohomology.\\

The proof given in the unpublished 
preprint \cite{BD1} is different from \cite{Jeu} in that it makes use of the 
full force of the theory of the {\it classical polylogarithm}. In particular,
one works in certain categories of mixed sheaves.
E.g., in its complex analytic incarnation, the polylog is a certain extension
of variations of Hodge structure.
One of the main advantages of this approach is that the 
complicated computations in
absolute Hodge cohomology are avoided. Instead, one needs to know the explicit 
shape of the polylog, say in
terms of the entries of its period matrix. Furthermore, one can 
treat regulators to other absolute cohomology theories in a completely 
analogous manner.\\

It should be noted that in both proofs of the weak version of Zagier's
conjecture for number fields, Borel's theorem plays a central technical 
role. In fact,
our Main Theorem~\ref{9bI} is the precise analogue of de Jeu's main
result (\cite{Jeu}, Thm.~3.12). In the number field case, he is able to
deduce the weak version of Zagier's conjecture,
because Borel's theorem implies the 
Beilinson--Soul\'e conjecture.\\ 

In the present paper, we follow the approach of \cite{BD1}, transferring
the geometric construction of loc.\ cit.\ 
to the elliptic setting, and replacing the classical
by the {\it elliptic polylogarithm}. 
As far as the explicit description of the objects is concerned,
our emphasis lies on the Hodge theoretic side. 
Let us stress however that
the formal setting to be created here is
equally applicable to the $\ell$-adic setting.
Given the recent results of Kings concerning 
the identification of the $\ell$-adic elliptic polylog 
(\cite{Ki2}, Section 4, 
in particular Thm.~4.1.3), 
a statement about the explicit 
shape of the $\ell$-adic Eisenstein symbol
seems within reach (see loc.~cit., Thm.~4.2.9, which identifies
the Eisenstein symbol on torsion points). In this
context, we expect a satisfactory treatment of Part 3 of the weak version of
the elliptic Zagier conjecture, which concerns the
{\it integrality criterion} (see \cite{W5}, 1.6.(B)).\\

Let us now describe 
how the strategy of the construction of explicit elements in motivic and 
absolute cohomology ``via polylogarithms'' looks
like. This strategy is again due to
Beilinson, and clearly present in papers like \cite{B1}, or \cite{BL}.\\

The zeroeth step is to introduce, as in \cite{BD1}, the concept of Yoneda
extensions in categories of mixed sheaves in order to reinterpret the target
of the regulator, i.e., absolute cohomology.\\

According to the motivic, i.e., largely conjectural 
philosophy, the source of the regulator 
admits a similar interpretation: there should be a $\Q$-linear 
ca\-te\-go\-ry of 
{\it mixed motivic sheaves} $\MM X$ on $X$, together with exact functors, 
called
{\it realizations}, for any absolute cohomology theory, into the 
corresponding
mixed sheaf category. There should be an isomorphism between motivic cohomology
and the Ext groups of Tate twists in $\MM X$
identifying the regulator with the morphism induced by the realization.\\

The polylogarithmic
strategy for the construction of special elements in motivic
cohomology, and the computation of their images in absolute cohomology
under the regulator, continues as follows:\\

\noindent Step 1. Construct the special elements in absolute cohomology
first. More precisely, construct the polylogarithmic extension on 
$\widetilde{\Eh}$, and
define the special elements by pulling back this extension via sections
$s: B \longrightarrow \widetilde{\Eh}$.
\begin{enumerate}
\item[(a)]
The construction is a priori sheaf theoretical, and uses concepts like
Leray spectral sequences. The polylog is characterized by a certain
universal property, one consequence of which is the so-called
{\it rigidity principle}. E.g., in the Hodge theoretic setting, rigidity 
states that the polylogarithmic 
extension of variations of Hodge structure is uniquely 
determined by the underlying extension of local systems. 
\item[(b)]
Via rigidity, it is possible to give explicit descriptions of the objects 
defined by
abstract theory.
\item[(c)]
Again via rigidity, it is possible to show that the abstract construction of
(a) is ``geometrically motivated'': the extensions occur as cohomology objects, 
with Tate coefficients, of certain
$B$-schemes.
\end{enumerate}
Step 2. Because of the present non-availability of a sheaf 
theoretical machinery
on the level of motives, step 1~(a) admits no straightforward motivic 
imitation. 
However, it
turns out that 1~(c) has a translation to motivic cohomology, yielding the
explicit elements we were looking for.
Compatibility with the construction on the level of realizations under the 
regulator is then a consequence of the definition.\\

In fact, one may expect this strategy to work in a context considerably more 
general than that of elliptic curves. We continue by indicating how, in the 
elliptic case, these steps are treated in the existing literature.\\

Step $0$ consists of the identification of absolute cohomology with
Ext groups in suitable categories of mixed sheaves. For the Hodge 
theoretic and
$\ell$-adic setting, this is done
in the appendices of \cite{HW}, in particular, loc.\ cit.,
A.1.10, A.2.7, B.4.\\

As for step 1, we refer to \cite{BL}, 1.3, or \cite{W2}, Section 1 for the 
construction of the sheaf theoretical incarnations of the elliptic
polylogarithm.
The explicit description in the setting of real Hodge structures is contained
in \cite{BL}, 3.3--3.6, and also in 
\cite{W4}, Section~4. It is recalled
in the present Subsection~\ref{2}. A specific geometric realization of the 
polylog
in absolute cohomology is given in \cite{BL}, 6.1.\\

Finally, loc.\ cit., 6.2--6.3 contains the construction of the class of the 
polylogarithm in motivic cohomology, corresponding to the geometric situation
established in loc.\ cit., 6.1.\\

While this geometric situation enables one to give a polylogarithmic
construction of the Eisenstein symbol on torsion (loc.\ cit., 6.4),
we found ourselves unable to use it for our purposes.\\ 

As a consequence,
a good part of this work
deals with an alternative solution of steps 1~(c), and hence 2
(which nonetheless in spirit follows faithfully the strategy of \cite{BL}).
The starting point is the so-called {\it residue spectral sequence}
associated to a 
stratification of a smooth scheme $X$ induced by a given $NC$-divisor.
It converges to motivic cohomology of the full scheme $X$. Its $E_1$-terms
are given by motivic cohomology of the locally closed subschemes associated
to the filtration by {\it open} subschemes, which comes about by considering
the complements of intersections of the divisor. 
In our situation, $X$ will be a power of the elliptic curve $\Eh$, and the
divisor will be associated to finitely many sections of the structural
morphism $\Eh \to B$. \\

The residue spectral sequence will be used twice:
\begin{enumerate}
\item[(A)] For the construction of the polylogarithm, we perform base
change via $\widetilde{\Eh} \to B$. We are able to control
the differentials of (a direct summand of) the spectral sequence. Consequently,
it is possible to define the polylog by an inductive process.
\item[(B)] The Bloch groups, the differentials, and the restrictions all
occur in a quotient of the residue spectral sequence. The elliptic symbol
$\{ s \}_k$ will come about as pullback of the polylogarithm by
\[
s: B \longto \widetilde{\Eh} \; .
\]
\end{enumerate}
(A) is a direct translation of the strategy used in \cite{BL}, 6.3.
The bulk of the present article consists of the detailed treatment of (B).\\

Let us comment our hypothesis $(DP)$:
Under one central aspect, our treatment of the elliptic polylogarithm is
less demanding than that of the classical situation: we do not have to work 
with simplicial schemes, and hence, do not need relative motivic cohomology 
(\cite{HW}, Appendix B). By contrast, if we intended to work in a situation 
where
sections were allowed to meet, simplicial techniques would enter at once.
As was pointed out to the author, $(DP)$ is a \emph{very} serious
restriction of generality: if the base $B$ contains a point whose residue
field is finite, then it implies that $\Eh(B)$ is torsion (in which
case none of the results mentioned in this article is new).\\

The coarse structure of this work is as follows: Section~\ref{I}
concerns the construction of the elliptic polylogarithm, Section~\ref{II}
establishes the formalism of elliptic Bloch groups and states the main
results mentioned above. Their proofs are presented in Section~\ref{III}.
For the description of the finer structure, we refer to the introductions
of the individual sections. \\

We need to mention that this division reflects 
the history of this article: Section~\ref{I} is a revision of 
the first half (Sections 1--6)
of the manuscript ``On the generalized 
Eisenstein symbol'', distributed in Spring 1997. By contrast, Sections~\ref{II}
and \ref{III} have almost no nontrivial intersection with the second
half of that manuscript. This is largely due to the author's discovery that one
of the central ideas of de Jeu's construction of the Bloch groups in the
setting of the classical Zagier conjecture (\cite{Jeu}) could be 
translated almost {\it verbatim} into the elliptic context, thus
solving the main technical problem encountered in the older version of
this work. In particular, we wish to point out that the main results no longer 
require an injectivity statement concerning the (product of the)
regulators. \\

\noindent {\bf Acknowledgements:} The present paper is a largely revised
version of
my 1997 {\it Habilitationsschrift}. It is an honour to have it included
in the proceedings of the conference commemorating the anniversary of
\cite{Bl}, which not only addressed the specific question of studying special
elements in $H^2_{\Mh} (\Eh , \BQ(2))$, but formed one of 
the motivations for Beilinson's conjectures in general. I am indebted
to the organizers for providing me with this opportunity.\\

I am grateful to the Centre Emile Borel, Paris, and
to the University of Chicago for their hospitality. During my stays in 
March 1997 and May 1999 respectively, most of the ideas for this work were put 
into a definitive form.\\

I wish to thank 
A.~Huber, F.~Lecomte, A.~Werner,
A.~Beilinson, S.~Bloch, A.~Deitmar, C.~Deninger, G.~Frey, 
T.~Gei\ss er, A.~Goncharov, G.~Kings, 
H.~Knospe, K.~K\" unnemann, E.~Landvogt, A.~Levin, B.~Moonen, 
O.~Patashnick, G.~Powell, 
A.J.~Scholl, and T.~Wenger for useful comments.\\ 

Special thanks go to G.\ Weckermann for \LaTeX ing large parts
of my manuscript.

\newpage

\tableofcontents

\newpage

%
%

\section{The elliptic motivic polylogarithm} \label{I}

In Subsection \ref{1}, we normalize the sheaf theoretic notations used 
throughout
the whole paper. In the $\ell$-adic setting, we are working with schemes of
finite type over $\Z[\frac{1}{\ell}]$. It is straightforward, however, to 
transfer
our main results to schemes of finite type over $\Q$. With \ref{1F},
we already encounter one of our main technical tools: the
spectral sequence associated to a filtration of a smooth scheme
by open subschemes, whose 
associated stratification is one by smooth schemes. Let us refer to it
as the {\it residue spectral sequence}.\\

Subsection \ref{2} gives a quick axiomatic 
description of the elliptic logarithmic sheaf 
$\Log$, and of the (small) elliptic polylogarithmic extension $pol$. 
We also give a more detailed discussion of some aspects of the
elliptic Zagier conjecture. The universal
property \ref{2A} of $\Log$ is needed only to connect the general definition of the 
logarithmic sheaf as a solution of a representability
problem to the somewhat ad hoc, but much more geometric definition of
Subsection~\ref{5}. A reader willing to accept the result on the shape
of the Hodge theoretic incarnation of the polylogarithm 
(\ref{2I}) may therefore take the constructions in Subsections~\ref{5}
and \ref{7} as a definition of both $\Log$ and $\pol$, and view 
Subsection~\ref{2}
as an extended introduction providing background material. \\

In Subsection \ref{3}, we establish the geometric situation used thereafter. 
The 
subsection is mainly intended for easier reference. Let us note that the 
definition was inspired by the construction of \cite{BD1}, Section~4, in the
classical case.\\

In Subsection~\ref{4}, we identify the residue spectral sequence in the 
special situation of open subschemes of a power $\Eh^n$, which are
given by removing sections of the structural morphism $\Eh \to B$. 
The main results are \ref{4F} and \ref{7B}.\\

Using the results of Subsection~\ref{4}, we then proceed to give the geometric
realization of $\Log$ and of $\pol$. Theorem~\ref{5E} identifies
$\Log$ (or rather, its pullback to $\widetilde{\Eh}$) with a projective
system of relative cohomology objects with coefficients in Tate twists of
certain subschemes of powers of $\Eh \times_B \widetilde{\Eh}$ over 
$\widetilde{\Eh}$. Together with
the Leray spectral sequence, this result suggests that extensions by $\Log$
should be described as elements of the projective limit of absolute cohomology
of these schemes. Corollary~\ref{5H} makes explicit this identification for
$\pol$. \\

Now that the zeroeth and first steps of the strategy, 
outlined in the introduction, have been accomplished,
the next thing to do is to translate the construction of Subsection \ref{5} to 
the setting of motivic cohomology, and give the definition of the class of
the elliptic polylogarithm. 
Here, we encounter a typical technical difficulty: the construction
of Subsection~\ref{5} was possible partly because of the vanishing of
absolute cohomology in a certain range of indices
(see e.g.\ the proof of Theorem~\ref{5B}). The analogous vanishing results 
in motivic cohomology are
at present not available.  
In order to circumvent this technical complication, we translate 
\cite{BL}, 6.2, 6.3 to our setting, and employ the action induced by the 
isogenies $[a]$ on
our motivic cohomology groups. 
The main result of Subsection~\ref{7} is Theorem~\ref{7E}. It states that a 
certain differential in the residue spectral sequence
is an isomorphism on the generalized eigenspace for the eigenvalue $a$.
As in \cite{BL}, this enables one to define the motivic elliptic
polylogarithm by an inductive process (\ref{7G}). By construction, everything 
is compatible with what was done before under the regulators. From its 
definition, we also get a
statement about the image of the polylog under the differential of our
spectral sequence (Proposition~\ref{7I}). This result will be responsible
for the special shape of our differential
\[
d_k : \Bl_{k,\Mh}  \longrightarrow \Bl_{k-1,\Mh} \otimes_{\Z} \, \Eh(B) 
\]
on the elliptic symbols 
$\{s\}_k$. The precise structure of our motivic realization 
of $\pol$ is explained by its {\it application}, in Sections~\ref{II} and
\ref{III}, to the formalism of elliptic Bloch groups. Readers interested
only in the construction of $\pol$ will find a concise presentation
in \cite{Ki}, where spectral sequences are completely avoided by using
motivic {\it homology}. In fact, we wish to mention that after studying
\cite{Ki}, it was possible to considerably simplify the original proof of  
Theorem~\ref{7E}, by observing that for relative
elliptic curves, the techniques used in \cite{Ki} generalize to base
schemes which are not necessarily smooth and quasi-projective over a field.\\

Subsection~\ref{8} discusses the specialization of the
elliptic polylogarithm to torsion sections. We get versions of the
results of \cite{BL}, 6.4 in our geometric setting.

\newpage
\subsection{Mixed sheaves} \label{1}

We start by defining the sheaf categories which will be relevant for us.
For our purposes, it will be necessary to work in the settings of \emph{mixed
$\ell$-adic perverse sheaves}, of \emph{algebraic mixed Hodge 
modules}, and of \emph{algebraic 
mixed Hodge modules over $\R$}. Since the procedures 
are entirely 
analogous, we introduce, for economical reasons, the following 
rules: whenever an area of paper is divided by a vertical bar\\
\geteilt{\vspace{1cm}}{\hspace{0cm}}

\noindent the text on the left of it will concern the Hodge theoretic 
setting, while the text on the right will deal with the $\ell$-adic 
setting. 
We let
\geteilt
{
\begin{gather*} 
A := \R \quad \text{or} \quad A := \C \; ,\\
F := \Q \quad \text{or} \quad F := \R \; .
\end{gather*}
}
{
\begin{gather*}
\ell:= \text{a fixed prime number} \; , \\ 
 A := \Z \left[ \frac{1}{\ell} \right] \; ,\\
 F := \Q_\ell \; .
\end{gather*}
}

For any reduced, separated and flat scheme $X$ of finite type over 
$A$, we let
\geteilt
{
\begin{gather*}
\Xt := X (\C) \; \text{as a topol. space} \; ,\\
\SXt := \Perv (\Xt , F) \; ,\\
\hspace{0cm}
\end{gather*}
}
{
\begin{gather*}
\Xt := X\otimes_A \overline{\Q} \; , \\
\SXt := \Perv (\Xt , \Q_\ell) \; ,\\
\hspace{0cm}
\end{gather*}
}

\noindent the latter categories denoting the respective categories of 
perverse sheaves on $\Xt$ (\cite{BBD}, 2.2).\\

Next we define the category $\SX$: in the $\ell$-adic setting, we fix 
a pair $(\bS , L)$ consisting of a \emph{horizontal stratification} $\bS$ 
of $X$ (\cite{H2}, \S\,2) and a collection $L = \{ L (S) \tei S \in 
\bS\}$, where each $L (S)$ is a set of irreducible lisse $\ell$-adic 
sheaves on $S$. For all $S \in \bS$ and $\Fh \in L (S)$, we require 
that for the inclusion $j: S \hookrightarrow X$, 
all higher direct images $R^n j_{\ast} \Fh$  are 
$(\bS,L)$--{\it constructible}, i.e., have lisse restrictions to all $S \in 
\bS$, which are extensions of objects of $L (S)$.
We assume that all $\Fh \in L(S)$ are pure.\\

Following \cite{H2}, \S\,3, we define $\DXQ$ as the full subcategory 
of $D^b_c (X, \Q_\ell)$ of complexes with $(\bS,L)$--constructible 
cohomology objects. Note that all objects will be mixed.
By \cite{H2}, \S\,3, $\DXQ$ admits a perverse t-structure, whose heart
we denote by 
$\Perv_{(\bS,L)} (X,\Q_\ell)$.
\geteilt
{
\begin{gather*}
\SX := \MH_{F} (X / A) \; .\\
\hspace{0cm}
\end{gather*}
}
{
\begin{gather*}
 \SX := \Perv_{(\bS,L)} (X, \Q_\ell) \; .\\
\hspace{0cm} 
\end{gather*}
}

For $A = \C$, we define $\MH_{F} (X / \C)$ as $\MH_F (X)$, the category
of algebraic mixed $F$-Hodge modules on $X$ (\cite{S1}, \cite{S2}). 
For $A =\R$, we refer to \cite{HW},
A.2.4.\\

Because of the horizontality requirement in the $\ell$-adic situation 
we have the full formalism of the six Grothendieck functors only on the 
direct limit $D^b_m (\eA_X , \Q_\ell)$ 
of the $D^b_{(\bS,L),m} (X_U , \Q_\ell)$, for $U$ open in $\spec A$, and 
$(\bS,L)$
as above (see \cite{H2}, \S\,2).
However, for a fixed morphism
\[
\pi : X \longrightarrow Y \; ,
\]
we have a notion of e.g.\ \emph{$\pi_{\ast}$--admissibility} for a pair 
$(\bS,L)$: this is the case if
\[
\DXQ \longinto D^b_m (\eA_X , \Q_\ell) 
\stackrel{\pi_{\ast}}{\longrightarrow} D^b_m (\eA_Y , \Q_\ell)
\]
factors through some $D^b_{(\bT,K),m} (Y, \Q_\ell)$. Our computations
will show, at least a posteriori, that $(\bS,L)$ can be chosen such that
all functors which appear are admissible. We will not stress these
technical problems and even suppress $(\bS,L)$ in our notation. \\

As in \cite{BBD}, we denote by $\pi_{\ast} , \pi^{\ast} , \Hom$ 
etc. the respective functors on the categories
\geteilt
{
\begin{gather*}
\DSX := D^b \MH_{F} (X / A) \, ,\\ 
\hspace{0cm}
\end{gather*}
}
{
\begin{gather*} 
\DSX := \DXQ \, ,\\
\hspace{0cm}
\end{gather*}
}

\noindent and $\Hh^q$ for the (perverse) \emph{cohomology functors}.\\

We refer to objects of $\SX$ as \emph{sheaves}, and to objects of $\SXt$ 
as \emph{topological sheaves}. Let us denote by
$$\V \longrightarrow \V_{\rtop}$$
the forgetful functor from $\SX$ to $\SXt$.\\

If $X$ is smooth, we let
\geteilt
{
\begin{gather*}
\SsX := \Var_{F} (X / A) \subset \SX \; , \\
\text{the category of admissible} \\[-0.3em]
\text{variations on $X$ / on $X / \R$,}\\
\SsXt := \; \text{the category of}\\[-0.3em]
\text{$F$--local systems on $\Xt$.}\\
\hspace{0cm}
\end{gather*}
}
{
\begin{gather*}
\SsX := \Et^{l,m}_{\Q_\ell} (X) \subset \SX \; ,\\ 
\text{the category of lisse} \\[-0.3em]
\text{mixed $\Q_\ell$--sheaves on $X$,}\\
\SsXt := \; \text{the category of}\\[-0.3em]
\text{lisse $\Q_\ell$--sheaves on $\Xt$.}\\
\hspace{0cm}
\end{gather*}
}

For the definition of admissible variations over $\R$, compare 
\cite{HW}, A.2.1.
We refer to objects of $\SsX$ as {\it smooth sheaves}, and to objects of
$\SsXt$ as {\it smooth topological sheaves}.

\begin{Rem} \label{1A}
Note that in the $\ell$-adic situation, the existence of a
\emph{weight filtration}, i.e., an 
ascending filtration $W_{\bullet}$ by subsheaves indexed by the
integers, such that $\Gr_m^{W}$ is of weight $m$,
is not incorporated in the definition of $\Sh^s$ -- compare the warnings in
\cite{H2}, \S\,3. In the Hodge theoretic setting, the 
weight filtration is part of the data.
\end{Rem}

We define $\SsWX \subset \SsX$
as the full subcategory of smooth sheaves with a weight filtration.
If we use the symbol $W_{\punkt}$, it will always refer to the 
weight filtration.

\begin{Rem} \label{1B}
We have to 
deal with a shift of the index when viewing e.g.\ a variation as a 
Hodge module, which occurs either in the normalization of the 
embedding
\[
\Var_{F} (X / A) \longrightarrow D^b \MH_{F} (X / A)
\]
or in the numbering of cohomology objects of functors induced by 
morphisms between schemes of different dimension. In order to 
conform with the conventions laid down in \cite{HW}, Appendix A,
and \cite{W1}, Chapter 4, we choose the second possibility: a variation
{\em is} a Hodge module, not just a shift of one such.
Similarly, a lisse mixed $\Q_\ell$-sheaf {\em is} a perverse mixed sheaf.
Therefore, if $X$ is of pure relative dimension $d$ over $A$, then
the embedding
$$\Et^{l,m}_{\Q_\ell} (X) \longrightarrow D^b_m (\eA_X , \Q_\ell)$$
associates to $\V$ the complex concentrated in degree $-d$, whose only
nontrivial cohomology object is $\V$.
\end{Rem}

As a consequence, the numbering of cohomology objects of the direct image
(say) will differ from what the reader might be used to: e.g., the
cohomology of a curve is concentrated in degrees $-1$, $0$, and $1$
instead of $0$, $1$, and $2$. Similarly, one has to distinguish between
the ``naive'' pullback $(\pi^s)^{\ast}$ of a smooth sheaf and the 
pullback $\pi^{\ast}$ on the level of $\DSX$: $(\pi^s)^{\ast}$ lands
in the category of smooth sheaves, while $\pi^{\ast}$ of a smooth
sheaf yields only a smooth sheaf up to a shift.

\begin{Rem} \label{1Ba}
In the $\ell$-adic situation, there is of course a ``non-perverse'' 
theory of sheaves. We have chosen the perverse numbering mainly
in order to get a coherent picture in all cohomology theories,
and because mixed motivic sheaves are expected to be perverse in nature. 
Similarly, one might want to consider schemes over
$\BQ$ or, more generally, a field of finite type
over its prime field. The results
of the present article continue to hold, \emph{mutatis mutandis}, in these
settings.
\end{Rem}

In the special situation of pullbacks, we allow ourselves one notational
inconsistency: if there is no danger of confusion,
we use the notation $\pi^{\ast}$ also for the naive pullback of smooth
sheaves.
Similar remarks apply for smooth topological sheaves.\\

For a scheme $a : X \to \spec A$, we define
\[
F (n)_X := a^{\ast} F (n) \in D^b \SX \; ,
\]
where $F (n)$ is the usual \emph{Tate twist} on $\spec A$.\\

If $X$ is smooth, we also have the naive Tate twist
\[
F (n) := (a^s)^{\ast} F (n) \in \SsX \subset \SX
\]
on $X$. If $X$ is of pure dimension $d$, then we have the equality
\[
F (n) = F(n)_X [d] \; .
\]

In order to keep our notation transparent, we have the following

\begin{Def}\label{1C}
For any morphism $\pi:X \to S$ of reduced, separated and flat $A$-schemes
we let 
\begin{align*}
\Rr_S(X,\argdot) := \pi_*:&D^b\Sh(X) \longrightarrow D^b\Sh(S) \; ,\\
\Hh^i_S(X,\argdot) := \Hh^i\pi_*: &D^b\Sh(X) \longrightarrow \Sh(S) \; . 
\end{align*}
\end{Def} 

\begin{Def}\label{1D}
For a separated, reduced, flat $A$-scheme $X$ of finite type, and an object
$M^{\punkt}$ of $D^b \Sh (X)$, define
\begin{align}
 &R\Gamma_{\abs} (X , M^{\punkt}) := R \Hom_{D^b \Sh (X)} (F 
(0)_X , M^{\punkt})\; ,\tag*{(a)}\\
&H^i_{\abs} (X , M^{\punkt}) := H^i R \Gamma_{\abs} (X , 
M^{\punkt})\; ,\notag\\
\intertext{the {\em absolute complex and absolute cohomology groups of $X$ with 
coefficients in $M^{\punkt}$}.}
& R\Gamma_{\abs} (X ,n) := R \Gamma_{\abs} (X , F (n)_X)\; , \tag*{(b)}\\
&H^i_{\abs} (X,n) := H^i_{\abs} (X , F(n)_X)\;.\notag 
\end{align}
\end{Def}

\begin{Rem} \label{1E}
If $X$ is a scheme over $S$, then we have the formulae 
\begin{align}
& R\Gamma_{\abs} (X , \argdot) = R\Gamma_{\abs} \left( S , 
       \Rr_S(X,\argdot) \right) \; , \notag \\
& H^i_{\abs}(X,\argdot)=H^i_{\abs}(S,\Rr_S(X,\argdot)) \; . \notag
\end{align}
\end{Rem}

We mention explicitly one aspect of the Grothendieck formalism, which we shall 
frequently employ: If
\[
Z \stackrel{i}{\longinto} X \stackrel{j}{\longleftinto} U
\]
are immersions of a reduced closed subscheme $Z$ and its complement $U$ in a 
reduced, separated and flat $A$-scheme $X$, then there is an exact triangle
\[
\vcenter{\xymatrix@R-10pt{ 
i_\ast i^! \ar[rr] && \id \ar[dl] \\
&j_\ast j^\ast \ar[ul]^{[1]} &\\}}
\]
of functors on $D^b \SX$, which we refer to as the {\it residue triangle}. 
If furthermore both $Z$ and $X$ are smooth over $A$, and $Z$ is of pure 
codimension $c$ in $X$, then there is a canonical isomorphism
\[
i^! \silo i^{\ast} (-c) [-2c]
\]
of functors on $\SsX$, referred to as {\it purity}. \\

Successive application of the cohomological functors
$H^\punkt_\abs$ resp.\ $\Hh^\punkt_S$ on $\DSX$ to certain residue triangles 
yields the following result, central to everything that will be done
from Section~\ref{II} onwards:

\begin{Thm} \label{1F}
Let $S$ be a smooth separated $A$-scheme of finite type, $X / S$ smooth, 
separated and of finite type. Let
\[
\emptyset = F_{-(n+1)} X \subset F_{-n} X \subset F_{-(n-1)} X 
\subset \ldots \subset F_0 X = X
\]
be a filtration of $X$ by open subschemes, such that the reduced locally 
closed subschemes
\[
j_k : G_kX := F_k X - F_{k-1} X \hookrightarrow X
\]
are smooth, and of pure codimension $c_k$. ($G_k X$ is closed in $F_kX$.) 
Let $M \in \DSX$ be a shift of a smooth sheaf on $X$.
\begin{enumerate}
\item[(a)] There is a natural spectral sequence
\[
E^{p,q}_1 = H^{-2c_p + p+q}_{\abs} (G_p X , j^{\ast}_p M (-c_p)) 
\Longrightarrow H^{p+q}_{\abs} (X,M)
\]
of $F$-vector spaces.
\item[(b)] There is a natural spectral sequence
\[
E^{p,q}_1 = \Hh^{-2c_p + p+q}_S (G_p X , j^{\ast}_p M (-c_p)) 
\Longrightarrow \Hh^{p+q}_S (X,M)
\]
of sheaves on $S$.
\item[(c)] The differentials $\partial^{p,q}_1$ from
\[
H^{-2c_p + p+q} (G_p X , j^{\ast}_p M (-c_p))
\]
to
\[
H^{-2c_{p+1} + p+q-1} (G_{p+1} X , j^{\ast}_{p+1} M (-c_{p+1})) \; ,
\]
for $H^\punkt \in \{ H^\punkt_\abs , \Hh^\punkt_S \}$,
are induced by the composition of the morphism $(i_p)_{\ast} i^!_p \to 
\id_{F_p}$ of the residue triangle associated to
\[
G_p X \stackrel{i_p}{\longinto} F_p X \longleftinto F_{p-1}X \; ,
\]
and the boundary morphism of the residue triangle associated to
\[
G_{p+1} X \stackrel{i_{p+1}}{\longinto} F_{p+1} X 
\longleftinto F_p X \; .
\]
\item[(d)] The edge morphisms
\begin{eqnarray*}
(p = -n:) \quad H^{-n+q} (X,M) & \longrightarrow & 
H^{-n+q} (G_{-n} X , j_{-n}^\ast M) \; , \\
(p = 0:) \quad H^{-2c_0+q} (G_0 X , j^{\ast}_0 M (-c_0)) & \longrightarrow & 
H^{q} (X ,M)
\end{eqnarray*}
for $H^{\cpunkt} \in \{ H^{\cpunkt}_{\abs} , \Hh^{\cpunkt}_S \}$, 
are the natural restriction to an open subscheme, and the Gysin map from a 
smooth closed subscheme of pure codimension $c_0$, respectively.
\end{enumerate}
\end{Thm}

\begin{Proof} 
This follows from the theory of exact 
$\delta$-couples (see \cite{Wei}, 5.9, or Example~\ref{13bE}). 
\end{Proof}

Let us refer to the spectral sequences of \ref{1F}~(a) and (b) as the 
{\it absolute} and {\it relative residue spectral sequence} respectively.\\

As in \cite{W1}, Chapter 3, we want to talk about {\it relatively unipotent} 
smooth sheaves:

\begin{Def} \label{1G}
Let $\pi : X \to Y$ be a morphism of smooth and separated $A$-schemes of finite 
type. $\pUSsX$ is defined as the full subcategory of {\it $\pi$-unipotent} 
objects of $ \SsX$, i.e., those sheaves admitting a filtration, whose graded 
parts are pullbacks of smooth sheaves on $\SsY$. Similarly, one defines 
$\pUSsWX$, and $\pUSsXt$. 
\end{Def}

\newpage
\subsection{Logarithmic sheaf and polylogarithmic extension} \label{2}

We aim at a sheaf theoretical description of the (small) elliptic polylogarithm. 
The first step is an axiomatic definition of the {\it logarithmic pro-sheaf}. 
We need the following result:

\begin{Thm} \label{2A} 
Let $\pi : X \to Y$ be a morphism of smooth and 
separated $A$-schemes of finite type, which identifies $X$ as the complement 
in a smooth, proper $Y$-scheme of an $NC$-divisor relative to $Y$, all of 
whose irreducible components are smooth over $Y$. Let $i \in X (Y)$. 
The functor
\[
i^{\ast} : \pUSsX \longrightarrow \SHs Y
\]
is pro-representable in the following sense:
\begin{enumerate}
\item[(a)] There is an object
\[
\Gen_i \in \proUSsX \; ,
\]
the {\em generic pro-$\pi$-unipotent sheaf with base point $i$ on $X$}, which 
has a weight filtration satisfying
\[
\Gen_i / W_{-n} \Gen_i \in \pUSsX \quad \mbox{for all} \; n \; .
\]
Note that this implies that the direct system
\[
(R^0 \pi_{\ast} \uHom (\Gen_i / W_{-n} \Gen_i , \V))_{n \in \Na} 
\]
of smooth sheaves on $Y$ becomes constant for any $\V \in \pUSsX$. 
This constant value is denoted by
\[
R^0 \pi_{\ast} \uHom (\Gen_i , \V) \; .
\]
\item[(b)] There is a morphism of sheaves on $Y$
\[
1 \in \Hom_{\SsY} (F(0) , i^{\ast} \Gen_i) \; .
\]
\item[(c)] The natural transformation of functors from $\pUSsX$ to $\SsY$
\begin{eqnarray*}
\ev : R^0 \pi_{\ast} \uHom (\Gen_i , \argdot ) & 
\longrightarrow & i^{\ast} \; , \\
\varphi &\longmapsto & (i^{\ast} \varphi) (1) 
\end{eqnarray*}
is an isomorphism.

\end{enumerate}
\end{Thm}

\begin{Proof} This is 
\begin{center}
\begin{tabular}{r|l}
\cite{W1}, Rem.~d) after Thm.~3.6. & \cite{W1}, Thm.~3.5.i).
\end{tabular}
\end{center}
\end{Proof}

By applying the functor $\Hom_{\SsY} (F (0) , \argdot )$ to the result in (c),
one obtains:

\begin{Cor} \label{2B}
The natural transformation of functors
\begin{eqnarray*}
\Hom_{\pi - \USsX} (\Gen_i , \argdot ) &\longrightarrow& 
\Hom_{\SsY} (F (0) , i^{\ast} \argdot ) \; ,\\
\varphi &\longmapsto & (i^{\ast} \varphi) (1)
\end{eqnarray*}
is an isomorphism.
\end{Cor}

Now let $B$ be a smooth, separated $A$-scheme of finite type,
\[
\pi : \Eh \longrightarrow B
\]
an elliptic curve with zero section $i \in \Eh (B)$,
\begin{eqnarray*}
\tEh := \Eh - i (B) \; , \; j : \tEh & \longinto & \Eh \; , \\
\tpi := \pi \verk j : \tEh & \longrightarrow & B \; .
\end{eqnarray*}
We may form the generic pro-$\pi$-unipotent sheaf with base point $i$ on $\Eh$. 

\begin{Def} \label{2C}
$\Log := \Gen_i \in$ $\proUSsE$ is called the {\em elliptic logarithmic 
pro-sheaf}.
\end{Def}

\begin{Rem} \label{2D} 
Our definition of $\Log$ coincides with that of \cite{BL}, 1.2 
(see in particular loc.\ cit., 1.2.8). There, the logarithmic pro-sheaf 
is denoted by $G$.
\end{Rem}

\begin{Def} \label{2E}
Define
\[
\fH := \Hh^0_B (\Eh , F (1)) \; .
\]
\end{Def}

\begin{Rem} \label{2F} 
The topological sheaf underlying $\fH$ is
\[
R^1 (\pi_{\rtop})_{\ast} F (1)_{\rtop} \; .
\]
\end{Rem}

By abuse of notation, we also denote the pullback of $\fH$ to $\Eh$ or 
$\tEh$ by $\fH$.\\

We need to know $\Hh^{\bullet}_B (\tEh , \Log (1) \, |_{\tEh})$:

\begin{Thm} \label{2G}
(a) $\Hh^q_B (\tEh , \Log (1) \, |_{\tEh}) = 0$ for $q \neq 0$.\\
(b) There is a canonical isomorphism
\[
\Hh^0_B (\tEh , \Log (1) \, |_{\tEh}) \silo W_{-1} (i^\ast \Log) \; .
\]
(c) The sheaf $i^\ast \Log$ is split:
\[
i^\ast \Log = i^\ast \Gr_\bullet^W \Log \; .
\]
(d) For every $k \ge 0$, there is an isomorphism
\[
\Gr_k^W \Log \silo \Sym^k \fH \; .
\]
Consequently, there is an isomorphism
\[
\Hh^0_B (\tEh , \Log (1) \, |_{\tEh}) \silo \prod_{k \ge 1} \Sym^k \fH \; .
\]
\end{Thm}

\begin{Proof}
For (a) and (b), we refer to
\cite{BL}, 1.3.3, or \cite{W2}, Thm.~1.3, together with loc.\ cit., 
Rem.~d) at the end of Chapter 2. (c) is \cite{BL}, 1.2.10~(vi), or
\cite{W2}, Prop.~6.1. Finally, (d) is \cite{BL}, 1.2.8~(a); alternatively,
see the discussion below.
\end{Proof}

The theorem enables one to define the {\em small elliptic polylogarithmic 
extension} as the extension
\[
\pol \in \Ext^1_{\pUSstE} 
(\Gr_{-1}^W \Log \; |_{\tEh} , \Log (1) \; |_{\tEh})
\]
mapping to the natural inclusion $\fH \hookrightarrow \prod_{k \ge 1} 
\Sym^k \fH$ under the isomorphism
\begin{eqnarray*}
&& \Ext^1_{\Sh \tEh} ((\tpi^s)^{\ast} (i^s)^\ast \Gr_{-1}^W \Log \; |_{\tEh} , 
\Log (1) \; |_{\tEh}) \\  
&=& \Hom_{D^b \Sh \tEh} \left( \tpi^{\ast} 
                                       (i^s)^\ast \Gr_{-1}^W \Log \; |_{\tEh} , 
\Log (1) \, |_{\tEh} \right)\\
&\silo & \Hom_{\Sh B} \left( \fH , \prod_{k \ge 1} \Sym^k \fH \right)
\end{eqnarray*}
induced by the isomorphism of \ref{2G}. 
Note that the definition of $\pol$ is independent of the choice of the 
isomorphisms
\[
\Gr_k^W \Log \silo \Sym^k \fH \; .
\]

For a description of $\Log$ and $\pol$ in the Hodge incarnation, we refer to 
\cite{BL}, Sections 3 and 4. The reader may find it useful to also consult 
\cite{W4}, Chapters 3 and 4, setting $N = 1$ in the notation of loc.\ cit.\\

In order to obtain one-extensions of sheaves on $B$ via the machinery described 
in \cite{W5}, Section 3, we need to fix an isomorphism
\[
\kappa : \Gr^W_{\bullet} \Log \silo \prod_{n \ge 0} \Sym^n \fH \; .
\]
We use the same isomorphism as in \cite{W4}, Chapter 2, 
and in \cite{W5}. We recall the definition:\\

By \ref{2B}, there is a canonical projection
\[
\varepsilon : \Log \longrightarrow F (0) \; .
\]
Furthermore, there is a canonical isomorphism
\[
\gamma : \Gr^W_{-1} \Log \silo \pi^{\ast} \Hh^0_B (\Eh , F (0))^{\vee}
\]
given by the fact that both sides are equal to $\pi^{\ast}$ of the mixed 
structure on the (abelianized) fundamental group sheaf (see \cite{W1}, 
Chapter 2). There is an isomorphism
\[
\alpha : \Hh^0_B (\Eh , F (0))^{\vee} = \fH^{\vee} (1) \silo \fH
\]
induced by the Poincar\'e pairing
\[
\langle \argdot , \argdot \rangle : 
\Hh^0_B (\Eh ,F (0)) \otimes_F \fH = (\fH(-1)) \otimes_F \fH \longrightarrow F (0)
\]
(see below).\\

Finally, both $\Gr^W_{\bullet} \Log$ and $\prod_{k \ge 0} \Sym^k \fH$ carry a 
canonical multiplicative structure: for $\Gr^W_{\bullet} \Log$, this is a formal 
consequence of
\begin{center}
\begin{tabular}{r|l}
\cite{W1}, Cor.~3.4.ii) & \cite{W1}, Cor.~3.2.ii)
\end{tabular}
\end{center}
(see Rem.~b) at the end of Chapter 3 of loc.\ cit.).\\

Our isomorphism
\[
\kappa : \Gr^W_{\bullet} \Log \silo \prod_{n \ge 0} \Sym^n \fH
\]
is the unique isomorphism compatible with $\varepsilon , \gamma, \alpha$, and the multiplicative structure of both sides.\\

In order to fix notation, we wish to make explicit the isomorphism $\alpha$ 
in the Hodge theoretic setting when $A = \C$ and $B = \spec \C$. 
We follow the normalization of \cite{W4}, p.~331:\\

Fix an isomorphism
\[
\theta : \Eh (\C) \silo \C / L \; .
\]
Then any basis $(e_1 , e_2)$ of $L$ gives an $F$-basis of the vector space
\[
L \otimes_{\Z} F = H_1 (\Eh (\C) , F) = H^1 (\Eh (\C) , F)^{\vee}
\]
underlying $\Hh^0_B (\Eh , F (0))^{\vee}$ (and denoted by $V_2$ in loc.\ cit.).\\

Assume that $\image (e_1 / e_2) > 0$, and denote by $(e^{\vee}_1 , e^{\vee}_2)$ 
the dual basis of $H^1 (\Eh (\C) , F)$. Then
\[
\alpha : \Hh^0_B (\Eh , F (0))^{\vee} \silo \fH
\]
is given by
\begin{eqnarray*}
e_1 & \longmapsto & -2 \pi i \cdot e^{\vee}_2 \; , \\
e_2 & \longmapsto & 2 \pi i \cdot e^{\vee}_1 \; .
\end{eqnarray*}
This description is independent of the choice of $\theta$, and of the choice
of $(e_1,e_2)$.\\

When $A = \C , \; B = \spec \C$, and furthermore $F = \R$, 
we wish, following \cite{W5}, 4.2, to describe the fibres of
\[
\pol \in \Ext^1_{\pUSstE} (\fH , \Log (1) \, |_{\tEh}) \; ,
\]
i.e., its pullbacks
\[
s^{\ast} \pol \in \Ext^1 (\fH , s^{\ast} \Log (1))
\]
via points $s \in \tEh (\C)$. The extensions take place in the category
\[
\MHSR = \Sh \spec \C
\]
of mixed polarizable $\R$-Hodge structures. 
As there are no non-trivial extensions in
$\MHSR$ of $\R(0)$ by $\fH$ (see \cite{Jn3}, Lemma 9.2), we have 
\[
\xymatrix@1{ 
s^{\ast} \Log (1) = \prod_{n \ge 0} s^{\ast} \Gr^W_{-n} \Log (1)
\ar[r]^-{\sim}_-{\kappa} & \prod_{n \ge 0} \Sym^n \fH (1) \; ,}
\]
which means that we have push-out maps
\begin{eqnarray*}
\Ext^1 (\fH , s^{\ast} \Log (1)) & \pfeil 
                     & \Ext^1 (\fH , \Sym^{k-1} \fH (1)) \\
& = & \Ext^1 (\R (0) , \fH^{\vee} \otimes_{\R} \Sym^{k-1} \fH (1))
\end{eqnarray*}
for any $k \ge 1$. For $k \ge 2$, there is an epimorphism
\[
\fH^{\vee} \otimes_F \Sym^{k-1} \fH \longonto \Sym^{k-2} \fH
\]
given by $\frac{1}{k}$ times ``derivation''. In any basis 
$(\varepsilon_1 , \varepsilon_2)$ of $\fH_{\rtop}$, it sends
\[
\varepsilon^{\vee}_j \otimes f (\varepsilon_1 , \varepsilon_2) \quad \mbox{to} 
\quad \frac{1}{k} \left( \frac{\partial}{\partial  \varepsilon_j} \right) 
f (\varepsilon_1 , \varepsilon_2) \; .
\]
We end up with an element
\[
s^{\ast} \pol_k \in \Ext^1_{\MHSR} (\R (0) , \Sym^{k-2} \fH (1))
\]
for $k \ge 2$. By \cite{Jn3}, Lemma 9.2, this group is isomorphic to
\[
\Sym^{k-2} H^1 (\Eh (\C) , 2 \pi i \R) \; .
\]
We think of it as being contained in 
\[
\Sym^{k-2} H^1_{\dR} (\Eh) \; ,
\]
the symmetric power of the de Rham-cohomology of $\Eh$. Via the isomorphism
$\theta$, we get a basis of this vector space:
\[
\left( (dz)^\alpha \overline{(dz)}^\beta \tei \alpha, \beta \ge 0, 
\alpha + \beta = k-2 \right) \; .
\]
Let us recall the following definition:

\begin{Def} \label{2H}
Denote by $covol(L)$ the covolume of the lattice $L$ in $\C$, and by
$\langle \argdot,\argdot \rangle_L$ the Pontryagin pairing between $\C / L$ and $L$.
\[
G_{\Eh , k} : \tEh (\C) \longrightarrow 
\Sym^{k-2} H^1_{\dR} (\Eh) 
\]
is defined as the map associating to $s \in \tEh (\CH)$ the element
\[
\frac{covol(L)}{2 \pi} \sum_{\alpha + \beta = k-2} 
(dz)^\alpha \overline{(dz)}^\beta 
\sum_{\gamma \in L-\{0\}} 
   \frac{\langle \theta(s),\gamma \rangle_L}
                    {\gamma^{\alpha+1} \bar{\gamma}^{\beta+1} }
\in \Sym^{k-2} H^1_{\dR} (\Eh) \; . 
\]
\end{Def}

It is known (see e.g.\ \cite{W5}, Prop.~1.3) that the map $G_{\Eh , k}$ lands
in fact in the subgroup $\Sym^{k-2} H^1 (\Eh (\C) , 2 \pi i \R)$. 
The connection to the
elliptic polylogarithm is given by the following result:

\begin{Thm} \label{2I}
For any elliptic curve $\Eh / \C$ and any point $s \in \tEh (\C)$, we have for 
$k \ge 2$:
\[
s^{\ast} \pol_k = G_{\Eh,k} (s) \in \Sym^{k-2} H^1 (\Eh (\C) , 2 \pi i \R) \; .
\]
\end{Thm}

\begin{Proof}
\cite{BL}, 3.3--3.6, or \cite{W4}, Cor.~4.10.a). 
\end{Proof}

\begin{Rem} \label{2J}
(a) In the above situation, there is a canonical isomorphism between
\[
\Ext^1_{\MHSR} (\R (0) , \Sym^{k-2} \fH (1)) = 
\Sym^{k-2} H^1 (\Eh (\C) , 2 \pi i \R)
\]
and a certain absolute Hodge cohomology group:\\

For an object $\V$ of a $\Q$-linear abelian category, which carries an action 
of $\fS_{n}$, denote by $\V^{\sgn}$ the sign character eigencomponent of 
$\V$.
Consider the action of $\fS_{k-1}$ on $\Eh^{(k-2)}$.
The K\"unneth isomorphism shows that 
\[
H^i (\Eh^{(k-2)} (\BC), 
(2 \pi i)^{k-1} \R)^{\sgn} = 0 
\quad \text{for} \quad i \neq k-2 \; ,
\]
and gives an identification, canonical up to sign, of
\[
H^{k-2} (\Eh^{(k-2)} (\BC) , (2\pi i)^{k-1} \R)^{\sgn}
\] 
and the vector space 
underlying $\Sym^{k-2} \fH (1)$ (for details, see the discussion in
\cite{W5}, 1.2). The Leray spectral sequence shows then that
\[
\Ext^1_{\MHSR} (\R (0) , \Sym^{k-2} \fH (1)) = 
\Ext^{k-1}_{\MHM_{\R} (\Eh^{(k-2)})} (\R (0) , \R (k-1))^{\sgn} \; .
\]
By \cite{HW}, Cor.~A.1.10, the latter group 
equals absolute Hodge cohomology
\[
H^{k-1}_{\abs} (\Eh^{(k-2)} , \R (k-1)_{\Eh^{(k-2)}})^{\sgn} \; ,
\]
and we interpret $s^{\ast} \pol_k = G_{\Eh , k} (s)$, for $s \in \tEh (\C)$, 
as an element of this group. Note that for these indices, the
absolute Hodge cohomology group above equals
Deligne cohomology
\[
H^{k-1}_{\Dh} (\Eh^{(k-2)} , \R (k-1)_{\Eh^{(k-2)}})^{\sgn} 
\]
(\cite{N}, (7.1)). \\[0.2cm]

(b) Because of the existence of the \emph{Eisenstein symbol on torsion} 
(\cite{B5}, Section 3; see also \cite{De1}, Section 8), and since its 
composition with the regulator to 
absolute Hodge cohomology has a description in 
terms of Eisenstein--Kronecker double series (\cite{De1}, Sections 9, 10), the 
element $s^{\ast} \pol_k$ lies in the image of the regulator
\[
r : H^{k-1}_{\Mh} (\Eh^{(k-2)} , \Q (k-1)_{\Eh^{(k-2)}})^{\sgn} \longrightarrow 
H^{k-1}_{\abs} (\Eh^{(k-2)} , \R (k-1)_{\Eh^{(k-2)}})^{\sgn}
\]
as soon as $s$ is a torsion point.\\

By contrast, $s^{\ast} \pol_k$ should not be expected to be of motivic origin 
if $s$ is non-torsion. \\[0.2cm]
(c) Returning to the general case of an elliptic curve
\[
\pi : \Eh \longrightarrow B
\]
and one of the sheaf theories of Subsection \ref{1}, it is quite wrong to expect 
the weight filtration of $s^{\ast} \Log$ to split for any $s \in \tEh (B)$. 
Still, the machinery developed in \cite{W5}, Section 3 allows to construct 
elements in
\[
\Ext^1_{\SsB} (F (0) , \Sym^{k-2} \fH (1)) = 
\Ext^1_{\SsWB} (F (0) , \Sym^{k-2} \fH (1))
\]
from \emph{specific} linear combinations of the $s^{\ast} \pol$.\\

Let us be more precise: in order to apply the results of loc.\ cit., we need 
to check that $\Sh^{s,W}$ satisfies the axioms of loc.\ cit., 3.1. For the 
\'etale setting, and that of admissible variations, this is the content of 
loc.\ cit., 3.2.a) and b). For variations over $\R$, the same proof as in 
loc.\ cit., 3.2.b) works, using the formalism of Grothendieck's functors for 
Hodge modules over $\R$ (\cite{HW}, Thm.~A.2.5), and \cite{W4}, pp.~270--273.\\

If we think of $\pol$ as a framed pro-sheaf satisfying 
\[
\Gr^W_{\bullet} \pol = \Gr_{-1}^W \Log \; |_{\tEh} \oplus 
\Gr^W_{\bullet} \Log (1) \; |_{\tEh} \; ,
\]
then the isomorphism
\[
\kappa : \Gr^W_{\bullet} \Log \silo \prod_{n \ge 0} \Sym^n \fH
\]
induces isomorphisms
\begin{eqnarray*}
\fH & \silo & \Gr^W_{-1} \pol \; , \\
\Gr^W_{-k-1} \pol & \silo & \Sym^{k-1} \fH (1) \; , \; k \ge 1 \; .
\end{eqnarray*}
Fix $s \in \tEh (B)$ and $ k\ge 1$, and define
\begin{eqnarray*}
x : \fH & \silo & s^{\ast} \Gr^W_{-1} \pol \; , \\
y : s^{\ast} \Gr^W_{-k-1} \pol & \silo & \Sym^{k-1} \fH (1)
\end{eqnarray*}
as the respective pullbacks via $s$ of the isomorphisms above. \\

If we denote by $\Lie_B$ the Lie algebra of the pro-unipotent part of the 
Tannakian dual of $\SsWB$, then we may interpret the ``coefficient'' 
$c_{y,x}$ of the representation $\pol$ of $\Lie_B$ as an element of
\[
\Gamma (B , \Lie^{\vee}_B \otimes_F \fH^{\vee} \otimes_F \Sym^{k-1} \fH (1))
\]
(see \cite{W5}, 3.3.).
We define $\{ s\}_k$ as follows:
\[
\{ s \}_1 := c_{y,x} \in \Gamma (B , \Lie^{\vee}_B \otimes_F 
\fH^{\vee} (1)) \; .
\]
For $k \ge 2$, 
\[
\{ s \}_k \in \Gamma (B , \Lie^{\vee}_B \otimes_F \Sym^{k-2} \fH (1))
\]
is the image of $c_{y,x}$ under the epimorphism
\[
\fH^{\vee} \otimes_F \Sym^{k-1} \fH (1) \longonto 
\Sym^{k-2} \fH (1) 
\]
given by $1/k$ times ``derivation'', i.e., $(k-1)/k$ times the symmetrization
of 
\[
\fH^{\vee} \otimes_F \fH \otimes_F \fH^{\otimes (k-2)}(1) \longonto
\fH^{\otimes (k-2)}(1) \; .
\]
The group $\Ext^1_{\SsB} (F (0) , \Sym^{k-2} (\fH (1))$ can be identified 
with the kernel of the differential
\[
d := d \otimes \id : \Gamma (B , \Lie^{\vee}_B \otimes_F \Sym^{k-2} \fH (1)) 
\pfeil \Gamma (B , \bigwedge^2 \Lie^{\vee}_B \otimes_F \Sym^{k-2} \fH (1)) \; ,
\]
where
\[
d : \Lie^{\vee}_B \longrightarrow \bigwedge^2 \Lie^{\vee}_B
\]
is minus the dual map of the commutator (see \cite{W5}, Section 2).\\

On elements of the shape $\{s \}_k$, the differential can be made explicit 
(\cite{W5}, Thm.~3.4), and consequently, we get a statement about which elements in
\[
\langle \{ s \}_k \tei s \in \tEh (B) \rangle_F
\]
lie in $\Ext^1_{\SsB} (F (0) , \Sym^{k-2} \fH (1))$ (\cite{W5}, 
Corollary~3.5).\\

Since we need to reprove \cite{W5}, 3.4 and 3.5 in a geometric 
(rather than sheaf-theoretical) way, we refer to Section~\ref{II}
for the precise statements.\\[0.2cm]
(d) If there is a category $\MMs$ of smooth motivic sheaves ``with 
elliptic polylogs'', i.e., satisfying the axioms of \cite{W5}, 3.1, then the 
same formal arguments as in (c) allow to deduce a statement on the intersection 
of
\[
\langle \{ s \}_k \tei s \in \tEh (B) \rangle_{\Q}
\]
and $\Ext^1_{\MMs\! B} (\Q (0) , \Sym^{k-1} \fH (1))$. If one assumes 
furthermore that the latter group is canonically isomorphic to 
\[
H^{k-1}_{\Mh} (\Eh^{(k-2)} , \Q (k-1)_{\Eh^{(k-2)}})^{\sgn} \; ,
\]
then it is possible to deduce Parts 1 and 2 of the {\em weak version of 
Zagier's conjecture for elliptic curves} (\cite{W5}, Conj.~1.6.B)) 
formally from the existence of a Hodge realization functor, and the result 
recalled in Theorem \ref{2I} (see \cite{W5}, 
Section 4 for the details). The purpose of the present article is to show
part of this conjecture without assuming any of the motivic folklore.
\end{Rem}

\newpage
\subsection{The geometric set-up} \label{3}

For easier reference, we assemble the notation used in the next subsections. Let
\geteilt
{\begin{gather*} A := \C \quad \text{or} \quad A := \R \; , \end{gather*}}
{\begin{gather*} \ell := \; \text{a fixed prime number,} \\
                 A := \Z \left[ \frac{1}{\ell} \right]  , \end{gather*}} 
                 
\noindent $B$ a smooth separated connected $A$-scheme of finite type, 
and of relative 
dimension $d (B)$, 
$\pi : \Eh \to B$ an elliptic curve with zero section $i$.\\

Set $\tEh := \Eh - i (B)$,
\begin{eqnarray*}
j : \tEh & \longinto & \Eh \; , \\
\tpi := \pi \verk j : \tEh & \longrightarrow & B \; ,
\end{eqnarray*}
$\fH := \Hh^0_B (\Eh , F (1))$.\\

By abuse of notation, we also denote the pullback of $\fH$ to $\Eh$
or $\widetilde{\Eh}$ by $\fH$.\\

Fix a subset $P \subset \Eh (B)$.
Consider an open $B$-subscheme $U$ of $\Eh$, 
which is complementary to 
the union of the images of finitely many sections 
$s_i \in P \subset \Eh (B)$. Define a reduced closed subscheme of $\Eh$ by
\[
U_{\infty} := \Eh - U \; .
\]
We thus have
\[
U = \Eh - \bigcup_{s \in U_{\infty} (B)} s (B) \; .
\]
Define $\Ch_{\pi , P}$ to be the set 
of those $U$ for which $U_{\infty} (B)$ 
consists of pairwise disjoint sections in $P$:
\[
s (B) \cap s' (B) = \emptyset \quad \mbox{for} \quad 
s , s' \in U_{\infty} (B) \; , \; s \neq s' \; .
\]
Write $\Ch_{\pi}$ for $\Ch_{\pi , \Eh (B)}$.\\

Fix an involution 
$\iota$ of $\Eh$ satisfying the following hypothesis: $\iota$ acts via 
multiplication by $-1$ on $\fH$.
$\iota$ will automatically act trivially on
\[
F (1) = \Hh^{-1}_B (\Eh , F (1)) \quad \mbox{and} \quad 
F (0) = \Hh^1_B (\Eh , F(1)) \; .
\]
If $P \subset \Eh(B)$ is stable under $\iota$, then $\iota$ acts on
$\Ch_{\pi , P}$. 
The subset of objects stable under $\iota$ will be denoted by
$\Ch_{\pi , P , \iota}$. Write $\Ch_{\pi , \iota}$ for 
$\Ch_{\pi , \Eh (B) , \iota}$.

\begin{Ex} \label{3A}
We will mostly be concerned with the case $\iota = [-1]$, where we set
\[
\Ch_{\pi , P , \iota} =: \Ch_{\pi , P , -} \; .
\]
\end{Ex} 

For $U \in \Ch_{\pi}$ and $n \ge 0$, define
\begin{eqnarray*}
\pi^n & : & \Eh^n \pfeil B \; , \\
j^n_U & : & U^n \longinto \Eh^n \; , \\
i^{(n)}_U & : & U^{(n)}_{\infty} := \Eh^n - U^n \longinto \Eh^n \; ,
\end{eqnarray*}
where $U^{(n)}_{\infty}$ carries the reduced scheme structure. 
(So $\pi^0 = j^0_U = \id_B$, and $U^{(0)}_{\infty} = \emptyset$.)\\

Let $U^{(n)}_{\infty , \reg}$ be the smooth part of $U^{(n)}_{\infty}$, 
and $U^{(n)}_{\infty, \sing}$ its complement. For any subscheme of $\Eh^n$, 
the subscript $\reg$ will denote the complement of 
$U^{(n)}_{\infty , \sing}$ in the subscheme. We shall work with the residue 
sequence associated to the following geometric situation:
\[
\xymatrix@1{ 
\; U^n \; \ar@{^{ (}->}[r]_{j^n_U} & \; \Eh^n_{\reg} \; & 
\; U^{(n)}_{\infty,\reg} \; 
\ar@{_{ (}->}[l]^{i^{(n)}_{U,\reg}} \; .}
\]
Observe that we have
\[
U^{(n)}_{\infty,\reg} = 
U_{\infty} \times_B U^{n-1} \coprod U \times_B U_{\infty} \times_B U^{n-2} 
                           \coprod \ldots \coprod U^{n-1} \times_B U_{\infty}
\]
for $n \ge 1$.\\

For the base change
\[
\pr_2 : \Eh \times_B \tEh \longrightarrow \tEh \; ,
\]
denote by $V \in \Ch_{\pr_2}$ the complement of the zero section $i$, 
and the diagonal $\Delta$. 

\begin{Ex} \label{3B}
On the base change $\Eh \times_B \tEh$, we may consider the involution
\[
\iota: (x,y) \longmapsto (y - x , y) \; .
\]
Note that this involution {\em does} act by 
$-1$ on $\fH = \fH \, |_{\tEh}$: it is the composition of $[-1]$ and the 
translation by $\Delta$. But translations by sections act trivially on $\fH$. 
While multiplication by $[-1]$ does not stabilize 
$V$, this involution does. It interchanges the sections 
$i , \Delta \in V_{\infty} (\tEh)$. Thus we have
\[
V \in \Ch_{\pr_2 , \iota} \; .
\]
\end{Ex}

Finally, we denote by $\vartheta$ the inclusion of the kernel $\Eh^{(n)}$
of the summation map
\[
\Eh^{n} \pfeil \Eh 
\] 
into $\Eh^n$.

\newpage
\subsection{Residue sequences. I} \label{4}

In this subsection, we are going to associate to $U \in \Ch_{\pi , \iota}$ a 
projective system $(\Gh^{(n)}_U)_{n \ge 0}$ of smooth sheaves on $B$. The
$\Gh^{(n)}_U$ will be constructed as direct summands of relative cohomology 
objects with coefficients in Tate twists of certain schemes over $B$ 
(\ref{4A}). 
For the transition map from $\Gh^{(n)}_U$ to $\Gh^{(n-1)}_U$, we use the 
boundary in the long exact sequence associated to the {\it residue at 
$U_{\infty}$} (\ref{4D}).

\begin{Prop} \label{4K}
(a) $\Rh_B (\Eh^n , F(n))^{\sgn}_{(- , \ldots , -)} = 
\Hh^0_B (\Eh^n , F (n))^{\sgn} [0]$.\\ 
The K\"unneth formula gives a 
canonical isomorphism
\[
\Hh^0_B (\Eh^n , F (n))^{\sgn} \silo \Sym^n \Hh^0_B (\Eh , F(1)) = \Sym^n \fH \; .
\]
(b) For $i \neq 0$, we have
\[
\Hh^{i + n+d (B)}_B (\Eh^n , F (n)_{\Eh^n} )^{\sgn}_{(- , \ldots , -)} = 
\Hh^i_B (\Eh^n , F (n))^{\sgn}_{(- , \ldots , -)} = 0 \; .
\]
\end{Prop}

\begin{Proof}
This is an immediate consequence of our hypothesis on the action of 
$\iota$ on the $\Hh^i_B (\Eh , F (1))$.
\end{Proof}

\begin{Def} \label{4A}
For $U \in \Ch_{\pi , \iota}$ and $n \ge 0$,
\[
\Gh^{(n)}_U := \Hh^0_B (U^n , F (n))^{\sgn}_{(- , \ldots , -)} = 
\Hh^{n+d (B)}_B (U^n , F (n)_{U^n})^{\sgn}_{(-, \ldots , -)} \; ,
\]
where the subscript refers to the intersection of the $(-1)$-eigenspaces of the 
componentwise application of $\iota$, and the superscript $\sgn$ refers to the 
sign-eigenspace under the natural action of the symmetric group $\fS_n$ on 
$U^n$.
\end{Def}

Observe in particular that $\Gh^{(0)}_U = F (0)$.\\

The following is an immediate consequence of the K\"unneth formula, and the 
graded-commutativity of the cup product:

\begin{Prop} \label{4B}
There is a canonical isomorphism
\[
\Gh^{(n)}_U \silo \Sym^n \Gh^{(1)}_U \; .
\]
\end{Prop}

For each 
$n \ge 1$, we want to construct a morphism
\[
\Gh^{(n)}_U \longrightarrow \Gh^{(n-1)}_U \otimes_F F [U_{\infty} (B)] (0) 
\]
via the residue. Assume that $U \in \Ch_{\pi}$.

\begin{Prop} \label{4C}
For any complex $M \in D^b \Sh \Eh^n_{\reg}$, 
which is a shift of a smooth sheaf on $\Eh^n_{\reg}$,
there is an exact triangle
\[
\vcenter{\xymatrix@R-10pt{ 
(i^{(n)}_{U,\reg})_{\ast} (i^{(n)}_{U,\reg})^{\ast} 
M (-1) [-2] \ar[rr] &&M \ar[dl] \ar @{} [drr] |{\stern} && \\
& (j^n_U)_{\ast} (j^n_U)^{\ast} M \ar[ul]^{[1]} &&& \\}}    
\]
\end{Prop}

\begin{Proof}
This is purity for the closed immersion
\[
i^{(n)}_{U,\reg} : U^{(n)}_{\infty , \reg} \longinto \Eh^n_{\reg} 
\]
of smooth $B$-schemes. 
\end{Proof}

We apply this to $M = F (j)_{\Eh^n_{\reg}}$, and evaluate the 
cohomological functor
\[
H^\punkt_{\abs} (\Eh^n_{\reg} , \ast)
\]
on the triangle $\stern$:
\begin{eqnarray*}
\ldots \stackrel{\res}{\longrightarrow} 
H^{i-2}_{\abs} (U^{(n)}_{\infty , \reg}, j-1) \!&\! \longrightarrow \!&\! 
H^i_{\abs} (\Eh^n_{\reg} , j) \longrightarrow H^i_{\abs} (U^n , j) \\
\stackrel{\res}{\longrightarrow} 
H^{i-1}_{\abs} (U^{(n)}_{\infty , \reg} , j-1) \!&\! \longrightarrow \!&\! 
\ldots
\end{eqnarray*}
We refer to this as the {\em absolute residue sequence}.\\

Application of the cohomological functor 
$\Hh^\punkt_B (\Eh^n_{\reg} , \ast)$ to the same exact triangle 
yields a long exact sequence of sheaves that we call the 
{\em relative residue sequence}:
\arraycolsep1mm 
\begin{eqnarray*}
\ldots \pfeil \Hh^{i-2}_B (U^{(n)}_{\infty,\reg} , F (j-1)_{U^{(n)}_{\infty,\reg}}) 
& \pfeil & \Hh^i_B (\Eh^n_{\reg} , F (j)_{\Eh^n_{\reg}}) \\
& & \hfill \pfeil \Hh^i_B (U^n , F (j)_{U^n}) \\
\pfeil \Hh^{i-1}_B (U^{(n)}_{\infty,\reg} , F (j-1)_{U^{(n)}_{\infty,\reg}}) 
& \pfeil & \ldots 
\end{eqnarray*}

In order to identify the terms
\[
\begin{array}{l}
H^{i-1}_{\abs} (U^{(n)}_{\infty , \reg} , j-1) \; , \\
\Hh^{i-1}_B (U^{(n)}_{\infty,\reg} , F (j-1)_{U^{(n)}_{\infty,\reg}}) \; , 
\quad n \ge 1 \; ,
\end{array}
\]
recall that we have
\[
U^{(n)}_{\infty,\reg} = U_{\infty} \times_B U^{n-1} \coprod U \times_B 
U_{\infty} \times_B U^{n-2} \coprod \ldots \coprod U^{n-1} \times_B 
U_{\infty} \; .
\]
The proof of the following is formally identical to the proof of 
\cite{HW}, Lemma 4.4.b), c):

\begin{Lem} \label{4D}
(a) We have
$$H^{i-1}_{\abs} (U^{(n)}_{\infty, \reg} , j-1) = 
\bigoplus^n_{k=1} H^{i-1}_{\abs} (U^{n-1} , j-1) \otimes_F 
F [U_{\infty} (B)] \; . $$
The eigenspace $H^{i-1}_{\abs} (U^{(n)}_{\infty , \reg} , j-1)^{\sgn}$ is 
isomorphic to
\[
H^{i-1}_{\abs} (U^{n-1} , j-1)^{\sgn} \otimes_F F [U_{\infty} (B)] \; ,
\]
where the last $\sgn$ refers to the action of $\fS_{n-1}$. 
The isomorphism is given by projection onto the components unequal to $k$, 
for some choice of an element $k \in \{ 1 , \ldots , n \}$. 
It is independent of the choice of $k$.\\[0.2cm]
(b) We have
$$\Rh_B (U^{(n)}_{\infty,\reg} , F (j-1)_{U^{(n)}_{\infty,\reg}}) = 
\bigoplus^n_{k=1} \Rh_B (U^{n-1} , F (j-1)_{U^{n-1}}) \otimes_F 
F [U_{\infty} (B)] \; .$$
As in (a), the eigenpart $\Hh^{i-1}_B (U^{(n)}_{\infty,\reg} , 
F (j-1)_{U^{(n)}_{\infty , \reg}})^{\sgn}$ is canonically isomorphic to
\[
\Hh^{i-1}_B (U^{n-1} , F (j-1)_{U^{n-1}})^{\sgn} \otimes_F 
F [U_{\infty} (B)] \; .
\]
For $i = n + d(B)$ and $j = n$, and $U \in \Ch_{\pi , \iota}$, 
the $(- , \ldots , -)$-eigenpart of the latter is
\[
\Gh^{(n-1)}_U \otimes_F F [U_{\infty} (B)]_{(-)} \; ,
\]
where $F [U_{\infty} (B)]_{(-)}$ denotes the $(-1)$-eigenspace of 
$F [U_{\infty} (B)]$ under the action of $\iota$.
\end{Lem}

Thus, the residue sequences define canonical {\em residue maps}
\begin{eqnarray*}
\res : H^i_{\abs} (U^n , j)^{\sgn} & \pfeil & 
H^{i-1}_{\abs} (U^{n-1} , j-1)^{\sgn} \otimes_F F [U_{\infty} (B)] \; , \\
\res : \Hh^i_B (U^n , F (j)_{U^n})^{\sgn} & \pfeil & 
\Hh^{i-1}_B (U^{n-1} , F (j-1)_{U^{n-1}})^{\sgn} \otimes_F F [U_{\infty} (B)]
\end{eqnarray*}
fitting into the relative and absolute residue sequences. In particular, observe that we have a residue map
\[
\res : \Gh^{(n)}_U \longrightarrow \Gh^{(n-1)}_U \otimes_F F [U_{\infty} (B)]_{(-)} (0)
\]
for $U \in \Ch_{\pi , \iota}$.\\

Note the following consequence of Lemma~\ref{4D}~(b):

\begin{Cor} \label{4Da}
For $i \ne n + d(B)$, the eigenparts 
\[
\Hh^{i-1}_B (U^{(n)}_{\infty,\reg}, 
F (j-1)_{U^{(n)}_{\infty , \reg}})^{\sgn}_{(- , \ldots , -)}
\]
and
\[
\Hh^{i-1}_B (U^{n-1} , F (j-1)_{U^{n-1}})^{\sgn}_{(- , \ldots , -)}
\]
are trivial.
\end{Cor}

\begin{Proof}
Lemma~\ref{4D}~(b) shows that the first part of the statement is implied
by the second part, which in turn follows from the K\"unneth formula,
together with the vanishing of
$\Hh^i_B(U,F(k))_{(-)}$ for $i \ne 0$. This vanishing is a consequence of 
our hypothesis on the action of $\iota$ on the $\Hh^i_B(\Eh,F(1))$.
\end{Proof}

For $n=1$, we preceding discussion can be summarized as follows:

\begin{Prop} \label{4Db}
The relative residue sequence for $n=1$ is canonically isomorphic to the
short exact sequence
\[
0 \pfeil \fH \pfeil \Gh^{(1)}_U \stackrel{\res}{\pfeil}
F [U_\infty (B)] _- (0) \pfeil 0 \; .
\]
\end{Prop} 

In particular, we get:

\begin{Cor} \label{4Dc}
(a) $\Gh^{(n)}_U \in \Sh^{s,W} \!\! B$.\\[0.2cm]
(b) The weights of $\Gh^{(n)}_U$ are contained in 
$\{ -n, -(n-1),\ldots,0 \}$, and there is a canonical isomorphism 
\[
\Sym^n \fH \isoto W_{-n} \Gh^{(n)}_U \; .
\] 
\end{Cor}

\begin{Proof}
This follows from \ref{4B} and \ref{4Db}.
\end{Proof}

We now apply the 
residue spectral sequence \ref{1F} associated to a filtration of $\Eh^n$ 
by open subschemes. For $U \in \Ch_{\pi}$, let
\[
F_k \Eh^n := \{ (x_1 , \ldots , x_n) \in \Eh^n \tei \mbox{at most $n+k$ 
coordinates lie in} \; U_{\infty} \} \; .
\]
We have
\[
\emptyset = F_{-(n+1)} \Eh^n \subset F_{-n} \Eh^n = U^n \subset 
F_{-(n-1)} \Eh^n = \Eh^n_{\reg} \subset \ldots \subset F_0 \Eh^n = \Eh^n \; .
\]
The ``graded pieces'' of this filtration are the reduced schemes
\begin{eqnarray*}
G_k \Eh^n & := & F_k \Eh^n - F_{k-1} \Eh^n \\
& = & \{ (x_1 , \ldots , x_n) \in \Eh^n \tei \mbox{precisely $n+k$ 
coordinates lie in} \; U_{\infty} \} \; .
\end{eqnarray*}
Observe that $G_k \Eh^n$ splits into several disjoint pieces, each isomorphic 
to 
\[
U^{n+k}_{\infty} \times_B U^{-k} \; .
\]
It follows that the $G_k \Eh^n$ are smooth over $B$, and of codimension
$n+k$ in $\Eh^n$. From Theorem \ref{1F}, we 
conclude:

\begin{Prop} \label{4E}
(a) For $H^{\punkt} \in \{ H^{\punkt}_{\abs} , \Hh^{\punkt}_B \}$, there is a 
natural residue spectral sequence
$$^n\stern \quad \quad \quad \quad ^nE^{p,q}_1 \Longrightarrow 
H^{2n+p+q} (\Eh^n , F (j)_{\Eh^n})^{\sgn} \; ,$$
\vspace*{0.3cm}

\noindent $^nE^{p,q}_1 = H^{-p+q} (U^{-p} , F (-n-p+j)_{U^{-p}})^{\sgn} 
\otimes_F \bigwedge^{n+p} F [U_{\infty} (B)], -n \le p \le 0$,\\[0.1cm] 
$^n E^{p,q}_1 = 0$ otherwise.\\[0.2cm]
(b) We have the equality $\res = \, ^n\partial^{-n,q}_1$ of maps
\[
H^{n+q} (U^n , F (j)_{U^n})^{\sgn} \pfeil 
H^{n+q-1} (U^{n-1} , F (j-1)_{U^{n-1}} )^{\sgn} \otimes_F 
                                              F [U_{\infty} (B)] \; .
\]
\end{Prop}

Assume that we are given the action of
a finite group $\FG$ via automorphisms on $\Eh$, 
covering a $\FG$-action on $B$, and respecting the open subscheme $U$. 
The preceding construction is equivariant with respect to the action of
$\FG^n \rtimes \fS_n$.
For any character
\[
\chi : \FG \longrightarrow F^{\ast} \; ,
\]
we have the $\uchi := (\chi , \ldots , \chi)$-eigenpart $^n\stern_{\chi}$ of 
the residue spectral sequence:

\begin{Prop} \label{4F}
(a) For $H^{\punkt} \in \{ H^{\punkt}_{\abs} , \Hh^{\punkt}_B \}$, there is a 
natural residue spectral sequence
$$^n\stern_{\chi} \quad \quad \quad \quad ^nE^{p,q}_{1,\chi} \Longrightarrow 
H^{2n+p+q} (\Eh^n , F (j)_{\Eh^n})^{\sgn}_{\uchi} \; ,$$
\vspace*{0.3cm}

\noindent $^nE^{p,q}_{1,\chi} = 
H^{-p+q} (U^{-p} , F (-n-p+j)_{U^{-p}})^{\sgn}_{\uchi} \otimes_F 
\bigwedge^{n+p} F [U_{\infty} (B)]_{\chi}, -n \le p \le 0$,\\[0.1cm] 
$^n E^{p,q}_{1,\chi} = 0$ otherwise.\\[0.2cm]
(b) We have the equality $\res = \, ^n\partial^{-n,q}_{1,\chi}$ of maps
\[
H^{n+q} (U^n , F (j)_{U^n})^{\sgn}_{\uchi} \pfeil 
H^{n+q-1} (U^{n-1} , F (j-1)_{U^{n-1}} )^{\sgn}_{\uchi} \otimes_F 
                                                F [U_{\infty} (B)]_{\chi} \; .
\]
\end{Prop}

Note that there are versions for motivic cohomology of the above constructions.
In particular:

\begin{Lem} \label{7A}
(a) There is a canonical long exact {\em motivic residue sequence}
\arraycolsep0.8mm
\begin{eqnarray*}
... & \stackrel{\res}{\to} & H^{i-2}_{\Mh} (U^{n-1} , j-1)^{\sgn} 
\otimes_F \! F [U_{\infty} (B)] \to H^i_{\Mh} (\Eh^n_{\reg} , j)^{\sgn} \to 
H^i_{\Mh} (U^n , j)^{\sgn} \\
& \stackrel{\res}{\to} & H^{i-1}_{\Mh} (U^{n-1} , j-1)^{\sgn} 
\otimes_F \!F [U_{\infty} (B)] \to ...
\end{eqnarray*}
\arraycolsep2mm
(b) Under the regulator, the motivic residue sequence maps to the 
absolute residue sequence.
\end{Lem}

\begin{Proof}
This is the localization sequence for motivic cohomology. The residue maps are 
the Gysin morphisms for inclusions of smooth closed subschemes.
\end{Proof}

\begin{Prop} \label{7B}
Proposition~\ref{4F} holds for $H^{\punkt} = H^{\punkt}_{\Mh}$. In other words, 
there is a residue spectral sequence in motivic cohomology.
Under the regulator, it maps to the residue spectral sequence in
absolute cohomology.
\end{Prop}

\begin{Proof}
The localization sequence was the only ingredient 
in the proof of the existence of the residue spectral sequence \ref{1F} 
associated to a filtration of a scheme by open subschemes.
\end{Proof}

\newpage
\subsection{Geometric origin of the logarithmic sheaf,\\ 
and of the polylogarithmic extension} \label{5}

The aim of this subsection is first to identify $\Log$, or rather, its 
restriction to the complement $\tEh$ of the zero section, as the projective 
limit of the smooth sheaves $\Gh^{(n)}_V$ of Subsection~\ref{4} 
(Theorem \ref{5E}). We then have (see \ref{5Ea})
\begin{eqnarray*}
\pol & \in & \Ext^1_{\SHs \tEh} (\fH , \Log (1) \, |_{\tEh} ) \\
& = & \li{n} \Ext^1_{\SHs \tEh} (\fH , \Gh^{(n)}_V (1)) \\
& = & \li{n} H^{n+2}_{\abs} (\Eh \times_B V^n , n+1)^{\sgn}_{- , ( - , 
                                                            \ldots , -)} \; .
\end{eqnarray*}
Our second task is a geometric 
(rather than sheaf theoretical) interpretation of $\pol$ in this projective 
limit (Corollary \ref{5H}), which will translate easily to the context of 
motivic cohomology.\\

So let
\[
\pr_2 : \Eh \times_B \tEh \longrightarrow \tEh \; ,
\]
$V \in \Ob (\Ch_{\pr_2})$ etc.\ as in Example~\ref{3B}.
We use the generator $(\Delta) - (i)$ to identify $F$ and 
$F [V_{\infty} (\tEh) ]_{(-)}$. 

\begin{Def} \label{5A}
For $n \ge 0$, set
\begin{eqnarray*}
\Gh^{(n)} := \Gh^{(n)}_V & = & \Hh^0_{\tEh} (V^n , F(n))^{\sgn}_{(- , 
                                                           \ldots , -)} \\
& = & \Hh^{n+1 + d(B)}_{\tEh} (V^n , F (n)_{V^n})^{\sgn}_{(- , \ldots , -)} 
\in \SHs \tEh \; .
\end{eqnarray*}
\end{Def}

\begin{Thm} \label{5B}
(a) There is a canonical isomorphism
\[
\Gh^{(n)} \silo \Sym^n \Gh^{(1)} \; .
\]
(b) $\Gh^{(n)} \in \Sh^{s,W} \tEh$.\\[0.2cm]
(c) Restriction from $\Eh^n \times_B \tEh$ to $V^n$, together with the
canonical isomorphism of \ref{4K}, induces an identification 
\[
W_{-n} \Gh^{(n)} = \Sym^n \fH
\]
of subobjects of $\Gh^{(n)}$.
There is an exact sequence
\[
0 \pfeil W_{-n} \Gh^{(n)} = \Sym^n \fH \pfeil \Gh^{(n)} 
\stackrel{\res}{\pfeil} \Gh^{(n-1)} \pfeil 0 \; .
\]
The surjection $\res : \Gh^{(n)} \longonto \Gh^{(n-1)}$ is given by the 
composition of
\[
\res : \Gh^{(n)} \longonto \Gh^{(n-1)} \otimes_F F [V_{\infty} (\tEh) ]_{(-)}
\]
and the identification
\begin{eqnarray*}
\Gh^{(n-1)} & \isoto & \Gh^{(n-1)} \otimes_F F [V_{\infty} (\tEh) ]_{(-)} \\
x & \longmapsto & x \otimes ((\Delta) - (i)) \; .
\end{eqnarray*}
(d) The exact sequences in (c) for $n \ge 0$ induce natural isomorphisms
\[
\Gr^W_{\punkt} \Gh^{(n)} \silo \bigoplus^n_{i=0} \Sym^i \fH \; ,
\]
which fit into commutative diagrams
\[
\vcenter{\xymatrix@R-10pt{ 
\Gr^W_{\punkt} \Gh^{(n)} \ar[r]^-{\cong}  \ar[d]_{\Gr^W_{\punkt} \res} & 
\bigoplus^n_{i=0} \Sym^i \fH \ar[d]^{\can} \\
\Gr^W_{\punkt} \Gh^{(n-1)} \ar[r]^-{\cong} &
\bigoplus^{n-1}_{i=0} \Sym^i \fH  \\}}    
\]
In particular, $\Gh^{(n)}$ is a unipotent sheaf relative to
\[
\tEh \longrightarrow B \; .
\]
(e) The diagrams
\[
\vcenter{\xymatrix@R-10pt{ 
\Gh^{(n)} \ar[r] \ar[d]^{\cong}_{(a)} & 
F (0) \\
\Sym^n \Gh^{(1)} \ar[r] &
\Sym^n F (0) \ar[u]^{\cong} \\}}    
\]
and
\[
\vcenter{\xymatrix@R-10pt{ 
\Sym^n \fH \ar[r] \ar@{=}[d] &
\Gh^{(n)} \ar[d]_{\cong}^{(a)} \\
\Sym^n \fH \ar[r] &
\Sym^n \Gh^{(1)} \\}}    
\]
commute.
\end{Thm}

\begin{Proof}
(a) and (b) are \ref{4B}, resp.\ \ref{4Dc} in our situation. 
Because of the vanishing results \ref{4K}
and \ref{4Da}, and since $F [V_{\infty} (\tEh) ]_{(-)}$ is one-dimensional,
the relative residue spectral sequence of \ref{4F}
degenerates into
a short exact sequence:
\[
0 \pfeil \Sym^n \fH \pfeil \Gh^{(n)} 
\stackrel{\res}{\pfeil} \Gh^{(n-1)} \pfeil 0 \; .
\]
This proves the remaining claims. 
\end{Proof}

Recall the isomorphism
\[
\kappa : \Gr^W_{\punkt} \Log \silo \prod_{i \ge 0} \Sym^i \fH \]
of Subsection~\ref{2}. Defining
\[
\Log^{(n)} := \Log / W_{-n-1} \Log
\]
for $n \ge 0$, we have
\[
\kappa : \Gr^W_{\punkt} \Log^{(n)} \silo \bigoplus^n_{i=0} \Sym^i \fH \; .
\]

\begin{Lem} \label{5C}
There is a unique isomorphism
\[
\Log^{(1)} \, |_{\tEh} \silo \Gh^{(1)}
\]
inducing the identity on the level of $\Gr^W_0 = F (0)$. On 
$\Gr^W_{-1} = \fH$, it induces the identity as well.
\end{Lem}

\begin{Proof}
Let $[\argdot]$ be the Abel--Jacobi map, associating to sections of
\[
\pr_2 : \Eh \times \Eh \longrightarrow \Eh
\]
extensions in
\[
\Ext^1_{\SHs \Eh} (F (0) , \fH)
\]
(see e.g.\ \cite{W4}, p.~272). By loc.\ cit., Prop.~2.5, we have
\[
\Log^{(1)} = [\Delta]
\]
in the framing given by $\kappa$. Since $[\argdot]$ is a homomorphism, it is 
also true that
\[
\Log^{(1)} = [\Delta] - [i] \; .
\]
From the definition of the isomorphism
\[
\Gr^W_0 \Gh^{(1)} \longrightarrow F (0) \; ,
\]
it follows that there is a (unique) isomorphism
\[
\Log^{(1)} \, |_{\tEh} \silo \Gh^{(1)} 
\]
inducing the identity on $\Gr^W_{\punkt}$.\\

Therefore, $\Gh^{(1)}$ extends to an object of $\Sh^{s,W} \Eh$, and admits a 
section over $i (B)$. But then, the universal property \ref{2B} ensures that 
the isomorphism is uniquely determined by its effect on $\Gr^W_0$.
\end{Proof}

\begin{Cor} \label{5D}
(a) $\Gh^{(n)}$ extends uniquely to a smooth sheaf $\uGh^{(n)}$ on $\Eh$. It 
has a weight filtration.\\[0.2cm]
(b) The properties of \ref{5B}~(a)--(e) carry over to $\uGh^{(n)}$. In 
particular, there is a natural isomorphism
\[
\eta^{(n)} : \Gr^W_{\punkt} \uGh^{(n)} \silo \bigoplus^n_{i=0} \Sym^i \fH \; .
\]
(c) The weight filtration of $i^{\ast} \uGh^{(n)}$ is split. In particular, 
there is a canonical monomorphism
\[
1 : F (0) \longrightarrow i^{\ast} \uGh^{(n)} \; .
\]
\end{Cor}

By the universal property of $\Log$ (\ref{2B}), there is a unique morphism
\[
\varphi : \Log \longrightarrow \uGh := \lim \uGh^{(n)}
\]
such that $\varphi \, |_{i (B)}$ sends $1$ to $1$.

\begin{Thm} \label{5E}
$\varphi$ is an isomorphism.
\end{Thm}

\begin{Proof}
By \cite{W4}, Thm.~2.6, we have a canonical isomorphism
\[
\varphi^{(n)}_0 : \Log^{(n)} \silo \Sym^n \Log^{(1)}
\]
respecting the sections 1 over $i (B)$. By Corollary~\ref{2B}, the morphism 
\[
\varphi^{(n)} : \Log \pfeil \underline{\Gh} \pfeil \underline{\Gh}^{(n)}
\]
equals
\[
\xymatrix@1{ 
\; \Log  \; \ar[r] & \; \Log^{(n)} \; \ar[r]^-{\cong}_-{\varphi_0^{(n)}} &
\; \Sym^n \Log^{(1)} \; \ar[r]^-{\cong}_-{\rm \ref{5C}} &
\; \Sym^n \underline{\Gh}^{(1)} \; }
\]
under our identifications of $\underline{\Gh}^{(n)}$ and 
$\Sym^n \underline{\Gh}^{(1)}$.
\end{Proof}

From now on, we use $\varphi$ to identify $\Log$ and $\underline{\Gh}$. 

\begin{Cor} \label{5Ea}
There is a canonical isomorphism
\begin{eqnarray*}
& & \Ext^1_{\SHs \tEh} (\fH , \Log (1) \, |_{\tEh} ) \\
& = & \li{n} \Ext^1_{\SHs \tEh} (\fH , \Gh^{(n)}_V (1)) \\
& = & \li{n} H^{n+2}_{\abs} (\Eh \times_B V^n , n+1)^{\sgn}_{- , ( - , 
                                                            \ldots , -)} \; .
\end{eqnarray*}
\end{Cor}

\begin{Proof}
This follows from \ref{5E}, \ref{4Da}, and \ref{4K}.
\end{Proof}

In 
order to compare the extension classes to be constructed in the Section~\ref{II}
to those described in Theorem \ref{2I}, we need to compare the isomorphism
\[
\eta := \li{n} \eta^{(n)} : \Gr^W_{\punkt} \underline{\Gh} \silo 
\prod_{n \ge 0} \Sym^n \fH
\]
of \ref{5D}~(b) to the isomorphism
\[
\kappa : \Gr^W_{\punkt} \underline{\Gh} = \Gr^W_{\punkt} \Log \silo 
\prod_{n \ge 0} \Sym^n \fH
\]
of Subsection~\ref{2}. Denote by $\eta_{-n}$ and $\kappa_{-n}$ the isomorphisms
\[
\Gr^W_{-n} \underline{\Gh} \silo \Sym^n \fH \; .
\]

\begin{Prop} \label{5F}
We have the equality
\[
\eta_{-n} = n! \cdot \kappa_{-n} \; .
\]
\end{Prop}

\begin{Proof}
For $n = 0,1$, this follows from Lemma \ref{5C}. For $n \ge 2$, let
\[
\varphi^{(n)}_0 : \uGh^{(n)} \silo \Sym^n \uGh^{(1)}
\]
be the isomorphism of \ref{5B}~(a) and \ref{5D}~(b).\\

By \ref{5B}~(e) and \ref{5D}~(b), the diagram
\[
\vcenter{\xymatrix@R-10pt{ 
\uGh^{(n)} \ar[r] \ar[d]^{\cong}_{\varphi_0^{(n)}} & 
F (0) \\
\Sym^n \uGh^{(1)} \ar[r] &
\Sym^n F (0) \ar[u]^{\cong} \\}}    
\]
commutes. By \cite{W4}, Thm.~2.6.a), the commutativity of this diagram 
characterizes $\varphi^{(n)}_0$ uniquely. From loc.\ cit., Thm.~2.6.b) and 
c), we know that the diagram
\[
\vcenter{\xymatrix@R-10pt{ 
\Sym^n \fH \ar[r]^{\frac{1}{n!} \cdot \kappa^{-1}} \ar@{=}[d] &
\uGh^{(n)} \ar[d]_{\cong}^{\varphi_0^{(n)}} \\
\Sym^n \fH \ar[r]^{\Sym^n \kappa^{-1}} &
\Sym^n \uGh^{(1)} \\}}    
\]
commutes. So our identity
\[
\eta_{-n} = n! \cdot \kappa_{-n}
\]
follows from \ref{5B}~(e) and \ref{5D}~(b).
\end{Proof}

The rest of this subsection will deal with a reinterpretation of the 
construction of $\pol$ recalled in Subsection~\ref{2} in terms of the absolute 
cohomology of the schemes
\[
\Eh \times_B V^n \; , \quad n \ge 0 \; .
\]
Recall the Abel--Jacobi morphism $[\argdot]$ associating elements of
\[
\Ext^1_{\SHs \Eh} (F (0) , \fH) = \Ext^1_{\SHs \Eh} (\fH , F (1))
\]
to sections of $\pr_2$. In particular, we have the 
elements
\[
[\Delta] \; , \; [- \Delta] = - [\Delta] \; , \; [i] = 0 \in 
\Ext^1_{\SHs \Eh} (\fH , F (1)) \; .
\]
Denote by $\pol^{(n)}$ the image of $\pol$ in
\[
\Ext^1_{\SHs \tEh} (\fH , \Gh^{(n-1)} (1)) \; .
\]
We have the following

\begin{Thm} \label{5G}
The map
\[
\li{n} \Ext^1_{\SHs \tEh} (\fH , \Gh^{(n)} (1)) \longrightarrow 
\Ext^1_{\SHs \tEh} (\fH , F (1))
\]
is injective. The image $\pol^{(1)}$ of $\pol$ is given by the image of
\[
[\Delta] = [\Delta] - [i] = \halb \left( [\Delta] - [- \Delta] \right)
\in \Ext^1_{\SHs \Eh} (\fH , F (1))
\]
under restriction from $\Eh$ to $\tEh$.
\end{Thm}

\begin{Proof}
The injectivity statement follows from Theorem \ref{2G}: The projective 
limit equals
\[
\End_{\Sh B} (\fH) = 
\Hom_{\Sh B} \left( \fH , \prod_{n \ge 1} \Sym^n \fH \right) \; .
\]
Our map is injective since it admits a left inverse given by the boundary in 
the Leray spectral sequence:
\[
\Ext^1_{\Sh \tEh} (\fH , F (1)) \longrightarrow 
\Hom_{\Sh B} (\fH , \fH) \; .
\]
Then our claim follows from \cite{W4}, Prop.~2.4.
\end{Proof}

The desired geometrical interpretation of the elliptic polylogarithm looks as 
follows:

\begin{Cor} \label{5H}
Under the identification 
of \ref{5Ea}, $\pol$ is the unique element of
\[
\li{n} H^{n+1}_{\abs} (\Eh \times_B V^{n-1} , n)^{\sgn}_{-,(-, \ldots , -)} \; ,
\]
whose $H^2_{\abs} (\Eh \times \tEh , 1)_-$-component $\pol^{(1)}$ is the 
image of
\[
[\Delta] - [i] = \halb ([\Delta] - [- \Delta])
\in H^2_{\abs} (\Eh \times_B \Eh , 1)_-
\]
under restriction from $\Eh$ to $\tEh$.
\end{Cor}

\newpage
\subsection{The polylog in motivic cohomology} \label{7}

Our aim here is to imitate the constructions of the Subsection~\ref{5}, and 
transfer them to the context of motivic cohomology. In particular, we are going 
to analyse the residue maps
\[
H^{n+2}_{\Mh} (\Eh \times_B V^n , n+1)^{\sgn}
\longrightarrow 
H^{n+1}_{\Mh} (\Eh \times_B V^{n-1} , n)^{\sgn} \; ,
\]
and to define a projective system
\[
\pol = \left( \pol^{(n)} \right)_n \in \li{n} 
H^{n+1}_{\Mh} (\Eh \times_B V^{n-1} , n)^{\sgn}_{- , (- , \ldots , -)}
\]
satisfying
\begin{eqnarray*}
\pol^{(1)} = [\Delta] - [i] = \halb ([\Delta] - [- \Delta]) & \in & 
H^2_{\Mh} (\Eh \times_B \tEh , 1)_- \\
& = & CH^1 (\Eh \times_B \tEh)_- \otimes_{\Z} \Q
\end{eqnarray*}
(Definition \ref{7G}). Compatibility with the element $\pol$ in absolute 
cohomology under the regulators (Theorem \ref{7H}) is then an immediate 
consequence of Corollary \ref{5H}.\\

Concerning notation, let us agree that as long as our statements concern only 
motivic cohomology, the restrictions on the base $B$ set up in Subsection~\ref{3} 
do not apply: $B$ is only assumed to be regular, noetherian and connected. 
If a statement concerns regulators, then it is assumed that $B$ is such that 
the absolute cohomology is defined, i.e., as in Subsection~\ref{3}.\\

Now let us recapitulate the main technical tools used in the construction of 
$\pol$ in absolute cohomology:\\

Firstly, we considered the complementary inclusions
\[
U^n \longinto \Eh^n_{\reg} \longleftinto U^{(n)}_{\infty ,\reg}
\]
in order to set up our absolute residue sequence, and in particular, the 
residue maps for absolute cohomology. The identification of the terms belonging 
to $U^{(n)}_{\infty , \reg}$ with absolute cohomology of $U^{n-1}$ was a 
combinatorial result (Lemma \ref{4D}~(a)). \\

Secondly, an analysis of the sheaf theoretic situation allowed to deduce the 
injectivity statement \ref{5G}.\\

While there is a motivic residue sequence (see
\ref{7A}), observe that
for motivic cohomology, we cannot use sheaf theory to 
deduce results of the type \ref{5G}.\\

The additional ingredient needed to handle the situation is the study of the 
action of the isogenies $[a]^n$ on motivic cohomology.

\begin{Rem} \label{7C}
(a) We translate the treatment in Sections 6.2, 6.3 of \cite{BL} to
our geometric context. \\[0.2cm]
(b) An analysis of the proof of \cite{BL}, Lemma~6.2.1 
shows that the construction below is independent of the 
choice of $a$ as long as the absolute value of $a$ is at least $2$.\\
\end{Rem}

For an open subscheme $W \subset \Eh^n$ satisfying
\[
([a]^n)^{-1} W \subset W \; ,
\]
define the endomorphism $\tr_{[a]^n}$ on $H^{\punkt}_{\Mh} (W , \ast)$ as the 
composition of the restriction to $([a]^n)^{-1} W$, and the norm
\[
[a]^n_{\ast} : H^{\punkt}_{\Mh} (([a]^n)^{-1} W , \ast) \longrightarrow 
H^{\punkt}_{\Mh} (W , \ast) \; .
\]
Now consider the motivic residue map \ref{7A} for 
$V \in \Ob (\Ch_{\pr_2})$:
\[
\res: H^i_{\Mh} (V^n , j)^{\sgn} \longto 
      H^{i-1}_{\Mh} (V^{n-1} , j-1)^{\sgn} \; .
\]
The operators $\tr_{[a]^{n+1}}$ resp.\ $\tr_{[a]^{n}}$ act on the source resp.\
the target of $\res$. We have:

\begin{Prop} \label{7D}
$\res$ is compatible with $\tr_{[a]^{n+1}}$ and $\tr_{[a]^{n}}$.
\end{Prop}

Now consider the generalized eigenspaces $(\argdot)_{(r)}$ of the eigenvalue 
$a^r$ for $\tr_{[a]^{n+1}}$. We get residue maps
\[
\res: H^i_{\Mh} (V^n , j)^{\sgn}_{(r)} \longto 
      H^{i-1}_{\Mh} (V^{n-1} , j-1)^{\sgn}_{(r)} \; ,
\]
which we need to understand in the following two cases:

\begin{Thm} \label{7E}
(a) For any $n \ge 1$, the residue
\[
\res 
: H^{\bullet}_{\Mh} (V^n , \argstern)^{\sgn}_{(1)} \longrightarrow 
H^{\bullet - 1}_{\Mh} (V^{n-1} , \argstern -1)^{\sgn}_{(1)}
\]
is an isomorphism.\\[0.2cm]
(b) For any $n \ge 2$, the residue
\[
\res 
: H^{\bullet}_{\Mh} (V^n , \argstern)^{\sgn}_{(2)} \longrightarrow 
H^{\bullet - 1}_{\Mh} (V^{n-1} , \argstern - 1)^{\sgn}_{(2)}
\]
is an isomorphism.
\end{Thm}

\begin{Proof}
First, we remark that the results of \cite{Ki}, 2.1--2.2 hold in the context of
elliptic curves over regular noetherian base schemes, rather than just over 
base schemes which are smooth and quasi-projective over a field: indeed,
the only reason for this restriction is the need to apply the main result 
of \cite{DeM} in the proof of \cite{Ki}, Thm.~2.2.3; 
for elliptic curves, the central vanishing result in the proof of  loc.~cit.\:
\[
H^\bullet_{\Mh} (\tEh, \argstern)_{(r)} = 0 \quad \text{if} \quad r \le 0 \; ,
\]
is the case $n=0$ of \cite{BL}, Corollary~6.2.4.\\

Our claim follows from \cite{Ki}, Cor.~2.1.4 and Thm.~2.2.3.
\end{Proof}

Again, we perform base change to
\[
\pr_1 : \Eh \times_B \Eh \longrightarrow \Eh \; .
\]
\begin{Cor} [cmp.~\cite{BL}, 6.3.4] \label{7F}
For any $n \ge 2$, the residue 
\[
H^{n+1}_{\Mh} (\Eh \times_B V^{n-1} ,n)^{\sgn}_{(1)} \longrightarrow 
H^n_{\Mh} (\Eh \times_B V^{n-2} ,n-1)^{\sgn}_{(1)}
\]
is an isomorphism.
\end{Cor}

We are now in a position to define the motivic
version of the elliptic polylog:

\begin{Def} [cmp.~\cite{BL}, 6.3.5] \label{7G}
Define the {\em elliptic motivic polylogarithmic classes}
\[
\pol^{(n)} \in 
H^{n+1}_{\Mh} (\Eh \times_B V^{n-1} , n)^{\sgn}_{(1)} 
\; , \quad n \ge 1 \; ,
\]
as the preimages of
\begin{eqnarray*}
[\Delta] - [i] = \halb ([\Delta] - [- \Delta]) & \in & 
H^2_{\Mh} (\Eh \times_B \tEh , 1)_{(1)} \\
& = & CH^1 (\Eh \times_B \tEh)_{-} \otimes_{\Z} \Q
\end{eqnarray*}
under the composition of the residue isomorphisms. Here, the subscripts
$(1)$ and $-$ concern the factor $\tEh$ only, and the last equality is 
a consequence of \cite{BL}, Cor.~6.2.4.
\end{Def}

From Theorem~\ref{7E}, and the behaviour of $[\Delta] - [i]$ under $\iota$,
we conclude formally:

\begin{Prop} \label{7Ga}
The class $\pol^{(n)}$ lies in fact in
\[
H^{n+1}_{\Mh} (\Eh \times_B V^{n-1} , n)^{\sgn}_{-,(1),(-,...,-)} \; ,
\]
where the first $-$ refers to the $(-1)$-eigenspace for the action of 
$[-1]$ on the first component $\Eh$.
\end{Prop}

\begin{Rem} \label{7Gb} 
In order to construct $\pol^{(n)}$, it would have been sufficient to
prove the statement of Theorem~\ref{7E} for the $(-,...,-)$-eigenspaces
of the objects considered there. This in turn could have been achieved by
employing the motivic residue spectral sequence \ref{7B}, together with
the following statement: on the $(-,...,-)$-eigenspace of
\[
H^{\bullet}_{\Mh} 
         (\Eh^n \times_B \tEh , \argstern)^{\sgn} \; ,
\]
the norm $\tr_{[a]^{n+1}}$ acts with eigenvalues $a^{k}$, where $k \ge n+1$.
This is proved in the same way as the results contained in \cite{BL}, 6.2;
for more details, we refer to Subsection~\ref{6}. 
Conversely, the motivic spectral sequence \ref{7B}, together with
(the proof of) \ref{7E} implies that the above statement actually holds
{\it before} passing to the $(-,...,-)$-eigenspace:
\[
H^{\bullet}_{\Mh} 
         (\Eh^n \times_B \tEh , \argstern)^{\sgn}_{(r)} = 0
         \quad \mbox{for} \quad r \le n \; .
\]
\end{Rem}

From \ref{5H}, we conclude:

\begin{Thm} \label{7H}
Under the regulators,
\[
\pol := (\pol^{(n)})_n \in \li{n} 
H^{n+1}_{\Mh} (\Eh \times_B V^{n-1} , n)^{\sgn}_{- , (1) , (- , \ldots , -)}
\]
is mapped to
\[
\pol \in \li{n} 
H^{n+1}_{\abs} (\Eh \times_B V^{n-1} , n)^{\sgn}_{- , (1) , (- , \ldots , -)} \; .
\]
\end{Thm}

Recall that we identified $F$ and $F [V_{\infty} (\tEh)]_{(-)}$ by fixing the 
generator
\[
(\Delta) - (i) \in F [V_{\infty} (\tEh)]_{(-)} \; . 
\]
For future reference, we note:

\begin{Prop} \label{7I}
The residue map from
\[
H^{n+1}_{\Mh} (\Eh \times_B V^{n-1} ,n)^{\sgn}_{-, (-, \ldots , -)}
\]
to
\[
H^n_{\Mh} (\Eh \times_B V^{n-2} ,n-1)^{\sgn}_{- , (- , \ldots , -)} 
                              \!\otimes_F F [V_{\infty} (\tEh)]_{(-)}
\]
maps $\pol^{(n)}$ to $\pol^{(n-1)} \otimes ((\Delta) - (i))$ .
\end{Prop}

\begin{Prop} [cmp.~\cite{BL}, 6.3.6] \label{7J} 
(a) The formation of the elliptic motivic polylog is compatible with change of 
the base $B$.\\[0.2cm]
(b) (Norm compatibility.) If $\psi : \Eh_1 \to \Eh_2$ is an \'etale isogeny of 
elliptic curves over $B$, then the norm map
\[
\tr_{\psi^{n+1}} : H^{n+1}_{\Mh} (\Eh_1 \times_B V^{n-1}_1 ,n) \longrightarrow 
H^{n+1}_{\Mh} (\Eh_2 \times_B V^{n-1}_2 ,n)
\]
(defined since $(\psi^n)^{-1} (V^{n-1}_2) \subset V^{n-1}_1)$ maps
\[
\pol^{(n)}_1 \quad \mbox{to} \quad \deg (\psi) \cdot \pol^{(n)}_2 \; .
\]
\end{Prop}

\begin{Proof}
(a) follows from our construction.
For (b), we use the formula
\[
\res \verk \tr_{\psi^{n+1}} = \tr_{\psi^n} \verk \res
\]
to reduce to the case $n = 1$. There,
\[
\psi^2_{\ast} : CH^1 (\Eh_1 \times_B \Eh_1) \otimes_{\Z} \Q \longrightarrow 
CH^1 (\Eh_2 \times_B \Eh_2) \otimes_{\Z} \Q
\]
maps $[\Delta_1]$ to $\deg (\psi) \cdot [\Delta_2]$.
\end{Proof}

We conclude by studying the behaviour of
\[
\pol^{(n)} \in 
H^{n+1}_{\Mh} (\Eh \times_B V^{n-1} , n)^{\sgn}_{(- , \ldots , -)}
\]
under the map
\begin{eqnarray*}
\pr_{\Sgn} : H^{n+1}_{\Mh} (\Eh \times_B V^{n-1} , n)^{\sgn}_{(- , \ldots , -)} 
& \pfeil & H^{n+1}_{\Mh} (V \times_{\tEh} V^{n-1} ,n) \\ 
& \longonto & H^{n+1}_{\Mh} ( V^{n} ,n)^{\Sgn}_{(- , \ldots , -)} \; .
\end{eqnarray*}
$\Sgn$ denotes the sign-eigenpart under the action of $\fS_n$. The first arrow 
is the natural restriction map, and the second arrow is the projection onto 
the $\Sgn$-$(- , \ldots , -)$-eigenspace.

\begin{Prop} \label{7K}
We have
\[
\pr_{\Sgn} (\pol^{(n)}) = 0 \; .
\]
\end{Prop}

\begin{Proof}
Recall that $\pol^{(n)}$ lies in
\[
H^{n+1}_{\Mh} (\Eh \times_B V^{n-1} , n)^{\sgn}_{- , (1), (- , \ldots , -)} \; ,
\]
where the subscript $(1)$ refers to the action of $[a]^n$ on $V^{n-1}$ only. 
By \cite{BL}, Cor.~6.2.2~(ii), 
\[
H^{\punkt}_{\Mh} (\Eh , \ast)_- \subset H^{\punkt}_{\Mh} (\Eh , \ast)_{(1)}
\]
for any elliptic curve. We conclude that we have
\[
\pr_{\Sgn} (\pol^{(n)}) \in 
H^{n+1}_{\Mh} (V^n , n)^{\sgn}_{(2) , (- , \ldots , -)} \; .
\]
Because the residue map and $\pr_{\Sgn}$ commute up to a factor $n$, and 
because of Theorem \ref{7E}~(b), we may assume $n = 1$. By definition,
\[
\pol^{(1)} = [\Delta] - [i] \in H^2_{\Mh} (\Eh \times_B \tEh , 1) \; .
\]
Its restriction to $V$ is therefore visibly trivial. 
\end{Proof}

\begin{OP} \label{7L}
Relate \ref{7G} to Beilinson's and Levin's definition of the elliptic
polylog (\cite{BL}, 6.3.5).
\end{OP}

\newpage
\subsection{The torsion case. I} \label{8}

At this point, we have assembled enough material to analyse the values
of the polylogarithm at torsion sections, i.e., the pull-backs of $\pol$
via sections $s \in \tEh (B)$ which are torsion in $\Eh (B)$. Since we
imitate faithfully the strategy of \cite{BL}~6.4, and since we will
recover the material of the present subsection as a special case of the
formalism to be developed in Section~\ref{II}, we content ourselves with
the \emph{statements} of the results.\\

Let $s \in \tEh (B)$ be torsion in $\Eh(B)$, and define
\[
U_s := \tEh -  s(B) \; .
\]
There is a natural cartesian diagram
\[
\vcenter{\xymatrix@R-10pt{ 
U_s \ar[r]^s \ar[d] & V \ar[d]^{\pr_2}\\
B \ar[r]^s & \tEh  \\}}    
\]
$U_s$ is an element of $\Ch_{\pi, (-)}$, where the subscript
refers to the $(-1)$-eigenspace of the action of
\[
\iota_s: \Eh \pfeil \Eh \; , \quad x \longto s-x \; .
\]
The residue
spectral sequence \ref{4F} associated to $U_s$ collapses into a long
exact sequence because of $\dim_F F [U_{s,\infty} (B) ]_{(-)} = 1$. 
If $[N]s = 0$, then $s$ is fixed by $[N+1]$. An analysis of the action of
$[N+1]$ on the residue sequence yields:

\begin{Prop} [cmp.~\cite{BL}, 6.4.1] \label{8A}
(a) There is an isomorphism,\\
canonical up to the choice of generator of
$F [U_{s,\infty} (B) ]_{(-)}$,
\[
H^{n+1}_{?} (\Eh \times_B U_s^{n-1} ,n)^{\sgn}_{-, (-, \ldots , -)} \isoto
\bigoplus_{k=1}^{n} 
H^{k+1}_{?} (\Eh \times_B \Eh^{k-1} ,k)^{\sgn}_{-, (-, \ldots , -)} \; .
\]
\noindent (b) The isomorphism of (a) is compatible with the residue: there
is a commutative diagram
\[
\vcenter{\xymatrix@R-10pt{ 
H^{n+1}_{?} (\Eh \times_B U_s^{n-1} ,n)^{\sgn}_{-, (-, \ldots , -)} 
\ar[r]^-{(a)}_-{\cong} \ar[d]_{\res} &
\bigoplus_{k=1}^{n} 
H^{k+1}_{?} (\Eh \times_B \Eh^{k-1} ,k)^{\sgn}_{-, (-, \ldots , -)} 
\ar@{>>}[d] \\
H^{n}_{?} (\Eh \times_B U_s^{n-2} ,n-1)^{\sgn}_{-, (-, \ldots , -)} 
\ar[r]^-{(a)}_-{\cong} &
\bigoplus_{k=1}^{n-1} 
H^{k+1}_{?} (\Eh \times_B \Eh^{k-1} ,k)^{\sgn}_{-, (-, \ldots , -)}  \\}}    
\]
\end{Prop}

\begin{Rem} \label{8B}
In fact, one can show that
\[
H^{k+1}_{?} (\Eh \times_B \Eh^{k-1} ,k)^{\sgn}_{-, (-, \ldots , -)}
= H^{k+1}_{?} (\Eh \times_B \Eh^{k-1} ,k)^{\sgn}_{-, -, \ldots , -}
\]
as subspaces of 
$H^{k+1}_{?} (\Eh \times_B \Eh^{k-1} ,k)^{\sgn}$. This is a direct consequence
of the invariance of $H^\bullet_? (\Eh, \argstern)_-$ 
under translations by sections
(which will be proved in \ref{14Ba}).
\end{Rem}

Let us \emph{fix} the generator $((s) - (i))$ of 
$F [U_{s,\infty} (B) ]_{(-)}$, and call an individual
projection
\[
H^{n+1}_{?} (\Eh \times_B U_s^{n-1} ,n)^{\sgn}_{-, (-, \ldots , -)} \longonto
H^{k+1}_{?} (\Eh \times_B \Eh^{k-1} ,k)^{\sgn}_{-, -, \ldots , -} 
\]
the \emph{$k$-component}.
For $k \ge 2$, let
\[
\pr: H^{k+1}_{?} (\Eh \times_B \Eh^{k-1} ,k)^{\sgn}_{-, -, \ldots , -} \pfeil
H^{k-1}_{?} (\Eh^{(k-2)} , k-1)^{\sgn}
\]
be the map $\vartheta^\ast  q_\ast$, where
\[
q: \Eh \times_B \Eh^{k-1} \pfeil \Eh^{k-1} 
\]
is the projection $(x,x_1,\ldots,x_{k-1}) \mapsto (x+x_1,\ldots,x+x_{k-1})$. 

\begin{Def} \label{8C}
(a) The class 
\[
\{ s \}_k' \in 
H^{k+1}_{?} (\Eh \times_B \Eh^{k-1} ,k)^{\sgn}_{-, -, \ldots , -} 
\]
is defined as the $k$-component of the pullback
\[
s^\ast \pol^{(n)} \in
H^{n+1}_{?} (\Eh \times_B U_s^{n-1} ,n)^{\sgn}_{-, (-, \ldots , -)} \; ,
\]
for any $n \ge k$. \\[0.2cm]
(b) Let $k \ge 2$. The class 
\[
\Eis_{?}^{k-2} (s) \in H^{k-1}_{?} (\Eh^{(k-2)} , k-1)^{\sgn}
\]
is defined to be the image under $\pr$ of $\{ s \}_k'$. 
\end{Def}

\begin{Rem} \label{8D}
Because of our normalization of the $k$-component,
and because of \ref{7I}, $\{ s \}_k'$ does not depend on the choice 
of $n \ge k$.
\end{Rem}

\begin{Def} [cmp.~\cite{BL}, 6.4.3] \label{8E}
Let $k \ge 2$. The map
\[
\Eis_{?}^{k-2}: \tEh(B) \cap \Eh(B)_{\tor} \pfeil
H^{k-1}_{?} (\Eh^{(k-2)} , k-1)^{\sgn}
\]
is called the \emph{Eisenstein symbol on torsion}.
\end{Def}

\begin{OP} \label{8F}
Relate this to Beilinson's
original (i.e., non-po\-ly\-lo\-ga\-rith\-mic)
definition of the \emph{Eisenstein symbol on torsion} (\cite{B5}, 3.1). 
\end{OP}

In the setting of absolute Hodge cohomology, 
the relation can be established. More
precisely,
Theorem~\ref{2I} allows to identify the image under
the regulator $r$ to absolute Hodge cohomology
of the Eisenstein classes when $A = \C$, $B = \spec \C$,
and $F = \R$. Let $k \ge 2$,
\[
\Eis_{\abs}^{k-2}: 
        \tEh(B) \cap \Eh(B)_{\tor} \pfeil
H^{k-1}_{\abs} (\Eh^{(k-2)} , k-1)^{\sgn} \; ,
\]
and identify the latter group, as in \ref{2J}~(a), with
\[
\Sym^{k-2} H^1 (\Eh (\C) , 2 \pi i \R) \; .
\]
Recall the map $G_{\Eh,k}$ introduced in \ref{2H}.

\begin{Thm} [cmp.~\cite{BL}, 6.4.5] \label{8G}
Let $s \in \tEh(B) \cap \Eh(B)_{\tor}$. Then
\[
\Eis_{\abs}^{k-2} (s) = \frac{k!}{k-1} \cdot G_{\Eh,k} (s) \; .
\]
\end{Thm} 

\begin{Proof} 
The factor $\frac{k!}{k-1}$ comes from normalizations, which differ from
those used in \ref{2H}. We refer to the proof of 
Proposition~\ref{13cF} for the explanation.
\end{Proof}

The description of the $\ell$-adic version of the
Eisenstein symbol on torsion is given in \cite{Ki2}, Thm.~4.2.9.

\newpage

%
%

\section{The formalism of elliptic Bloch groups} \label{II}

This section contains the main results of the present article.
Let us insist that we work under the 
hy\-po\-the\-sis $(DP)$ (Definition~\ref{9aA}), which we impose on a subset
$P \subset \Eh(B)$, relative to which our constructions are done. 
They will be performed simultaneously in the motivic and the absolute setting,
and they will be compatible under the respective regulators. \\

In Subsection~\ref{9a}, we give
a variant of the residue spectral sequences \ref{4F} and \ref{7B}
(Theorem~\ref{9aC}). This spectral sequence 
$^k\tilde{\stern}$ is central to everything to
follow. With the exception of the elliptic symbols $\{ s \}_k$, \emph{all}
the data constituting the formalism of elliptic Bloch groups are
constructed directly from the $E_1$-terms, differentials, and edge morphisms
of $^k\tilde{\stern}$ 
(Definitions~\ref{9aD}, \ref{9aE}, and \ref{9aF}). The elliptic
polylogarithm, and hence the results of Section~\ref{I}, enter the definition
of the elliptic symbols (\ref{9aG}), and the identification of
their behaviour under the differentials (Theorem~\ref{9aI}).\\

Subsection~\ref{9b} contains the results mentioned in the introduction.
We analyze the restrictions
\[
\varrho_k:
H^{k-1}_? (\Eh^{(k-2)}, k-1)^{\sgn} \pfeil \kerr(d_k) \; .
\]
In the setting of absolute cohomology, they are isomorphisms 
(Theorem~\ref{9bA}). Furthermore, in Hodge theory, we are able to establish
compatibility with \ref{2I}, thereby giving an explicit description
in terms of Eisenstein--Kronecker series of the elements in $\kerr(d_k)$
(Theorem~\ref{9bE}). For motivic cohomology, our picture is incomplete. 
If the set $P$ is contained in the torsion subgroup $\Eh(B)_{\tor}$, 
or if $k \le 3$, then we are able to
show that the $\varrho_k$ are isomorphisms (Corollary~\ref{9bC},
Corollary~\ref{9bK}). In the general case, a detailed analysis of 
the spectral sequence $^k\tilde{\stern}$ yields our Main Theorem~\ref{9bI},
which states that the $\varrho_k$ are isomorphisms if the groups
$H^i_\Mh(\Eh^m , m)_{- , \ldots , -}$
vanish in a certain range of indices. This (and the corresponding phenomenon
in absolute cohomology) motivates Conjecture~\ref{9bL}, which we see as
the elliptic analogue of the Beilinson--Soul\'e conjecture.\\

In Subsection~\ref{9c}, we discuss norm compatibility
of our formalism. \\

The wish to present our main results as concisely as possible lead us
to the decision to give the proofs of a number of results in 
Section~\ref{III}. These results are \ref{9aZ},
\ref{9aC}, \ref{9aEb}, \ref{9bA},
\ref{9bE}, and \ref{9bI}. The proofs of
\ref{9aZ} and \ref{9aEb} will constitute Subsection~\ref{6}.
The construction of the spectral sequence $^k\tilde{\stern}$ 
of \ref{9aC} will be done in
Subsections~\ref{13a} and \ref{13b}. The proofs of \ref{9bA} and
of \ref{9bI} will be given simultaneously 
(Subsection~\ref{15}), 
which is possible since the analogues of the
vanishing assumptions of \ref{9bI} are satisfied in absolute cohomology.
Finally, \ref{9bE} will be proved in \ref{13c}.

\newpage
\subsection{Construction of the elliptic Bloch groups} \label{9a}

The aim of this subsection is the geometric construction of the
data which appear in Zagier's conjecture for elliptic curves 
\[
\pi: \Eh \pfeil B \; .
\]
Let us fix a subscript $? \in \{ \Mh, \abs \}$. We shall define:
\begin{itemize}
\item[(a)] The \emph{elliptic Bloch groups} $\Bl_{k,?} = \Bl_{k,?} (\Eh)$,
$k \ge 1$.
\item[(b)] The \emph{elliptic symbols}
\[
\{ \argdot \}_k : \widetilde{\Eh} (B) \longrightarrow
\Bl_{k,?} \; .
\]
\item[(c)] The \emph{differentials}
\[
d_k : \Bl_{k,?} \longrightarrow \Bl_{k-1,?}
\otimes_{\Q} \Bl_{1,?} \quad (:= 0 \; \mbox{for} \; k=1)    
\]
mapping $\{ s \}_k$ to $\{ s \}_{k-1} \otimes \{ s \}_1$.
\item[(d)] The \emph{restrictions}
\[
\varrho_k: H^{k-1}_{?} ( \Eh^{(k-2)} , k-1 )^\sgn
                                                        \pfeil \kerr (d_k)
\]for $k \ge 2$, and
\[
\varrho_1: H^2_{?} ( \Eh , 1)_- = \Eh(B) \otimes_\BZ \BQ \pfeil 
\kerr (d_1) = \Bl_{1,?} \; .
\]
\end{itemize}
More precisely, these data will depend on the choice of a subset $P$
of the Mordell--Weil group $\Eh (B)$ satisfying the hypothesis
$(DP)$, to be introduced in \ref{9aA}.
The data $\Bl_k$, $d_k$, and $\varrho_k$ will occur in a certain
variant of the residue spectral sequences \ref{4F} and \ref{7B}, which we
describe first. Consider the group
\[
H^{\bullet}_?(\Eh^2 \times_B \Eh^{n},\argstern)_{-, \ldots ,-}^{+, \, \sgn} \; ,
\]
where the subscript refers to the actions of $[-1]$ 
on all components,
the superscript $+$ to the $(+1)$-eigenpart for the action of
$\fS_2$, and the superscript $\sgn$ to the action of $\fS_n$.
We need the following result:

\begin{Prop} \label{9aZ}
Let $? \in \{ \Mh, \abs \}$. Consider the morphisms
\[
\vartheta: \Eh^{(n)} \longinto \Eh^{n+1} 
\]
and 
\[
\Sigma: \Eh^2 \times_B \Eh^n \longonto \Eh \times_B \Eh^n \; , \;
(x_1,x_2,y) \longmapsto (x_1+x_2,y) \; .
\]
Define the map $\jmath$ as the composition of
\[
\vartheta^* \Sigma_*: 
H^{\bullet}_?(\Eh^2 \times_B \Eh^{n},\argstern)_{-, \ldots ,-}^{+, \, \sgn} 
\pfeil
H^{\bullet - 2}_? (\Eh^{(n)}, \argstern - 1) 
\]
and the projection
\[
H^{\bullet - 2}_? (\Eh^{(n)}, \argstern - 1) \longonto
H^{\bullet - 2}_? (\Eh^{(n)}, \argstern - 1)^\sgn
\]
onto the $\sgn$-eigenspace under the action of $\fS_{n+1}$. Then $\jmath$ 
is an isomorphism. 
\end{Prop}

The proof of this result will be given in Subsection~\ref{6}. \\

\begin{Def} \label{9aA} 
(a) A subset $P \subset \Eh (B)$
is said to satisfy the {\em disjointness property $(DP)$}, if any 
two unequal
sections
in $P$ have disjoint support:
\[
s , s' \in P \; ; \; s \ne s' \; \Longrightarrow \; s-s' \in \tEh(B) \; .
\]
\noindent (b) The elliptic curve $\Eh$ is said to satisfy the 
disjointness property $(DP)$ if $\Eh(B)$ satisfies $(DP)$.
\end{Def}

Of course, the condition $(DP)$ is empty if $B$ is the spectrum of
a field. Furthermore,
if an integer $N \ge 2$ is invertible on $B$, i.e., if
$[N]: \Eh \to \Eh$
is \'etale, then
\[
\Eh [N] (B) \le \Eh (B)
\]
satisfies $(DP)$.\\

Fix a subset $P \subset \Eh (B)$ satisfying $(DP)$, and 
$? \in \{ \Mh, \abs \}$. The convention of Subsection~\ref{7} will continue
to be used: In the case $?=\Mh$, the base $B$ is only assumed to be regular, 
noetherian and connected. If $?=\abs$, or if 
a statement concerns regulators, then it is
assumed that $B$ is such that absolute cohomology is defined, i.e., as in
Subsection~\ref{3}.\\

The following should be seen as the elliptic analogue 
of the groups $L_{(p)}$ of
\cite{Jeu}, page~226:

\begin{Def} \label{9aB}
Let $U \in \Ch_{\pi , P}$. The group
\[
\tH^i_?(U^r \times_B \Eh^s,j)
\]
is defined as the quotient of $H^i_?(U^r \times_B \Eh^s,j)$ 
by the sum of the images
of the cup product
\[
H^1_?(U,1) \otimes_F H^{i-1}_?(U^{r-1} \times_B \Eh^s,j-1) \pfeil 
H^i_?(U^r \times_B \Eh^s,j)
\]
in all coordinate directions.
\end{Def}

Recall the Abel--Jacobi morphism
\[
[\argdot]: \Eh(B)_F \pfeil H^2_?(\Eh,1)_- \; , \quad 
s \longmapsto [s] \; .
\]
Here we have applied the convention
\[
\ast_F := \ast \otimes_\BZ F \; ,
\]
which will be used from now on. 
More generally, for a subset $S$ of an abelian group $\FG$, we shall write
$S_F$ for the subspace of $\FG_F$ generated by $S$.
We use the notation
$[P_F]$ for the image of $P_F$ under $[\argdot]$ in $H^2_?(\Eh,1)_-$.
When $?=\Mh$, then the map
$[\argdot]$ is an isomorphism, and we identify $P_\BQ$ and $[P_\BQ]$. \\

Our variant of the residue sequence reads as follows: 

\begin{Thm} \label{9aC}
Let $k \ge 2$. 
\begin{enumerate}
\item[(a)]
There is a 
natural spectral sequence
\[
^k\tilde{\stern} \quad \quad \quad \quad ^k\tE^{p,q}_1 \Longrightarrow 
H^{2k-2+p+q}_? (\Eh^{(k-2)}, k-1)^{\sgn} \; ,
\]
where the terms $^k\tE^{p,q}_1$ are given as follows:
\begin{enumerate}
\item[(1)] If $-k \le p \le -2$, then
\[
^k\tE^{p,q}_1 = \ls{U \in \Ch_{\pi , P , -}} 
\tH^{-p+q}_?(U^2 \times_B U^{-p-2},-p)_{-, \ldots ,-}^{+, \, \sgn} \otimes_F 
\bigwedge^{k+p} [P_F] \; ,
\]
where the superscript $+$ refers to the $(+1)$-eigenspace for the action of
$\fS_2$, and the superscript $\sgn$ to the action of $\fS_{-p-2}$.
The subscript $-, \ldots ,-$ refers to the action of multiplication by $-1$
on all $-p$ components of $U^2 \times_B U^{-p-2}$.
\item[(2)] If $p = -1$, then
\[
^k\tE^{p,q}_1 = 
[P_F] \otimes_F \ls{U \in \Ch_{\pi , P , -}} 
\tH^{2k-3+q}_?(U \times_B \Eh^{k-2}, k-1)_{-, \ldots ,-}^\sgn \; .
\]
\item[(3)] If $p = 0$, then 
\[
^k\tE^{p,q}_1 = 
\Sym^2 [P_F] \otimes_F 
H^{2k-4+q}_?(\Eh^{k-2}, k-2)_{-, \ldots ,-}^\sgn \; .
\]
\item[(4)] In the cases not covered by (1)--(3), we have
$^k\tE^{p,q}_1 = 0$. 
\end{enumerate}
\item[(b)]
The differential 
$^2\partial^{-2,q}_1$ on $^2\tE^{-2,q}_1$
is induced by the direct limit of the residue maps between
\[
H^{2+q}_?(U^2 ,2)_{-,-}^+ 
\]
and 
\[
F [U_\infty(B)]_- \otimes_F 
H^{1+q}_?(U , 1)_- \; ,
\]
the summation maps $F [U_\infty(B)]_- \to P_F$, and the Abel--Jacobi map.
\item[(c)]
For $k \ge 3$, the differential 
$^k\partial^{-k,q}_1$ on $^k\tE^{-k,q}_1$
is induced by the direct limit of the residue maps between
\[
H^{k+q}_?(U^2 \times_B U^{k-2},k)_{-, \ldots ,-}^{+, \, \sgn} 
\]
and 
\[
H^{k-1+q}_?(U^2 \times_B U^{k-3},k-1)_{-, \ldots ,-}^{+, \, \sgn} \otimes_F 
F [U_\infty(B)]_- \; ,
\]
the summation maps $F [U_\infty(B)]_- \to P_F$, and the Abel--Jacobi map.
\item[(d)]
For $k \ge 2$, the edge morphism
\[
^k\varrho^q: H^{k-2+q}_? (\Eh^{(k-2)}, k-1)^{\sgn} \pfeil
\ls{U \in \Ch_{\pi , P , -}}
H^{k+q}_?(U^2 \times_B U^{k-2},k)_{-, \ldots ,-}^{+, \, \sgn}
\]
is given by the composition of 
the direct limit of the restrictions
\[
H^{k+q}_?(\Eh^2 \times_B \Eh^{k-2},k)_{-, \ldots ,-}^{+, \, \sgn} \pfeil
H^{k+q}_?(U^2 \times_B U^{k-2},k)_{-, \ldots ,-}^{+, \, \sgn} 
\]
and the inverse of the isomorphism
\[
\jmath: 
H^{k+q}_?(\Eh^2 \times_B \Eh^{k-2},k)_{-, \ldots ,-}^{+, \, \sgn} \isoto
H^{k-2+q}_? (\Eh^{(k-2)}, k-1)^{\sgn} 
\]
of \ref{9aZ}.
\item[(e)] The spectral sequence $^k\tilde{\stern}$ is compatible with
the regulators. It is 
covariantly functorial with respect to change of the set $P$,
and contravariantly functorial with respect to change of the base $B$.
\end{enumerate}
\end{Thm}

The proof of this result will be given in Subsections~\ref{13a}
and \ref{13b}. Our Main 
Theorem~\ref{9bI} will be a consequence of a 
detailed analysis of the differentials
in the spectral sequence $^k\tilde{\stern}$, in particular, of the
differential $^k\partial^{-k,1}_1$ on
\[
^k\tE^{-k,1}_1 =
\ls{U \in \Ch_{\pi , P , -}} 
\tH^{k+1}_?(U^2 \times_B U^{k-2},k)_{-, \ldots ,-}^{+, \, \sgn} \; .
\]

\begin{Def} \label{9aD}
Let $k \ge 0$.
Define the {\em $k$-th elliptic Bloch group} 
\[
\Bl_{k,P,?} := \Bl_{k,P,?} (\Eh)
\]
as follows:
\begin{enumerate}
\item[(0)] $\Bl_{0,P,?} := 0$.
\item[(1)] $\Bl_{1,P,?} := H^2_?(\Eh,1)_-$.
\item[(2)] If $k \ge 2$, then 
\[
\Bl_{k,P,?} \subset {^k\tE^{-k,1}_1} =
\ls{U \in \Ch_{\pi , P , -}} 
             \tH^{k+1}_?(U^2 \times_B U^{k-2},k)_{-, \ldots ,-}^{+, \, \sgn}
\]
is defined as the image of the composition $p$ of the following two maps: first,
the restriction from the kernel of
\[
\ls{U \in \Ch_{\pi , P , -}} 
H^{k+1}_?(\Eh \times_B U^{k-1},k)^\sgn_{-,-,\ldots,-}
\pfeil
\ls{U \in \Ch_{\pi , P , -}} H^{k+1}_?(U^k ,k)^\Sgn_{-,\ldots,-}
\]
to the kernel of 
\[
\ls{U \in \Ch_{\pi , P , -}} 
\tH^{k+1}_?(U \times_B U^{k-1},k)^\sgn_{-,-,\ldots,-}
\pfeil
\ls{U \in \Ch_{\pi , P , -}} \tH^{k+1}_?(U^k ,k)^\Sgn_{-,\ldots,-} \; ,
\]
where $\sgn$ refers to the action of $\fS_{k-1}$, and $\Sgn$ to the action of
$\fS_k$; second, the symmetrization with respect to the first two coordinates
on $\ls{U \in \Ch_{\pi , P , -}} 
\tH^{k+1}_?(U \times_B U^{k-1},k)^\sgn$.    
\end{enumerate}
\end{Def}

In order to prepare the definition of the differentials $d_k$,
consider the direct limit of the residue sequences
\[
0 \pfeil H^1_?(U,1)_- \stackrel{\res}{\pfeil} F [U_\infty]_- 
\stackrel{[\argdot]}{\longrightarrow} 
H^2_?(\Eh,1)_- \pfeil H^2_?(U,1)_- \pfeil 0 \; .
\]
Observe that $H^1_?(\Eh,1)_-$ is trivial. For $?=\Mh$, this is a consequence
of the existence 
of the canonical isomorphism
\[
\CO^*(\argdot) \otimes_\BZ \BQ \isoto H^1_\Mh(\argdot,1) \; .
\]
For $? = \abs$, one has:
\[
H^1_{\abs} (\Eh , 1)_{-} = 
\Hom_{\Sh B} (F (0) , \fH) \; ,
\]
which is trivial because of weight reasons.\\

The following is an immediate consequence of the exactness of the
above sequence:

\begin{Lem} \label{9aDa}
(a) There is a canonical isomorphism
\[
[P_F] \isoto 
\cokerr (\res: \ls{U \in \Ch_{\pi , P , -}} H^1_?(U,1)_- \pfeil F [P]_-) \; .
\]
\noindent (b) There is a canonical isomorphism
\[
\ls{U \in \Ch_{\pi , P , -}} H^2_?(U,1)_- \isoto 
\cokerr ([\argdot]: P_F \pfeil 
                    H^2_?(\Eh,1)_-) = H^2_?(\Eh,1)_- / [P_F] \; .
\]
\end{Lem}

By definition, and because of Theorem~\ref{9aC}~(c), the differential
$^k\partial^{-k,1}_1$ respects the Bloch groups if $k \ge 3$. 

\begin{Def} \label{9aE}
(a) Define the \emph{differential}
\[ 
d_2: \Bl_{2,P,?} \pfeil 
\Sym^2 \Bl_{1,P,?} = \Sym^2 H^2_?(\Eh, 1)_-
\]
as the symmetrization of the map induced (see \ref{9aEa} below)
by the composition of
\[
\ls{U \in \Ch_{\pi , P , -}} H^{3}_?(\Eh \times_B U,2)_{-, -}
\stackrel{\res}{\pfeil} H^2_?(\Eh, 1)_- \otimes_F F[P]_-
\]
and of
\[
H^2_?(\Eh, 1)_- \otimes_F F[P]_-
\stackrel{[\argdot]}{\pfeil} H^2_?(\Eh, 1)_-^{\otimes 2} \; .
\]
\noindent
(b) Let $k \ge 3$. Define the \emph{differential}
\[
d_k:= {^k\partial^{-k,1}_1}: \Bl_{k,P,?} \pfeil 
                                    \Bl_{k-1,P,?} \otimes_F [P_F] \; .
\]
\end{Def}

We still need to show: 

\begin{Lem} \label{9aEa}
The map
\[
d_2: \Bl_{2,P,?} \pfeil \Sym^2 H^2_?(\Eh, 1)_- 
\] 
is well defined.
\end{Lem}

\begin{Proof} Let us consider the direct
limit of the residue sequences
\begin{eqnarray*}
\ldots \! & \! \stackrel{\res}{\pfeil} \! & \! 
                          F [U_\infty]_- \otimes_F H^1_?(U,1)_- \\
\pfeil H^3_? (\Eh \times_B U,2)_{-,-} \pfeil H^3_? (U \times_B U,2)_{-,-} 
\! & \! \stackrel{\res}{\pfeil} \! & \! 
                          F [U_\infty]_- \otimes_F H^2_?(U,1)_- \; .
\end{eqnarray*}
Together with Lemma~\ref{9aDa}, it shows that the map 
$([\argdot] \circ \res_1, [\argdot] \circ \res_2)$
from
\[
H^3_? (U \times_B \Eh,2)_{-,-} \oplus H^3_? (\Eh \times_B U,2)_{-,-} 
\quad \text{to} \quad
\Sym^2 H^2_?(\Eh,1)_-
\]
factors first through the image of 
$H^3_? (U \times_B \Eh,2)_{-,-} \oplus H^3_? (\Eh \times_B U,2)_{-,-}$
in $H^3_? (U \times_B U,2)_{-,-}$, and then also through the image in
$\tH^3_? (U \times_B U,2)_{-,-}$. One obtains the claim by symmetrization.
\end{Proof}

Next, we have:

\begin{Prop} \label{9aEb}
The edge morphism
$^k\varrho^1$ on $H^{k-1}_? (\Eh^{(k-2)}, k-1)^{\sgn}$
lands in the Bloch group.
\end{Prop}

The proof of this result will be given in Subsection~\ref{6}.

\begin{Def} \label{9aF} 
Let $k \ge 2$. Define the \emph{restriction}
\[
\varrho_k:= {^k\varrho^1}: 
H^{k-1}_? (\Eh^{(k-2)}, k-1)^{\sgn} \pfeil \kerr(d_k) \; .
\]
\end{Def}

Now for the elliptic symbols. Their definition involves the elliptic
polylogarithm. Fix a section $s$ in $P$ unequal to (hence disjoint from) 
$i$, and form the pullback
\[
\{ s \}'_{k} := s^{\ast} \pol^{(k)} \in 
H^{k+1}_? (\Eh \times_B U^{k-1}_s , k)^{\sgn} \; ,
\]
where we set $U_s := \tEh -  s (B)$.
According to \ref{7K}, we have trivial image of $\{ s \}'_{k}$ 
under the map
\[
H^{k+1}_? (\Eh \times_B U^{k-1}_s , k)^{\sgn} \longrightarrow 
H^{k+1}_? (U^k_s , k)^{\Sgn} \; .
\]

\begin{Def} \label{9aG}
(a) The \emph{elliptic symbol}
$\{ s \}_{1} \in \Bl_{1,P,?} = H^2_?(\Eh, 1)_-$ is defined as 
\[
[s] = \halb (\{ s \}'_{1} - \{ -s \}'_{1}) = \halb([s] - [-s]) \in 
H^{2}_? (\Eh , 1)_- \; .
\]
\noindent (b) For $k \ge 2$, the \emph{elliptic symbol}
$\{ s \}_{k} \in \Bl_{k,P,?}$
is defined in two steps: first, form the projection 
$\{ s \}''_{k}$ of $\{ s \}'_{k}$ 
onto the $-, -,\ldots ,-$-eigenpart of
\[
H^{k+1}_? (\Eh \times_B (U_s \cap U_{-s})^{k-1} , k)^{\sgn} 
\]
under the action of multiplication by $-1$ on all $k$ components;
second, take the image $\{ s \}_{k}$ of $\{ s \}''_{k}$
in
\[
\ls{U \in \Ch_{\pi , P , -}} 
             \tH^{k+1}_?(U^2 \times_B U^{k-2},k)_{-, \ldots ,-}^{+, \, \sgn}
\]
under the map $p$ of Definition~\ref{9aD}.
\end{Def}

\begin{Rem} \label{9aGa}
In \cite{W5}, 4.3, we indicated why (from a sheaf theoretical point
of view) one should expect the elliptic Zagier 
machinery to produce only trivial elements in
\[
H^{k+1}_? (\Eh^k , k)^{\Sgn}_{- , \ldots , -} \; .
\]
As we see, this prediction is consistent with our present definition:
\[
\{ s \}'_k \longmapsto 0 \quad \text{in} \quad H^{k+1}_? (U^k_s , k)^{\Sgn}
\]
for any single $s$.
\end{Rem}

By construction of our data, we have:

\begin{Prop} \label{9aH}
The Bloch groups $\Bl_{k,P,?}$ (\ref{9aD}), the 
elliptic symbols $\{ \argdot \}_k$ (\ref{9aG}), the
differentials $d_k$ (\ref{9aE}) and the restrictions $\varrho_k$ (\ref{9aF})
are compatible with
the regulators, and functorial with respect to change of the set $P$,
and of the base $B$.
\end{Prop}

The identification of the images of the elliptic symbols under the
differentials is a direct consequence of their definition, and of the
construction of the elliptic polylogarithm:

\begin{Thm} \label{9aI}
Assume that $k \ge 2$, and $s \in P - \{ i \}$.
\begin{enumerate}
\item[(a)] For $k=2$,
\[ 
d_2: \Bl_{2,P,?} \pfeil 
\Sym^2 \Bl_{1,P,?} = \Sym^2 H^2_?(\Eh, 1)_-
\]
maps $\{ s \}_2$ to $\{ s \}_1 \otimes \{ s \}_1$.
\item[(b)] For $k \ge 3$,
\[
d_k: \Bl_{k,P,?} \pfeil \Bl_{k-1,P,?} \otimes_F [P_F] 
\]
maps $\{ s \}_k$ to $\{ s \}_{k-1} \otimes \{ s \}_1$.
\end{enumerate}
\end{Thm}

\begin{Proof} This follows from Proposition~\ref{7I}.
\end{Proof}

\newpage
\subsection{Statement of the main results} \label{9b}

Fix a subset $P \subset \Eh(B)$ satisfying $(DP)$.
The aim of this subsection is to study the nature of the restrictions
\[
\varrho_k:
H^{k-1}_? (\Eh^{(k-2)}, k-1)^{\sgn} \pfeil \kerr(d_k) \; .
\]
We expect them to be isomorphisms. Let us state
what we will actually be able to \emph{prove}:\\

For $? = \abs$, the sheaf theoretical interpretation of the absolute
cohomology groups will imply:

\begin{Thm} \label{9bA}
Assume that $k \ge 2$. The morphism
\[
\varrho_k: 
H^{k-1}_\abs (\Eh^{(k-2)}, k-1)^{\sgn} \pfeil \kerr(d_k) 
\]
is an isomorphism.
\end{Thm}

For the proof, we refer to Subsection~\ref{15}. 

\begin{Cor} \label{9bAa}
There is a canonical map from $\kerr (d_{k,\Mh})$, the kernel of the
differential in the motivic setting, to
\[
H^{k-1}_{\abs} (\Eh^{(k-2)} , k-1)^{\sgn} \; .
\]
Its composition with
\[
\varrho_k: 
H^{k-1}_\Mh (\Eh^{(k-2)} , k-1)^{\sgn} \longrightarrow \kerr (d_{k,\Mh})
\]
is the regulator.
\end{Cor}

One might imagine applications of this result in numerical experiments
in connection with Beilinson's conjecture on
\[
L(\Sym^{k-2} h^1 (\Eh), k-1) \; .
\]

There is a particularly easy case, where the \emph{existence} of 
the spectral sequence $^k\tilde{\stern}$ alone implies the desired
property of $\varrho_k$:

\begin{Thm} \label{9bB}
Let $k \ge 2$, and $? \in \{ \Mh,\abs \}$. Assume that $[P_F]$ is
trivial.
\begin{enumerate}
\item[(a)] The natural inclusions
\[
\kerr(d_k) \subset \Bl_{k,P,?} \subset \ls{U \in \Ch_{\pi , P , -}} 
\tH^{k+1}_?(U^2 \times_B U^{k-2},k)_{-, \ldots ,-}^{+, \, \sgn} 
\]
are equalities.
\item[(b)] The morphism
\[
\varrho_k: 
H^{k-1}_? (\Eh^{(k-2)}, k-1)^{\sgn} \pfeil \Bl_{k,P,?} 
\]
is an isomorphism.
\end{enumerate}
\end{Thm}

\begin{Proof}
Indeed, the spectral sequence $^k\tilde{\stern}$ is
concentrated in the column $p = -k$.
\end{Proof}

\begin{Cor} \label{9bC}
The conclusions of Theorem~\ref{9bB} hold if $P$ is contained in
the torsion subgroup $\Eh(B)_{\tor}$.
\end{Cor} 

A somewhat less obvious application of \ref{9bB} reads as follows:

\begin{Cor} \label{9bD}
The conclusions of Theorem~\ref{9bB} hold in the Hodge 
theoretic setting, for
\[
B = \spec A \quad \text{and} \quad F = \BR \; .
\]
\end{Cor}

\begin{Proof}
Under these conditions, the target space of the Abel--Jacobi map is
trivial (see \cite{Jn3}, Lemma~9.2).
\end{Proof}

In Subsection~\ref{13c}, we will
establish compatibility, up to scaling,
of the present construction and the one used for Theorem~\ref{2I}:

\begin{Thm} \label{9bE}
Let $k \ge 2$, and assume we are in the Hodge setting, 
with $A = \C$, $B = \spec \C$, and $F = \R$. 
Then for any $\sum_{\alpha} \lambda_{\alpha} \{ s_{\alpha} \}_k$ 
in $\kerr(d_k) = \Bl_{k,\Eh,\abs}$, we have:
\[
\varrho_k^{-1} 
\left( \sum_{\alpha} \lambda_{\alpha} \{ s_{\alpha} \}_{k} \right) = 
\frac{k!}{k-1} \cdot
\sum_{\alpha} \lambda_{\alpha} G_{\Eh ,k} (s_{\alpha})
\]
under the identification of
\[
H^{k-1}_{\abs} (\Eh^{(k-2)} , k-1)^{\sgn} = 
\Ext^1_{\MHSR} (\R (0) , \Sym^{k-2} \fH (1))
\]
and
\[
\Sym^{k-2} H^1 (\Eh (\C) , 2 \pi i \R) \; ,
\]
and with the map $G_{\Eh , k}$ of \ref{2H}.
\end{Thm}

\begin{Rem} \label{9bF}
Corollary~\ref{9bC} yields in particular
a construction of ele\-ments in motivic cohomology
\[
H^{k-1}_\Mh (\Eh^{(k-2)}, k-1)^{\sgn} 
\]
from formal linear combinations 
$\sum_{\alpha} \lambda_{\alpha} \{ s_{\alpha} \}$ of torsion sections.
For every $\C$-valued point $b$ of the base $B$, their image under
the associated regulator
map to 
\[
H^{k-1}_\abs (\Eh^{(k-2)}_b, k-1)^{\sgn} 
\]
is explicitly described by the formula given in \ref{9bE}. When $P$ is 
a finite subgroup of $\Eh(B)$,
then on
the level of absolute Hodge 
cohomology, we recover $(-\frac{k!}{k-1} \cdot (\sharp P)^{-(k-2)})$ times
Beilinson's Eisenstein symbol on torsion defined and studied
in \cite{B5} and \cite{De1} (for the
factor $(-(\sharp P)^{-(k-2)})$, see \cite{De1}, (10.9)). 
The purpose of Subsection~\ref{12}
will be to establish the relation between 
this construction and the one of Subsection~\ref{8}.
\end{Rem}

We note the following consequence of \ref{9bC}:

\begin{Cor}\label{9bG}
Let $k \ge 2$, and $? \in \{ \Mh,\abs \}$. Assume that $P$ is a
finite subset of $\Eh(B)_{\tor}$, 
and define $U \subset \Eh$ as the complement of the images
of the sections in $P$. Then the morphism
\[
\varrho_k: 
H^{k-1}_? (\Eh^{(k-2)}, k-1)^{\sgn} \pfeil 
\tH^{k+1}_?(U^2 \times_B U^{k-2},k)_{-, \ldots ,-}^{+, \, \sgn} 
\]
is an isomorphism.
\end{Cor}

\begin{Rem} \label{9bH} Let $P$ and $U$ be as in \ref{9bG}.
\begin{enumerate}
\item[(a)] It can be shown
that the
restriction induces isomorphisms
\[
H^{i}_? (\Eh^{n}, j)_{- , \ldots , -} \isoto 
\tH^{i}_?(U^n, j)_{-, \ldots ,-}
\]
for arbitrary $i,j,n$. 
This statement should be compared to \cite{DeW}, (5.2).
\item[(b)] In a similar spirit, Beilinson's definition of the Eisenstein 
symbol on torsion makes use of the identification of 
$H^{i}_? (\Eh^{n}, j)^\sgn$
with the co-invariants of $H^{i}_?(U^n, j)^\sgn$ under the action
of a certain group (see \cite{B5}, 3.1.1~(c) and
\cite{De1}, (8.16)). 
\item[(c)] Quite amusingly, the next step in Beilinson's original construction consists
in generating elements in $H^{k-1}_?(U^{k-1}, k-1)^\sgn$ 
by taking cup products of
invertible functions on $U$. By contrast, the present construction relies
heavily on dividing out the images of such cup products in
\[
H^{k+1}_?(U^2 \times_B U^{k-2},k)_{-, \ldots ,-}^{+, \, \sgn} \; .
\]
(The two $k$ correspond to one another.)
\end{enumerate}
\end{Rem}

Our main result concerns the case $? = \Mh$:

\begin{MainThm} \label{9bI}
Assume that $k \ge 2$. The morphism 
\[
\varrho_k: 
H^{k-1}_\Mh (\Eh^{(k-2)}, k-1)^{\sgn} \pfeil \kerr(d_k) 
\]
is an isomorphism if
\[
H^i_\Mh(\Eh^m , m)_{- , \ldots , -} = 0
\]
for $2 \le m \le k-2$ and $-k+2m+2 \le i \le m$:
\begin{eqnarray*}
H^i_{\Mh} (\Eh^2 , 2)_{-,-} & = & 0 \; , \quad -k + 6 \le i \le 2 \; , \\
H^i_{\Mh} (\Eh^3 , 3)_{-,-,-} & = & 0 \; , \quad -k + 8 \le i \le 3 \; ,\\
& \vdots & \\
H^{k-2}_{\Mh} (\Eh^{k-2} , k-2)_{- , \ldots , -} & = & 0 \; .
\end{eqnarray*}
\end{MainThm}

The proof of this result will constitute Subsection~\ref{15}. Since the 
conditions given above are empty for $2 \le k \le 3$, we obtain:

\begin{Cor} \label{9bK}
(a) The morphism 
\[
\varrho_2: 
\CO^* (B)_\Q \isoto H^1_\Mh (B, 1) \pfeil \kerr(d_2) 
\]
is an isomorphism.\\[0.2cm]
(b) The morphism
\[
\varrho_3:
H^2_\Mh (\Eh, 2)_- \cong H^2_\Mh (\Eh^{(1)}, 2)^\sgn \pfeil \kerr(d_3) 
\]
is an isomorphism.
\end{Cor} 

These particular cases will be discussed in more detail
in Subsections~\ref{10} and \ref{11}, respectively. 
More precisely, these subsections will contain \emph{direct} proofs
of \ref{9bK}~(a) and \ref{9bK}~(b) respectively.\\

The following should be seen 
as the \emph{elliptic analogue of the Beilinson--Soul\'e conjecture}.

\begin{Conj} \label{9bL} 
Let $\Eh \to B$ be an elliptic curve over a regular noetherian base. Then
\[
H^i_{\Mh} (\Eh^k , j)_{- , \ldots , -} = 0 
\]
for any $k \ge 0$ and $i \le \min (k , 2j-1)$.
\end{Conj}

As its ``classical'' counterpart, it is motivated by the situation
in sheaf theory: we have
\[
H^i_{\abs} (\Eh^k , j)_{- , \ldots , -} = 
\Ext^{i-k}_{\Sh B} (F (0) , \fH^{\otimes k} (j-k)) \; ,
\]
which vanishes for $i < k$. For $i=k\le 2j -1$,
\[
\Hom_{\Sh B} (F (0) , \fH^{\otimes k} (j-k)) = 0
\]
since $\fH^{\otimes k} (j-k)$ is of strictly negative weight $k-2j$.\\

From \ref{9bI}, we deduce:

\begin{Prop} \label{9bM}
If Conjecture~\ref{9bL} holds for the elliptic curve $\Eh$, then
\[
\varrho_k: 
H^{k-1}_\Mh (\Eh^{(k-2)}, k-1)^{\sgn} \pfeil \kerr(d_k) 
\]
is an isomorphism for any $k \ge 2$.
\end{Prop}

\begin{Rem} \label{9bMa}
This is the elliptic analogue of the case $p=1$ of \cite{Jeu}, 
Thm.~3.12. In fact, his condition on the nonexistence of
``low weight $K$-theory''
corresponds to Conjecture~\ref{9bL}, while the disjointness
property $(DP)$ corresponds to his condition ``if $u,v \in U$, then
$u - v$ [...] is in $\CO^*$''.
\end{Rem}

The explicit shape of the pre-images under $\varrho_k$ of elements of
\[
\kerr(d_k) \subset \Bl_{k,P,\abs}
\]
in the setting of absolute Hodge cohomology (see \ref{9bE}), and its comparison
with Beilinson's definition of the Eisenstein symbol on torsion
(see \ref{9bF}) motivates the following:

\begin{Def} \label{9bN}
Let $k \ge 2$ and  $? \in \{ \Mh , \abs \}$, and assume that
\[
\varrho_k: 
H^{k-1}_? (\Eh^{(k-2)}, k-1)^{\sgn} \pfeil \kerr(d_k) 
\]
is an isomorphism.
\[
\Eis^{k-2}_? := \varrho_k^{-1} : \kerr (d_k) \longrightarrow 
H^{k-1}_? (\Eh^{(k-2)}, k-1)^{\sgn}
\]
is called the {\em Eisenstein symbol}.
\end{Def}

Whenever it is defined, the Eisenstein symbol is (by \ref{9aH})
compatible with the regulators,
with change of the set $P$, and of the base $B$. In the setting of
absolute Hodge cohomology, its explicit shape is given by \ref{9bE}.

\begin{Rem} \label{9bO}
(a) View
$$C(P)^{\le 2}: \;\; \BQ[P] \stackrel{\delta}{\pfeil} \BQ[P] 
\otimes_{\BQ} P_\BQ \; ,$$
where $\delta((s)) := (s) \otimes s$, as a complex with entries in degrees $1$
and $2$.\\

Loosely speaking, \ref{9bI} gives a connection between $H^1$ of a quotient of
$C(P)^{\le 2}$ and a certain part of
$$H^{k-1}_{\Mh} (\Eh^{(k-2)} , k-1)^{\sgn} \; .$$
\noindent (b) A formally similar picture features in work of 
Goncharov (\cite{G}). Let $B$
be the spectrum of a field $K$, and $P = \Eh(K)$. 
If $K$ is algebraically closed, one considers the full complex 
(in degrees $\ge 1$)
$$C(\Eh (K))^{\bullet}: \;\; \Z[\Eh(K)] \stackrel{\delta}{\pfeil} 
\Z[\Eh(K)] \otimes_{\Z} \, \Eh(K) \stackrel{\delta}{\pfeil} 
\Z[\Eh(K)] \otimes_{\Z} \bigwedge^2 \Eh(K) \stackrel{\delta}{\pfeil} \ldots$$
The differential is defined as
$$\delta ((s) \otimes t_1 \wedge \ldots \wedge t_m) :=
(s) \otimes s \wedge t_1 \wedge \ldots \wedge t_m \; .$$
By dividing out certain $n$-fold convolution products of divisors of functions
on $\Eh$, one gets quotients $B_n(\Eh)$ of $\Z[\Eh(K)]$, which organize into
quotient complexes of $C(\Eh (K))^{\bullet}\,$:
$$B(\Eh,k)^{\bullet}: \;\; B_k(\Eh) \stackrel{\delta}{\pfeil}
B_{k-1}(\Eh)  \otimes_{\Z} \Eh(K) \stackrel{\delta}{\pfeil} 
B_{k-2}(\Eh)  \otimes_{\Z} \bigwedge^2 \Eh(K) \stackrel{\delta}{\pfeil}
\ldots$$
It is possible to define a subcomplex $\Dh^{\bullet}_{(k-1)}$ of the Gersten
complex for $\Eh^{(k-2)}$, and a morphism of complexes
$$\Dh^{\bullet}_{(k-1)} \pfeil B(\Eh,k)^{\bullet} \; .$$
In particular, one gets, like in our formalism, a morphism from a certain
subgroup of $K_{k-1}(\Eh^{(k-2)})^{(k-1)}$ (namely, the pre-image of 
$H^1(\Dh^{\bullet}_{(k-1)})$) to the cohomology group
$H^1(B(\Eh,k)^{\bullet})$. 
Because of Galois descent, this last statement continues to hold for arbitrary
fields $K$, if one defines
$$B(\Eh,k)^{\bullet} := 
\left( B(\Eh \times_K \bar{K},k)^{\bullet} \right)^{G_K} \; ,$$
and if one works modulo torsion.
For the details, we refer to \cite{G}, Section 7.\\[0.2cm]
(c) In fact, \cite{G} is not only concerned with the special motivic cohomology
groups
$$H^{k-1}_{\Mh} (\Eh^{(k-2)} , k-1)^{\sgn} \; .$$
Goncharov conjectures (loc.\ cit., Conjecture 9.5.b)) that there are canonical
isomorphisms
$$H^i(B(\Eh,k)^{\bullet} \otimes_{\Z} \Q) \silo 
H^{k+i-2}_{\Mh} (\Eh^{(k-2)} , k-1)^{\sgn} \; $$
for any $i$.
\end{Rem}

\begin{OP} \label{9bP}
Relate the two approaches.
The construction of $B_k(\Eh) \otimes_{\Z} \Q$ 
as a quotient of the space of $K$-rational
divisors on $\Eh$, and the shape of the differential $\delta$ suggest a 
relation between the vector spaces $B_k(\Eh) \otimes_{\Z} \Q$ 
and the fixed part under the Galois group of the
limit over all Galois extensions $K'/K$ of $\Bl_{k,\Eh(K'),\Mh} 
(\Eh \times_K K')$. In fact, one might hope for a relation between the whole
complex $B(\Eh,k)^{\bullet}$, and the row $q=1$ of the spectral sequence 
$^{k} \tilde{\stern}$ of \ref{9aC}.
The connection might be provided by work of Levin (\cite{L}):
there, an explicit description of the elliptic motivic polylogarithm 
is given in terms
of collections of divisors, and symbols on the divisors.
\end{OP}

\newpage
\subsection{Norm compatibility} \label{9c}

In order to reduce Parts~1 and 2 of the weak version of
Zagier's conjecture for elliptic curves
to Conjecture~\ref{9bL} (see \ref{9cE}), 
we have to study a last structural property
of our construction, the so-called \emph{norm compatibility}.\\

We fix the following situation:
\[
\psi : \Eh_1 \pfeil \Eh_2
\]
is an \'etale isogeny, and $P_2 \subset \Eh_2 (B)$ a subset with the 
disjointness property $(DP)$. The pre-image 
$\psi^{-1} (P_2) \subset \Eh_1 (B)$ of 
$P_2$ under $\psi : \Eh_1 (B) \to \Eh_2 (B)$ also satisfies disjointness.
Finally fix a subset $P_1 \subset \psi^{-1} (P_2)$.\\

Recall that the motivic and absolute cohomology theories are covariant
under finite morphisms $f$. If $f$ is a finite Galois covering with group
$G$, then $f^*$ identifies the target space with the $G$-invariants of the
source, and under this identification, the \emph{norm map} $\tr_f = f_*$ equals
the sum of the $g_*$ over $g \in G$. Thus, $f_* f^*$ equals multiplication
by the order of $G$.

\begin{Thm} \label{9cA}
(a) The norm map $\tr_{\psi^k}$ induces a morphism
\[
\tr_{\psi^k} : \Bl_{k,P_1 ,?} (\Eh_1) \longrightarrow \Bl_{k,P_2,?} (\Eh_2) \; .
\]
It is compatible with the regulators, and functorial with respect to change of 
$P_2$ and $P_1 \subset \psi^{-1} (P_2)$, 
and of the base $B$.\\[0.2cm]
(b) For $k \ge 2$, the diagram
\[
\vcenter{\xymatrix@R-10pt{ 
H^{k-1}_? (\Eh^{(k-2)}_1 , k-1)^{\sgn} 
        \ar[r]^-{\varrho_k} \ar[d]_{\deg (\psi) \cdot \psi^{k-2}_{\ast}} & 
\Bl_{k,P_1,?} (\Eh_1) \ar[d]^{\tr_{\psi^k}} \\
H^{k-1}_? (\Eh^{(k-2)}_2 , k-1)^{\sgn} \ar[r]^-{\varrho_k} &
\Bl_{k,P_2,?} (\Eh_2) \\}}    
\]
commutes.\\[0.2cm]
(c) The norm is compatible with the differential $d_{k}$ 
in the following sense: 
\begin{enumerate}
\item[(1)]
The diagram
\[
\vcenter{\xymatrix@R-10pt{ 
\Bl_{2,P_1,?} (\Eh_1) \ar[r]^-{d_2} \ar[d]_{\tr_{\psi^k}} & 
\Sym^2 \Bl_{1,P_1,?} = \Sym^2 H^2_?(\Eh_1, 1)_- 
                          \ar[d]^{\Sym^2 \psi} \\
\Bl_{2,P_2,?} (\Eh_2) \ar[r]^-{d_2} &
\Sym^2 \Bl_{1,P_2,?} = \Sym^2 H^2_?(\Eh_2, 1)_- \\}}    
\]
commutes.
\item[(2)]
For $k \ge 3$, the diagram
\[
\vcenter{\xymatrix@R-10pt{ 
\Bl_{k,P_1,?} (\Eh_1) \ar[r]^-{d_k} \ar[d]_{\tr_{\psi^k}} & 
\Bl_{k-1,P_1,?} (\Eh_1) \otimes_F [P_{1,F}]  
                          \ar[d]^{\tr_{\psi^{k-1}} \otimes \psi} \\
\Bl_{k,P_2,?} (\Eh_2) \ar[r]^-{d_k} &
\Bl_{k-1,P_2,?} (\Eh_2) \otimes_F [P_{2,F}] \\}}    
\]
commutes.
\end{enumerate}
\end{Thm}

\begin{Proof} (a) follows from the definitions, and
(c) follows from the compatibility of norms with residues.
It remains to show (b). We claim that the diagram
\[
\vcenter{\xymatrix@R-10pt{ 
H^{k+1}_?(\Eh^2_1 \times_B \Eh^{k-2}_1 , k)_{-, \ldots ,-}^{+, \, \sgn} 
             \ar[r]^-{\jmath} \ar[d]_{\psi^k_\ast} &        
H^{k-1}_? (\Eh^{(k-2)}_1 , k-1)^{\sgn}
          \ar[d]^{\deg (\psi) \cdot \psi^{k-2}_\ast}  \\
H^{k+1}_?(\Eh^2_2 \times_B \Eh^{k-2}_2 , k)_{-, \ldots ,-}^{+, \, \sgn} 
             \ar[r]^-{\jmath} &
H^{k-1}_? (\Eh^{(k-2)}_2 , k-1)^{\sgn} \\}}    
\]
commutes. Here, the morphism $\jmath$ is the one occurring
in the definition of $\varrho_k$ (see \ref{9aC}~(d) and 
\ref{9aZ}). It was defined as the anti-symmetrization of
$\vartheta^* \Sigma_*$, for
\[
\Sigma: \Eh^{k}_j \longonto \Eh^{k-1}_j \; , \quad
(x_1,x_2,y) \longmapsto (x_1+x_2,y) \; ,
\]
and 
\[
\vartheta: \Eh^{(k-2)}_j \longinto \Eh^{k-1}_j \; .
\]
Note that Galois descent allows us to assume that
$\kerr(\psi)$  consists of sections of $\Eh_1$.
The first map $\Sigma_*$ commutes with the norms. The factor $\deg (\psi)$
comes from $\vartheta^*$: define $\Eh^{[k-2]}_1 \subset \Eh^{k-1}_1$ by the
cartesian diagram
\[
\vcenter{\xymatrix@R-10pt{ 
\Eh^{[k-2]}_1 \ar@{^{(}->}[r] \ar[d] &        
\Eh^{k-1}_1 \ar[d]^{\psi^{k-1}}   \\
\Eh^{(k-2)}_2 \ar@{^{(}->}[r]^-{\vartheta} &
\Eh^{k-1}_2  \\}}    
\]
Because of base change for
cohomology, the claim follows from the observation that $\Eh^{[k-2]}_1$
consists of $\deg (\psi)$ copies of $\Eh^{(k-2)}_1$.
\end{Proof}

In order to study the behaviour of $\kerr (d_k)$ under norms, 
we need the following:

\begin{Prop} \label{9cB}
Assume that $k \ge 2$, and that the 
$\varrho_l (\Eh_i)$ are isomorphisms for $4 \le l \le k$ and $i = 1,2$.
Then the norm
\[
\tr_{\psi^k} : \Bl_{k,P_1,?} (\Eh_1) \longrightarrow \Bl_{k,P_2,?} (\Eh_2)
\]
is injective.
\end{Prop}

\begin{Proof}
Use induction on $k$, \ref{9cA}~(b) and (c), and the fact 
that $\psi^{k-2}_{\ast}$ is an isomorphism on 
$H^{k-1}_? (\Eh^{(k-2)} , k-1)^{\sgn}$ (\cite{BL}, Lemma~5.1.2).
\end{Proof}

The behaviour of the elliptic symbols under norms is as follows:

\begin{Thm} \label{9cC}
Assume that $\kerr (\psi)$ consists of sections of $\Eh_1$, and that the
subset $P_1$ of $\Eh_1 (B)$ is invariant under translation by
elements of $\kerr (\psi) (B)$. Let
$s_1 \in P_1$ be disjoint from $\kerr (\psi)$.
Then under the norm map
\[
\tr_{\psi^k} : \Bl_{k,P_1,?}(\Eh_1) \longrightarrow \Bl_{k,P_2,?}(\Eh_2)
\]
of \ref{9cA}, the element
\[
\sum_{t \in \kerr (\psi) (B)} \{ s_1 + t \}_k
\]
is mapped to $\deg (\psi) \cdot \{ \psi (s_1) \}_k$. 
\end{Thm}

\begin{Proof}
This follows from the construction of the elliptic symbols, and from
Proposition~\ref{7J}~(b).
\end{Proof}

From the preceding results, we deduce formally what we called 
\emph{norm compatibility for $\{\argdot \}_k$ with respect to $d_k$} in 
\cite{W5}, 1.6:

\begin{Cor} \label{9cD}
Assume that $k \ge 2$, and that the 
$\varrho_l (\Eh_i)$ are isomorphisms for $4 \le l \le k$ and $i = 1,2$.
Furthermore, assume that 
$\kerr (\psi)$ consists of sections of $\Eh_1$, and that the
subset $P_1$ of $\Eh_1 (B)$ is invariant under translation by
elements of $\kerr (\psi) (B)$.
Let $\lambda_{\alpha} \in F$, and $s_{1,\alpha} \in P_1$ disjoint from 
$\kerr (\psi)$. Then
\[
d_k \left( \sum_{\alpha} \lambda_{\alpha} 
\{ \psi (s_{1,\alpha}) \}_k \right) = 0 
\]
if and only if
\[
d_k \left( \sum_{\alpha} \lambda_{\alpha} 
\sum_{t \in \kerr (\psi) (B)} \{ s_{1,\alpha} + t \}_k \right) = 0 \; ,
\]
and if this is the case, then the equality
\[
\psi^{k-2}_{\ast} \circ \Eis^{k-2}_? \left( \sum_{\alpha} \lambda_{\alpha} 
\sum_{t \in \kerr (\psi) (B)} \{ s_{1,\alpha} + t \}_k \right) = 
\Eis^{k-2}_? \left( \sum_{\alpha} \lambda_{\alpha} 
\{ \psi (s_{1,\alpha}) \}_k \right)
\]
holds in $H^{k-1}_? (\Eh^{(k-2)}_2 , k-1)^{\sgn}$.
\end{Cor}

Let us stress that by \ref{9bC}, the conclusion of \ref{9cD} holds
if $P$ is contained in $\Eh(B)_{\tor}$. 
For Beilinson's original Eisenstein symbol on torsion
(\cite{B5}, 3.1), the corresponding statement is already known: see e.g.\
the proof of \cite{W5}, Thm.~1.9.1.\\

Furthermore, the above result, together with \ref{9bM}, \ref{9aH}, \ref{9aI}
and \ref{9bE} yields:

\begin{Thm} \label{9cE}
(a) Conjecture~\ref{9bL} implies Parts~1 and 2 of the weak version of
Zagier's conjecture for elliptic curves (\cite{W5}, Conj.~1.6.(B))
satisfying the disjointness property $(DP)$.\\[0.2cm]
(b) For $k \le 3$, Parts~1 and 2 of the weak version of
Zagier's conjecture hold for elliptic curves 
satisfying the disjointness property $(DP)$.
\end{Thm}

\newpage

%
%

\section{The proofs} \label{III}

This section contains the proofs of the results stated in Section~\ref{II},
as well as amplifications of these results in the special cases
$k=2$ (Subsection~\ref{10}), $k=3$ (Subsection~\ref{11}), and the 
torsion case (Subsection~\ref{12}); 
the material contained in \ref{10}, \ref{11}, and \ref{12} 
will however not be needed in the rest of the article, and relies
only on the subsections with smaller cardinal numbers. \\

Subsection~\ref{6} contains the proofs of Propositions~\ref{9aZ}
and \ref{9aEb}. The main ingredient is the analysis of the
decomposition of a motive of (a power of) an elliptic curve $\Eh$ relative to
a base $B$, which was carried out in \cite{BL}, Section~5.
Corollary~\ref{6Ea} yields an identification of the $n$-th symmetric 
power of $h^1 (\Eh)$ inside some 
$h^\bullet (\Eh^k)$, which so far seems not to have been used in the
literature, but which will serve best our purposes.\\

We then proceed to construct 
the spectral sequence $^k\tilde{\stern}$ of \ref{9aC}. It comes about
as a quotient of a residue spectral sequence of the type discussed in
Subsection~\ref{4}. 
Since the category of spectral
sequences is not abelian, we need to have abstract criteria which ensure that
such a quotient exists. Subsection~\ref{13a} gives two such criteria:
Proposition~\ref{13aA} and Corollary~\ref{13aG}. In order to apply
\ref{13aA}, we have to construct and recognize acyclic complexes. 
This is where we use
one of the central ideas of de 
Jeu's construction of the Bloch groups in the
setting of the classical Zagier conjecture (\cite{Jeu}, pages 222--226).
We should mention that in the old version of the
present article (``On the generalized 
Eisenstein symbol'', Spring 1997), we tried to form the quotient by passing to
the coinvariants under the action of the Mordell--Weil group on the 
residue spectral sequence. This is a process which a priori destroys all
exactness properties, and this is why we had to impose an injectivity
hypothesis on the (product of the) regulators. 
Thanks to Proposition~\ref{13aB}, which
generalizes de Jeu's construction, this hypothesis is no longer needed.
Subsection~\ref{13b} contains the application to residue spectral
sequences.\\

In Subsection~\ref{14}, we continue the study of motives 
and motivic cohomology of
elliptic curves $\Eh$ started in Subsection~\ref{6}. The main result 
(Theorem~\ref{14Ba}) states that translations by sections act trivially on
the eigenspace
\[
H^\bullet_? (\Eh, \argstern )_- \; .
\]
As a consequence, we obtain that
the exterior cup product
\[
\Eh(B) \otimes_\BZ H^\bullet_? (\Eh,\argstern)_{-}
\pfeil
H^{\bullet +2}_? (\Eh^{2} , \argstern +1)_{-,-} 
\]
lands in the $\sgn$-eigenspace of 
$H^{\bullet +2}_? (\Eh^{2} , \argstern +1)_{-,-}$ (Corollary~\ref{14D}).
This result will be responsible for the vanishing of 
many edge morphisms in the spectral sequence $^k\tilde{\stern}$.
(This in turn ensures that the vanishing assumptions on the groups
\[
H^i_\Mh(\Eh^m , m)_{- , \ldots , -} 
\]
are used in a more economical way than in the old version of this
article.) \\

We now have at our disposal all the ingredients needed for the proofs
of Theorem~\ref{9bA} and of Main Theorem~\ref{9bI}. They are given in
Subsection~\ref{15}.\\

Finally, Subsection~\ref{13c} connects the present \emph{geometric}
construction of the formalism of elliptic Bloch groups to the
\emph{sheaf theoretical} one of \cite{W5}, thereby
proving Theorem~\ref{9bE}.

\newpage
\subsection{Motives generated by elliptic curves. I} \label{6}

The aim of this subsection is to provide proofs of Proposition~\ref{9aZ}
and Proposition~\ref{9aEb}.\\

Let $B$ be a regular and noetherian scheme. In \cite{BL}, Section 5, the 
authors introduce, for a fixed family of elliptic curves $\pi : \Eh \to B$, the 
$\Q$-tensor category $\Mh^{(\Eh)}$
of \emph{relative motives for abelian schemes isogenous to a power of $\Eh$}. 
Morphisms are defined by so-called {\it linear cycles} 
only (loc.\ cit., 5.3). These correspond to classes of abelian subschemes.
By loc.\ cit., 5.4.3, the category $\Mh^{(\Eh)}$ is semisimple and abelian.

\begin{Def}\label{6A}
Let $a \neq 0$ and $n \ge 0$ be integers. Denote by 
$[a]^n$ the isogeny given by 
multiplication by $a$ on $\Eh^n$. 
\begin{itemize}
\item[(a)] Denote by $h^\bullet(\Eh^n) \in \Mh^{(\Eh)}$ the 
relative motive associated 
to $\Eh^n$.
\item[(b)] Define the endomorphism $\tr_{[a]^n}$ on $h^\bullet(\Eh^n)$ by
\[
\tr_{[a]^n} := [a]^n_{\ast} \; .
\]
\end{itemize}
\end{Def}

The functor $h^\bullet$ is \emph{contravariant}.
The morphisms in $\Mh^{(\Eh)}$ 
are enough to get the standard results about the K\"unneth 
decomposition (loc.\ cit., 5.4.5)
-- see \cite{DeM}, \cite{K} for the case when the base $B$ is smooth and 
quasi-projective over a field. 
In particular, one has:

\begin{Prop} \label{6B}
Let $n \ge 0$. 
\begin{itemize}
\item[(a)] The direct summand $h^r(\Eh^n) \subset h^\bullet(\Eh^n)$ can
be characterized as the eigenpart for the value $a^r$ of the action of
$\tr_{[a]^n}$, for any integer $a$, whose absolute value is at least $2$.
\item[(b)] There is a canonical isomorphism of ring objects in $\Mh^{(\Eh)}$
\begin{eqnarray*}
h^\bullet(\Eh^n) 
& \isoto & \bigwedge^\bullet h^1 (\Eh^n) =  
\bigwedge^\bullet \left( h^1 (\Eh)^n \right) \\
& = & \left( \bigwedge^\bullet h^1 (\Eh) \right)^{\! \!\otimes n} =  
\left( \BQ(0) \oplus h^1 (\Eh) \oplus \BQ(-1) \right)^{\otimes n} \; .
\end{eqnarray*}
\item[(c)] The isomorphism in (b) is compatible with the actions 
of $\tr_{[a]^n}$
and $\tr_{[a]^1}$, respectively.
In particular, multiplication by a non-zero 
integer $a$
on $\Eh^n$
induces an isomorphism on the relative motive
$h^\bullet(\Eh^n)$.
\end{itemize}
\end{Prop}

\begin{Proof}
(a) is part of the construction of \cite{BL}, 5.4. (b) is loc.~cit.,
5.4.5. For (c), we refer to loc.~cit., Lemma~5.1.1.
\end{Proof}

We note another direct consequence of \cite{BL}, 5.4.5:

\begin{Prop} \label{6C}
We have
\[
h^1 (\Eh) = h^\bullet (\Eh)_- \; ,
\]
the $(-1)$-eigenspace for the action of $[-1]$.
\end{Prop}

Combining \ref{6B}~(b) and \ref{6C}, we get:

\begin{Cor} \label{6D}
There is a canonical isomorphism 
\[
h^\bullet(\Eh^n)_{-, \ldots , -} 
\isoto h^1 (\Eh)^{\otimes n} \; ,
\]
where the subscript refers to the actions of $[-1]$ on the components.
\end{Cor}  

Let us identify the $n$-th symmetric power of $h^1 (\Eh)$ inside some 
$h^\bullet (\Eh^k)$.
The classical way of doing this is to define
\[
\Eh^{(n)} \subset \Eh^{n+1}
\]
as the kernel of the summation map. One has an action of $\fS_{n+1}$, and an 
isomorphism, canonical up to sign,
\[
\Sym^n h^1 (\Eh) \silo h^\bullet (\Eh^{(n)})^{\sgn}
\]
(see \cite{De1}, 8.7 if $B$ is the spectrum of a field, the proof of \cite{BL},
6.2.1 in the general case, and the discussion in \cite{W5}, 1.2 for the
possible normalizations of the isomorphism).\\

For our purposes, another realization of $h^\bullet (\Eh^{(n)})^{\sgn}$
will be useful:

\begin{Cor} \label{6Ea}
(a) There is a canonical isomorphism 
\[
\jmath: h^\bullet(\Eh^2 \times_B \Eh^n)^{+, \, \sgn}_{-, \ldots , -} 
\otimes \BQ(1)
\isoto  h^\bullet (\Eh^{(n)})^\sgn \; ,
\]
where the subscript refers to the actions of $[-1]$ on the $n+2$ components,
the superscript $+$ to the $(+1)$-eigenpart for the action of
$\fS_2$, and the superscripts $\sgn$ to the actions of $\fS_n$
and $\fS_{n+1}$, respectively.\\[0.2cm]
(b) The isomorphism $\jmath$ equals the composition of $\vartheta^* \Sigma_*$,
where
\[
\Sigma: \Eh^2 \times_B \Eh^n \longonto \Eh \times_B \Eh^n \; , \;
(x_1,x_2,y) \longmapsto (x_1+x_2,y) \; ,
\]
and the projection onto the $\sgn$-eigenpart of $h^\bullet (\Eh^{(n)})$ under
the action of $\fS_{n+1}$. \\[0.2cm]
(c) There is an isomorphism, canonical up to sign,
\[
h^\bullet(\Eh^2 \times_B \Eh^n)^{+, \, \sgn}_{-, \ldots , -}
\otimes \BQ(1) 
\isoto \Sym^n h^1 (\Eh) \; .
\]
\end{Cor}

\begin{Proof}
This is a direct consequence of \ref{6D}.
\end{Proof}
\quad \\

\begin{Proofof}{Proposition~\ref{9aZ}}
This is the statement on the level of absolute
and motivic cohomology corresponding to
Corollary~\ref{6Ea}.
\end{Proofof}

Another consequence of \ref{6B}~(b) and \ref{6D} reads as follows:

\begin{Cor} \label{6Eb}
Consider the $\sgn$-eigenpart
\[
h^\bullet(\Eh \times_B \Eh^{n+1})^{\sgn}_{-, -, \ldots , -} = 
h^1(\Eh) \otimes \Sym^{n+1} h^1(\Eh)
\]
for the action of $\fS_{n+1}$ on 
$h^\bullet(\Eh \times_B \Eh^{n+1})_{-, -, \ldots , -}$, and the 
$\Sgn$-eigenpart
\[
h^\bullet(\Eh^{n+2})^{\Sgn}_{-, \ldots , -} = \Sym^{n+2} h^1(\Eh)
\]
for the action of $\fS_{n+2}$. Then the motives
\[
h^\bullet(\Eh^2 \times_B \Eh^n)^{+, \, \sgn}_{-, \ldots , -} 
\cong \BQ(-1) \otimes \Sym^n h^1 (\Eh) 
\]
(where the $\sgn$ here refers to the action of $\fS_n$) and 
\[
\kerr \left( \pr_\Sgn: 
h^\bullet(\Eh \times_B \Eh^{n+1})^{\sgn}_{-, -, \ldots , -} 
\pfeil h^\bullet(\Eh^{n+2})^{\Sgn}_{-, \ldots , -}
\right)
\]
coincide as submotives of $h^\bullet(\Eh^{n+2})_{-, \ldots , -}$.
\end{Cor}

\begin{Proofof}{Proposition~\ref{9aEb}}
This follows from the statement on the level of absolute
and motivic cohomology corresponding to Corollary~\ref{6Eb}.
\end{Proofof}

\newpage
\subsection{The case $k=2$} \label{10}

The purpose of this subsection is to give proofs of Theorem~\ref{9aC},
Theorem~\ref{9bA} and Main Theorem~\ref{9bI} for $k=2$. First,
we hope that the reader will thus get a flavour of the idea of the proof
in the general case. More concretely, one of the results 
(\ref{10K}) will actually
enter the proof of \ref{9bI} for $k = 3$.\\

We keep the conventions of the Sections~\ref{I} and \ref{II}. In particular,
$P \subset \Eh(B)$ satisfies $(DP)$ (see \ref{9aA}), and
$? \in \{ \Mh , \abs \}$.\\

Let $U \in \Ch_{\pi , P}$. Let us consider the filtration introduced
before \ref{4E}, in the case $n=2$:
\[
F_p \Eh^2 := \{ (x_1 , x_2) \in \Eh^2 \tei \mbox{at most $2+p$ 
coordinates lie in} \; U_{\infty} \} \; .
\]
Using the notation of Subsection~\ref{3}, we have
\[
\emptyset = F_{-3} \Eh^2 \subset F_{-2} \Eh^2 = U^2 \subset 
F_{-1} \Eh^2 = \Eh^2_{\reg} \subset F_0 \Eh^2 = \Eh^2 \; .
\]
There is a version ``before sign eigenspaces'' of 
the residue spectral sequences \ref{4F} and \ref{7B}.
Its direct limit over all $U \in \Ch_{\pi , P , -}$ looks as follows:

\begin{Prop} \label{10A}
There is a 
natural residue spectral sequence
$$^2\stern \quad \quad \quad \quad ^2E^{p,q}_1 \Longrightarrow 
H^{4+p+q}_? (\Eh^2 , 2)_{-,-} \; ,$$
$$^2E^{p,q}_1 = \bigoplus
\ls{U \in \Ch_{\pi , P , -}} H^{-p+q}_? (U^{-p} , -p)_{-,\ldots,-} \otimes_F 
F [P]_-^{\otimes (2+p)} \; ,$$
where the sum $\oplus$ runs through all subsets of $\{ 1,2 \}$ of cardinality
$-p$. 
\end{Prop}

For fixed $q$, let us consider the complex $^2E^{\bullet,q}_1$. It is
the direct limit of complexes $^2E^{\bullet,q}_1(U)$ of the following form,
concentrated in degrees $-2, -1, 0$:
\[
H^{2+q}_? (U^2 , 2)_{-,-} 
\! \pfeil \!\!\!\!
\begin{array}{c}
H^{1+q}_? (U , 1)_- \otimes_F F [U_\infty (B)]_- \\
\oplus \\
F [U_\infty (B)]_- \otimes_F H^{1+q}_? (U , 1)_- 
\end{array}
\!\!\!\! \pfeil \!
H^{q}_? (B , 0) \otimes_F 
F [U_{\infty} (B)]_-^{\otimes 2} 
\]
By construction, the differentials are given by the residue maps.\\

Next, let us consider a certain
quotient complex $^2\tE^{\bullet,q}_1$
of $^2E^{\bullet,q}_1$, namely, the direct limit of
the complexes $^2\tE^{\bullet,q}_1(U)$:
\[
\tH^{2+q}_? (U^2 , 2)_{-,-} 
\! \pfeil \!\!\!\!
\begin{array}{c}
\tH^{1+q}_? (U , 1)_- \otimes_F [U_\infty(B)_F] \\
\oplus \\

[U_\infty(B)_F] \otimes_F \tH^{1+q}_? (U , 1)_-
\end{array}
\!\!\!\! \pfeil \!
H^{q}_? (B , 0) \otimes_F 
[U_\infty(B)_F]^{\otimes 2} \; .
\]
Recall that we use the notation $U_\infty(B)_F$ for the subspace of $P_F$
generated by $U_\infty(B)$. Thus, $[U_\infty(B)_F]$ is the image of
$F [U_{\infty} (B)]_-$ under the Abel--Jacobi map $[ \argdot ]$.\\

The statement and the proof of the following
is modeled after \cite{Jeu}, pages 222--226:

\begin{Prop} \label{10B}
The canonical projection
\[
^2E^{\bullet,q}_1 \longonto {^2\tE^{\bullet,q}_1}
\]
is a quasi-isomorphism.
\end{Prop}

\begin{Proof}
We work on finite level, i.e., we fix $U \in \Ch_{\pi , P , -}$.
We want do define a certain subcomplex $F^{\bullet,q}$
of $^2E^{\bullet,q}_1(U)$. Let
\[
F^{-2,q} \subset {^2E^{-2,q}_1(U)}
\]
be defined as the sum of the images of the cup products
\[
H^1_?(U,1)_- \otimes_F H^{1+q}_?(U,1)_-
\]
(in the two possible directions). Note that each individual cup product
is a monomorphism, a canonical left inverse being given by the residue --
recall that
\[
\res: H^1_?(U,1)_- \pfeil F [U_\infty (B)]_-
\]
is injective.
Define
\[
C^{-1,q} \subset {^2E^{-1,q}_1(U)}
\]
as the direct sum 
\[
\begin{array}{c}
H^1_? (U , 1)_- \otimes_F H^q_?(B,0) \otimes_F F [U_\infty (B)]_-  \\
\oplus \\
F [U_{\infty} (B)]_- \otimes_F H^1_? (U , 1)_- \otimes_F H^q_?(B,0)
\end{array}
\]
Finally, define $F^{\bullet,q}$ to be the smallest subcomplex of
$^2E^{\bullet,q}_1(U)$ containing $F^{-2,q}$ and $C^{-1,q}$.\\

Let us introduce some abbreviations: let
\[
K' := \kerr 
\left( [\argdot]: F[U_\infty(B)]_- \pfeil H^2_?(\Eh,1)_- \right) \; .
\]
Because of the localization sequence, 
the residue induces an isomorphism
\[
H^1_?(U,1)_- \isoto K' \subset K := F[U_\infty(B)]_- \; .
\]
If we write $H^q(B)$ for $H^q_?(B,0)$ and $H^{1+q}(U)$ for $H^{1+q}_?(U,1)_-$,
then the complex $F^{\bullet,q}$ has the following shape:
\[
\begin{array}{c}
K' \otimes_F H^{1+q}(U) \\
+ \\
H^{1+q}(U) \otimes_F K'
\end{array} 
\stackrel{\partial}{\pfeil}
\begin{array}{c}
K' \otimes_F H^q(B) \otimes_F K\\
\oplus \\
K \otimes_F H^q(B) \otimes_F K'
\end{array}
+ \imm(\partial)
\stackrel{\partial}{\pfeil}
\imm(\partial)
\]
It is easy to see that we have 
\[
F^{\bullet,q} = 
\kerr \left( ^2E^{\bullet,q}_1 \longonto {^2\tE^{\bullet,q}_1} \right) \; .
\]
All that remains to be shown is acyclicity of $F^{\bullet,q}$.
As for injectivity at $F^{-2,q}$, one sees from the above representation
first that 
\[
\kerr(\partial) \subset K' \otimes_F H^q(B) \otimes_F K' \; ,
\]
and then that $\kerr(\partial) = 0$. So consider $F^{-1,q}$. Again, one
shows directly that
\[
\kerr(\partial) \cap 
\left(
\begin{array}{c}
K' \otimes_F H^q(B) \otimes_F K\\
\oplus \\
K \otimes_F H^q(B) \otimes_F K'
\end{array}
\right)
\]
is contained in
\[
\Delta(K' \otimes_F H^q(B) \otimes_F K') \subset 
\left(
\begin{array}{c}
K' \otimes_F H^q(B) \otimes_F K\\
\oplus \\
K \otimes_F H^q(B) \otimes_F K'
\end{array}
\right) \quad .
\]
\end{Proof}

\begin{Cor} \label{10C}
There is a 
natural spectral sequence
$$^2\tilde{\stern} \quad \quad \quad \quad ^2\tE^{p,q}_1 \Longrightarrow 
H^{4+p+q}_? (\Eh^2 , 2)_{-,-} \; ,$$
$$^2\tE^{p,q}_1 = \bigoplus
[P_F]^{\otimes (2+p)} \otimes_F 
\ls{U \in \Ch_{\pi , P , -}} \tH^{-p+q}_? (U^{-p} , -p)_{-,\ldots,-} \; ,$$
where the sum $\oplus$ runs through all subsets of $\{ 1,2 \}$ of cardinality
$-p$. 
\end{Cor}

\begin{Proof}
Indeed, from the $E_2$-terms onwards, this is precisely the
spectral sequence of \ref{10A}.
\end{Proof}

In particular, we get a proof of Theorem~\ref{9aC} for $k=2$ by passing
to the $\fS_2$-invariants of the spectral sequence
$^2\tilde{\stern}$, and applying the isomorphism
\[
\jmath: H^{4+p+q}_? (\Eh^2 , 2)_{-,-}^+ \isoto H^{2+p+q}_? (B , 1)
\]
of \ref{9aZ}.

\begin{Rem} \label{10D}
In fact, exactly the same method works for other twists. Although we
shall not need the result in this generality, we mention it for the
sake of completeness:\\

Let $j \in \BZ$. There is a 
natural spectral sequence
$$E^{p,q}_1 \Longrightarrow 
H^{4+p+q}_? (\Eh^2 , j)_{-,-} \; ,$$
$$E^{p,q}_1 = \bigoplus
[P_F]^{\otimes (2+p)} \otimes_F 
\ls{U \in \Ch_{\pi , P , -}} \tH^{-p+q}_? (U^{-p} , -2-p+j)_{-,\ldots,-} \; ,$$
where the sum $\oplus$ runs through all subsets of $\{ 1,2 \}$ of cardinality
$-p$. 
\end{Rem}

Let us study the lines $q=1$ and $q=0$ of $^2\tilde{\stern}$. As before,
the columns are indexed by $-2,-1,0$:
\[
\vcenter{\xymatrix@R-10pt{ 
\ld \tH^3_? (U^{2} , 2)_{-,-}
        \ar[r]^-{\partial_1} \ar@{-->}[drr]_{\partial_2} & 
[P_F] \otimes_F \left( H^2_? (\Eh , 1)_-/[P_F] \right)^{\oplus 2} \ar[r] & 
\quad \quad \quad \\
\quad \quad \quad \quad \quad \quad \quad \; \; \ar[r] &
0 \ar[r] &
[P_F]^{\otimes 2} \\}}    
\]
Here, we used \ref{9aDa}~(b) to identify the term $^2\tE^{-1,1}_1$.
The term $^2\tE^{-1,0}_1$ is trivial because $\tH^1_?(U,1)_-$ is.
The differential $\partial_2$ is only
defined on the kernel of $\partial_1$.

\begin{Prop} \label{10E}
Restriction to the $U^2$ induces a canonical isomorphism
\[
H^3_? (\Eh^2 , 2)_{-,-} \isoto \kerr(\partial_2) \; .
\]
\end{Prop}

\begin{Proof}
This follows from the spectral sequence, together with the fact that
$^2\tE^{0,-1}_1 = 0$.
\end{Proof}

Let us consider the situation after passing to the $\fS_2$-invariants,
and together with the Bloch groups:
\[
\vcenter{\xymatrix@R-10pt{ 
\Bl_{2,P,?} \ar[r]^-{d_2} \ar[d]_{=} & 
\Sym^2 H^2_? (\Eh , 1)_- \ar[d]^{\can} & \\
\Bl_{2,P,?} \ar[r] \ar@{_{(}->}[d] & 
H^2_? (\Eh , 1)_- \otimes_F \left( H^2_? (\Eh , 1)_-/[P_F] \right) & \\
\ld \tH^3_? (U^{2} , 2)_{-,-}^+
        \ar[r]^-{\partial_1} \ar@{-->}[drr]_{\partial_2} & 
[P_F] \otimes_F \left( H^2_? (\Eh , 1)_-/[P_F] \right) \ar[r] \ar@{^{(}->}[u] & 
\quad \quad \quad \quad \\
\quad \quad \quad \quad \quad \quad \quad \; \; \ar[r] &
0 \ar[r] &
\Sym^2 [P_F] \\}}    
\]

By definition of the differential $d_2$, the diagrams are commutative.

\begin{Prop} \label{10F}
(a) $\kerr(d_2) = ( \kerr(\partial_2) \cap \Bl_{2,P,?} ) = 
\kerr(\partial_2)$.\\[0.1cm]
(b) The sequence
\[
0 \pfeil H^1_? (B , 1) \stackrel{\varrho_2}{\pfeil} \Bl_{2,P,?}
\stackrel{d_2}{\pfeil} \Sym^2 H^2_? (\Eh , 1)_-
\]
is exact.\\[0.1cm]
(c) The image of $d_2$ contains $\Sym^2 [P_F]$.
\end{Prop}

\begin{Proof}
The above diagram shows that 
\[
\kerr(d_2) = ( \kerr(\partial_2) \cap \Bl_{2,P,?} ) \: .
\]
By \ref{10E}, the kernel of $\partial_2$ equals 
\[
H^3_? (\Eh^2 , 2)_{-,-}^+ \stackrel{\ref{9aZ}}{=} H^1 (B , 1) \; ,
\]
which in turn maps to $\kerr(d_2)$ via $\varrho_2$.
This shows (a) and (b). Part (c) follows from the existence of the
elliptic symbols, and from \ref{9aI}~(a).
\end{Proof}

This proves Theorem~\ref{9bA} and Main Theorem~\ref{9bI} for $k=2$.  

\begin{Rem} \label{10G}
(a) The sequence in \ref{10F}~(b) should be compared to 
the sequence (6) of \cite{GL}.\\[0.2cm]
(b) In Subsection~\ref{13c}, we will connect the
\emph{geometric} construction of the present article to the
\emph{sheaf theoretical} one of \cite{W5}, Sections 2 and 3
(Proposition~\ref{13cF}).
Consequently, the isomorphism
\[
\Eis^0_?: \kerr(d_2) \isoto H^1_? (B , 1)
\]
is compatible, up to a factor $2$, 
with the construction of \cite{W6} (see page 296 of the
introduction of loc.~cit.). In particular, \cite{W6}, Thm.~3.11 
gives an explicit description of
\[
\halb \cdot \Eis^0_\Mh: \kerr(d_2) \isoto \CO^\ast (B)_\BQ = H^1_\Mh (B , 1)
\]
in terms of complex analytic functions.
This description, together with \cite{W5}, Prop.~1.9.2
shows \emph{a posteori}, that the construction is also equivalent
to the one sketched in \cite{W5}, 1.9, which makes use of the Poincar\'e
line bundle, and is independent of the rest of \cite{W5}.
\end{Rem}

\begin{OP} \label{10Ga}
Find a \emph{conceptual} explanation of the connection to \cite{W5},
1.9. We expect
\cite{Ki2}, Thm.~4.1.3 to play a central role; there, the 
elliptic polylogarithm is related to the generalized Picard scheme.
As a consequence (loc.~cit., 4.2.3), one shows that the restriction
of $\Eis^0$ to torsion points can be expressed in terms of elliptic units.
\end{OP}

We note a formal consequence of \ref{10F}~(c):

\begin{Cor} \label{10H}
The exterior cup product
\[
\Sym^2 [P_F] \pfeil H^4_?(\Eh^2 , 2)_{-,-}^+
\]
is trivial.
\end{Cor}

\begin{Proof}
In fact, the exterior cup product is the edge morphism
\[
^2\tE^{0,0,+}_1 \pfeil H^4_?(\Eh^2 , 2)_{-,-}^+
\]
in the spectral sequence $^2\tilde{\stern}^+$. But by \ref{10F}~(c),
the term $^2\tE^{0,0,+}_3$ is trivial.
\end{Proof}

For a more general statement concerning the vanishing of exterior cup
products, we refer to Theorem~\ref{14A}.\\

We end this subsection by collecting all the information
we have con\-cer\-ning the
terms occurring in the spectral sequence $^2\tilde{\stern}^+$.
For the purposes of Subsection~\ref{11}, all we need
will be the vanishing of
\[
\ls{U \in \Ch_{\pi , P , -}} \tH^2_? (U^{2} , 2)_{-,-}^+ \; .
\]

\begin{Prop} \label{10I}
If $? = \abs$, then assume that we are in the setting of Hodge theory.
\begin{enumerate}
\item[(a)] The groups
\[
H^i_?(\Eh^2 , 2)_{-,-}^+ \quad \text{and} \quad
\ls{U \in \Ch_{\pi , P , -}} \tH^i_? (U^{2} , 2)_{-,-}^+
\]
are trivial for $i \ne 3 , 4$.
\item[(b)] The restriction map
\[
H^4_?(\Eh^2 , 2)_{-,-}^+ \pfeil
\ls{U \in \Ch_{\pi , P , -}} \tH^4_? (U^{2} , 2)_{-,-}^+
\]
is bijective.
\item[(c)] The differential
\[
\partial_1: \ls{U \in \Ch_{\pi , P , -}} \tH^3_? (U^{2} , 2)_{-,-}^+
\pfeil [P_F] \otimes_F \left( H^2_? (\Eh , 1)_-/[P_F] \right)
\]
is surjective.
\end{enumerate}
\end{Prop}

\begin{Proof}
First, we claim that the restriction induces an isomorphism
\[
H^i_?(\Eh^2 , 2)_{-,-}^+ \pfeil
\ls{U \in \Ch_{\pi , P , -}} \tH^i_? (U^{2} , 2)_{-,-}^+
\]
for $i \ne 3,4$, and a surjection for $i = 4$.
Using the spectral sequence $^2\tilde{\stern}^+$,
this follows formally from the fact that
\[
^2\tE^{-1,q,+}_1 = 0 \quad \text{for} \quad q \ne 1 \; ,
\]
\[
^2\tE^{0,q,+}_1 = 0 \quad \text{for} \quad q \ne 0 \; .
\]
This is because of 
\[
H^n_? (U , 1) = 0 \quad \text{for} \quad n \ne 1,2 \; ,
\]
\[
H^n_? (B , 0) = 0 \quad \text{for} \quad n \ne 0 \; .
\]
For $? = \Mh$, this is \cite{Sou}, Prop.~1 and
Thm.~4 (and because there is no nontrivial $K$-theory in the
range of negative indices). In the Hodge theoretic
setting, this follows e.g.\ from the interpretation of 
$H_\abs$ as Yoneda-Ext groups in Saito's category
of mixed Hodge modules (\cite{HW}, Cor.~A.1.10 and
Thm.~A.2.7), and the theory of weights (see \cite{S2},
(4.5.3)). Alternatively,
compute $H_\abs$ via the Leray spectral sequence, use
the analogue of the generalized Weil conjecture (\cite{S2},
(4.5.2)), and the following facts about the category $\MHS$
of mixed polarizable $F$-Hodge structures: (1)~It is of cohomological
dimension $1$; (2)~There are no nontrivial one-extensions of $F(0)$
by objects of weights $\ge 0$ (see e.g. \cite{HW}, Thm.~1 on page 297). \\

For the same reasons,
\[
H^i_?(\Eh^2 , 2)_{-,-}^+ = H^{i-2}_? (B , 1)
\]
is trivial for $i \ne 3,4$. It remains to show that 
\[
^2\tE^{0,0,+}_\infty = {^2\tE^{-1,1,+}_2} = 0 \; .
\]
For the first term, we refer to \ref{10H}. For the second, observe that
the map
\[
^2\tE^{-1,1,+}_\infty = {^2\tE^{-1,1,+}_2} \pfeil H^4_?(\Eh^2 , 2)_{-,-}^+
\]
is again induced by the symmetrization of the exterior cup product on
\[
[P_F] \otimes_F \left( H^2_? (\Eh , 1)_-/[P_F] \right) \; .
\]
The claim thus follows from \ref{14A}.
\end{Proof}

In the situation considered in \ref{10I}, there are very few non-zero
$E_1$-terms in $^2\tilde{\stern}^+$, and all differentials are surjective.
We describe the situation by the following diagram, which is concentrated in
the range $-2 \le p \le 0$, $0 \le q \le 2$. All terms which do not occur
are trivial. 
\[
\vcenter{\xymatrix@R-10pt{ 
\ld \tH^4_? (U^{2} , 2)_{-,-}^+ & & \\
\ld \tH^3_? (U^{2} , 2)_{-,-}^+
        \ar@{->>}[r]^-{\partial_1} \ar@{-->>}[drr]_{\partial_2} & 
[P_F] \otimes_F \left( H^2_? (\Eh , 1)_-/[P_F] \right) & 
\quad \quad \quad \quad \\
\quad \quad \quad \quad \quad \quad \quad  &  &
\Sym^2 [P_F] \\}}    
\]
In particular, the only nontrivial $E_3$-terms occur in the column
$p= -2$.\\

In the $\ell$-adic setting, the proof of \ref{10I}
does not work: The $\ell$-adic
co\-ho\-mo\-lo\-gy groups admit no interpretation in terms of Yoneda-Ext groups,
except for small indices. However, part of the argument runs through,
yielding the following:

\begin{Prop} \label{10K}
(a) The groups
\[
H^i_?(\Eh^2 , 2)_{-,-}^+ \quad \text{and} \quad
\ls{U \in \Ch_{\pi , P , -}} \tH^i_? (U^{2} , 2)_{-,-}^+
\]
are trivial for $i \le 2$.\\[0.2cm]
(b) The restriction map
\[
H^4_?(\Eh^2 , 2)_{-,-}^+ \pfeil
\ls{U \in \Ch_{\pi , P , -}} \tH^4_? (U^{2} , 2)_{-,-}^+
\]
is injective.
\end{Prop}

\begin{Proof}
This is because of
\[
H^n_?(\argdot , 1) = 0 \quad \text{for} \quad n \le 0 \; ,
\]
\[
H^{n}_?(\argdot , 0) = 0 \quad \text{for} \quad n \le -1 \; .
\]
\end{Proof}

\newpage
\subsection{Quotient spectral sequences} \label{13a}

Suppose given a spectral sequence
\[
E^{p,q}_1 \Longrightarrow E^{p+q} \; .
\]
The aim of this subsection is to give conditions under which
one may divide out subcomplexes $F^{\bullet,q}_1$
of the $E^{\bullet,q}_1$, so as to
obtain a quotient spectral sequence
\[
\tE^{p,q}_1 = (E^{p,q}_1 / F^{p,q}_1) \Longrightarrow \tE^{p+q} \; .
\]
One obvious condition was already used in the proof of \ref{10C}:

\begin{Prop} \label{13aA}
In the above situation, if the complexes $F^{\bullet,q}_1$ are acyc\-lic,
then there is a natural spectral sequence
\[
\tE^{p,q}_1 \Longrightarrow E^{p+q} \; .
\]
\end{Prop}

\begin{Proof}
Indeed, from the $E_2$-terms onwards, this coincides with the spectral
sequence $E^{p,q}_1 \Longrightarrow E^{p+q}$.
\end{Proof}

In order to apply \ref{13aA} in our geometric context, the
following criterion will be useful:

\begin{Prop} \label{13aB}
Let $n \ge 0$ an integer, and $K$ a field. Assume given the following data:
\begin{enumerate}
\item[(1)] a $K$-vector space $A$, together with a subspace $\Fa$.
\item[(2)] for any subset $J \subset \{ 1,\ldots,n \}$, a
$K$-vector space $B^J$.
\item[(3)] for any subset $J \subset \{ 1,\ldots,n \}$, and any $m \in J$,
a morphism
\[
\partial_m^J: B^J \pfeil B^{J - \{ m \} }\otimes_K A 
\]
and a morphism
\[
\epsilon_m^J: B^{J - \{ m \} }\otimes_K {\Fa} \pfeil B^J \; .
\]
\end{enumerate}
For a subset $I \subset \{ 1,\ldots,n \}$, we write $A^{\otimes I}$
for
\[
\bigotimes_{n \in I} A^{\otimes \{n\}} \; ,
\]
where $A^{\otimes \{n\}}$ is a copy of $A$. 
The same convention applies to $\Fa$.
Suppose that the data (1)--(3) satisfy the following hypotheses:
\begin{enumerate}
\item[(a)] $\partial_m^J \circ \epsilon_m^J = 
                       \id_{ B^{J - \{ m \} }\otimes_K {\Fa}^{\otimes \{m\}} }$.
\item[(b)] The $\partial_m^J$ commute; 
more precisely, for $m \ne m' \in J$, the diagram
\[
\vcenter{\xymatrix@R-10pt{ 
B^J \ar[rrr]^-{\partial_m^J} \ar[dd]_{\partial_{m'}^J} &&& 
B^{J-\{m\}} \otimes_K A^{\otimes \{m\}} 
              \ar[dd]^{\partial_{m'}^{J-\{m\}} \otimes \id_{A^{\otimes \{m\}}} } \\
&&& \\
B^{J-\{m'\}} \otimes_K A^{\otimes \{m'\}} 
    \ar[rrr]_-{\partial_{m}^{J-\{m'\}} \otimes \id_{A^{\otimes \{m'\}}} } &&&
B^{J-\{m,m'\}} \otimes_K A^{\otimes \{m,m'\}}   \\}}
\]
commutes.
\item[(c)] The $\partial_m^J$ and $\epsilon_m^J$ commute up to a sign;
more precisely, for $m \ne m' \in J$, there exists $c \in \{ \pm 1 \}$
such that the diagram
\[
\vcenter{\xymatrix@R-10pt{ 
B^{J-\{m\}} \otimes_K \Fa^{\otimes \{m\}} 
          \ar[rrr]^-{\partial_{m'}^{J-\{m\}} \otimes \id_{\Fa^{\otimes \{m\}}} }
          \ar[dd]_{\epsilon_m^J} &&& 
B^{J-\{m,m'\}} \otimes_K \Fa^{\otimes \{m\}} \otimes_K A^{\otimes \{m'\}} 
          \ar[dd]^{\epsilon_{m}^{J-\{m'\}} \otimes \id_{A^{\otimes \{m'\}}} } \\
&&& \\
B^J \ar[rrr]_-{\partial_{m'}^J} &&&
B^{J-\{m'\}} \otimes_K A^{\otimes \{m'\}}   \\}}
\]
commutes up to a factor $c$.
\item[(d)] The $\epsilon_m^J$ commute up to a sign;
more precisely, for $m \ne m' \in J$, there exists $c \in \{ \pm 1 \}$
such that the diagram
\[
\vcenter{\xymatrix@R-10pt{ 
B^{J-\{m,m'\}} \otimes_K \Fa^{\otimes \{m,m'\}}
      \ar[rrr]^-{\epsilon_{m'}^{J-\{m\}} \otimes \id_{\Fa^{\otimes \{m\}}} }
      \ar[dd]_{\epsilon_{m}^{J-\{m'\}} \otimes \id_{\Fa^{\otimes \{m'\}}} } &&&
B^{J-\{m\}} \otimes_K \Fa^{\otimes \{m\}} \ar[dd]^{\epsilon_m^J}  \\
&&& \\
B^{J-\{m'\}} \otimes_K \Fa^{\otimes \{m'\}} \ar[rrr]_-{\epsilon_{m'}^J} &&&
B^J  \\}}
\]
commutes up to a factor $c$.
\end{enumerate}
Define the complex $E^\bullet$ as follows:
\[
E^i := \bigoplus \left( B^J \otimes_K A^{\otimes I} \right) \; ,
\]
where the sum $\oplus$ runs through all decompositions 
$I \coprod J = \{ 1,\ldots,n \}$ such that $\sharp I = i$. 
The differential
$\partial^i$ on $E^i$ is defined as follows:
\[
\partial^i : B^J \otimes_K A^{\otimes I} \pfeil
           \bigoplus_{m \in J} 
           \left( B^{J-\{m\}} \otimes_K A^{\otimes (I \coprod \{m\})} \right)
\]
equals 
\[
\left( (-1)^{pos_J(m)} \partial_m^J \otimes \id_{A^{\otimes I}} 
                          \right)_{m \in J} \; ,
\]
where $pos_J(m) \ge 0$ denotes the position of $m$ in $J$, in other words,
the image of $m$ under the unique order-preserving bijection
\[
J \silo \{ 0,\ldots,\sharp J - 1 \} \; .
\]                          
Define ${\Fb}^J$ as the subspace of $B^J$ generated by the
images of the $\epsilon^J_m$, and the subcomplex $F^\bullet$ of $E^\bullet$
as the smallest subcomplex containing the 
${\Fb}^J \otimes_K A^{\otimes I}$:
\begin{eqnarray*}
F^i & := & \imm(\partial^{i-1} \tei F^{i-1}) + 
\bigoplus \left( {\Fb}^J \otimes_K A^{\otimes I} \right) \\
    & \stackrel{(c)}{=} & \bigoplus \left( 
          {\Fb}^J \otimes_K A^{\otimes I}
           + B^J \otimes_K \kerr (A^{\otimes I} \pfeil (A/\Fa)^{\otimes I})
                     \right) \; .
\end{eqnarray*}
Then the following holds:
\begin{enumerate}
\item[(i)] $E^i / F^i = \bigoplus \left( (B^J / {\Fb}^J) \otimes_K 
(A / {\Fa})^{\otimes I} \right)$.
\item[(ii)] The complex $F^\bullet$ is acyclic.
\end{enumerate}
\end{Prop}

\begin{Rem} \label{13aC} 
This is the abstract principle behind \cite{Jeu}, 
pages 222--226.
\end{Rem}

For the proof of \ref{13aB}~(ii), the following will be needed:

\begin{Lem} \label{13aBa}
Let $n \ge 0$ an integer,
\[
0 \pfeil V' \pfeil V \pfeil V'' \pfeil 0
\]
an exact sequence of $K$-vector spaces. Consider the complex $C_n^\bullet$
given by
\[
C_n^i := \bigoplus \left(
                       (V')^{\otimes J} \otimes_K V^{\otimes I}
                   \right) \; ,
\]
where the sum $\oplus$ runs through all decompositions 
$I \coprod J = \{ 1,\ldots,n \}$ such that $\sharp I = i$. 
The map
\[
C_n^n=V^{\otimes n} \pfeil (V'')^{\otimes n} \]
induces a quasi-isomorphism
\[
C_n^\bullet \pfeil (V'')^{\otimes n}[-n] \; .
\]
\end{Lem}

\begin{Proof}
This holds for $n \le 1$. For $n \ge 2$, observe that
\[
C_n^\bullet \cong (C_1^\bullet)^{\otimes n} \; .
\]
\end{Proof}

\begin{Proofof}{Proposition~\ref{13aB}}
(i) follows directly from the definition. For (ii), we consider $A$ and
$\Fa$ as being fixed, while we would like to study subsystems
\[
( (B')^\ast ) \subset (B^\ast)
\]
which are respected by the $\partial_m^J$ and $\epsilon_m^J$,
as well as quotients of such systems. We think of the
data $\Fb^J$, $E^\bullet$, and $F^\bullet$ as being functions of $(B^\ast)$:
\[
\Fb^J = \Fb^J(B^\ast) \; , \quad E^\bullet = E^\bullet(B^\ast) \; , \quad
                                 F^\bullet = F^\bullet(B^\ast) \; .
\]
We may assume that the $B^J$ are finite dimensional, and proceed by induction on
\[
\sum_J \dim_K B^J \; .
\]
Choose $J_0$ in such a way that $B^{J_0} \ne 0$, but $B^J = 0$ for all
proper subsets $J$ of $J_0$. Thus we have $\Fb^{J_0} = 0$. Choose a non-zero
element $b$ of $B^{J_0}$, and consider the subsystem 
$\langle b \rangle ^\ast \subset B^\ast$
generated by $b$. We have $\langle b \rangle^J =0$ if $J_0 \not\subset J$, and
\[
\langle b \rangle^J =
   \epsilon_{j_1} \circ \ldots \circ \epsilon_{j_k} 
  (b \otimes_K \Fa^{\otimes (J - J_0) }) \cong \Fa^{\otimes k}
\]
if $J = J_0 \coprod \{j_1,\ldots,j_k\}$.
Observe that we have the equality 
\[
\Fb^J(\langle b \rangle^\ast) = \langle b \rangle ^J
\]
except for $J = J_0$, where $\Fb^{J_0} = 0$. In any case, we have
\[
\Fb^J(\langle b \rangle^\ast) = \langle b \rangle ^J \cap \Fb^J \; .
\]
Similarly,
\[
F^i(\langle b \rangle ^\ast) = E^i(\langle b \rangle ^\ast)
\]
except for $i = i_0 := n - \sharp J_0$, where 
\begin{eqnarray*}
E^{i_0}(\langle b \rangle ^\ast) & = & b \otimes_K A^{\otimes i_0} \; , \\
F^{i_0}(\langle b \rangle ^\ast) & = & b \otimes_K 
                      \kerr (A^{\otimes i_0} \pfeil (A/\Fa)^{\otimes i_0}) \; .
\end{eqnarray*}
In any case, we have
\[
F^\bullet(\langle b \rangle ^\ast) = 
                   E^\bullet(\langle b \rangle ^\ast) \cap F^\bullet \; .
\]
We obtain an exact sequence of complexes 
\[
0 \pfeil F^\bullet(\langle b \rangle ^\ast) \pfeil F^\bullet \pfeil
F^\bullet(B^\ast / \langle b \rangle ^\ast) \pfeil 0 \; .
\]
It remains to show that $F^\bullet(\langle b \rangle ^\ast)$
is acyclic. But this follows from
the above description of $F^i(\langle b \rangle ^\ast)$, and from
Lemma~\ref{13aBa}.    
\end{Proofof}

We will need another condition on the existence of quotient
spectral sequences. For this,
it turns out to be convenient to use the framework of exact couples,
of which we recall the basic notions
and results. We follow the treatment and notation
of \cite{Wei}, 5.9, except that we use cohomological numbering.
In particular, the roles of the super- and 
the subscripts will be interchanged.

\begin{Def} \label{13aD}
(a) A \emph{(bigraded) exact couple} consists of a quintuple
$(D,E;i,j,k)$, where 
\[
D = D^{\bullet,\bullet} = \oplus_{p,q \in \BZ} D^{p,q} \quad \text{and} \quad
E = E^{\bullet,\bullet} = \oplus_{p,q \in \BZ} E^{p,q}
\]
are bigraded abelian groups, and $i$, $j$, and $k$ are morphisms of bidegree
$(-1,1)$, $(0,0)$, and $(1,0)$ respectively, which fit into a triangle
\[
\vcenter{\xymatrix@R-10pt{ 
D \ar[rr]^i && D \ar[dl]^j \\
& E \ar[ul]^k &\\}}
\]
which is exact in the obvious sense.\\[0.1cm]
(b) The \emph{differential} associated to a bigraded exact couple
$(D,E;i,j,k)$ is the morphism of bidegree $(1,0)$
\[
\partial := j \circ k : E \pfeil E \; .
\]
\end{Def}

It is easy to check that $\partial^2=0$. 
In an obvious fashion, one defines morphisms between exact couples.\\

There is an operation which assigns
to an exact couple $\CC = (D,E;i,j,k)$ another exact couple
\[
\CC' = (D',E';i',j',k')
\]
called the \emph{derived couple} of $\CC$ (see \cite{Wei}, p.~154): $D'$ is
defined as the image of $i$, and $E'$ as the cohomology group with respect
to $\partial$. The morphisms $i'$ and $k'$ in the triangle
\[
\vcenter{\xymatrix@R-10pt{ 
D' \ar[rr]^{i'} && D' \ar[dl]^{j'} \\
& E' \ar[ul]^{k'} &\\}}
\]
are directly induced by $i$ and $k$. The definition of $j'$ is more complicated:
for $x \in D'$, choose $y \in D$ such that $i(y) = x$. Then the class of
$j(y)$ in $E'$ is independent of the choice of $y$. One defines
\[
j'(x) := j(y) \; .
\]
This process can be iterated. We obtain a sequence of exact couples
\[
\CC_n = (D_n,E_n;i_n,j_n,k_n) \; , \; n = 1,2,\ldots 
\]
(with $\CC_1 = \CC$ and $\CC_2 = \CC'$).
Note that $i_n$ and $k_n$ are still of bidegree 
$(-1,1)$ and $(1,0)$ respectively, while $j_n$ is of bidegree $(n-1,-n+1)$.
Consequently, the differential $\partial_n$ is of bidegree $(n,-n+1)$.

\begin{Def} \label{13aE}
Let $\CC = (D,E;i,j,k)$ be an exact couple.
\begin{enumerate}
\item[(a)] We say that $\CC$ is \emph{bounded
below} if for each $n$ there is an integer $f(n)$ such that $D^{p,q} = 0$
whenever $p < f(p+q)$.
\item[(b)] We say that $\CC$ is \emph{regular} if for each $p$ and $q$
the differentials $\partial_n^{p,q}$ are zero for all large $n$.
\end{enumerate}
\end{Def}

\begin{Prop} \label{13aF}
Let $(D_1,E_1;i_1,j_1,k_1)$
be an exact couple, which is 
bounded below and regular. Assume that the morphisms
\[
i : D^{p,q}_1 \pfeil D^{p-1,q+1}_1
\]
are isomorphisms if $p \ge 1$. Then there is a natural spectral sequence
\[
E^{p,q}_1 \Longrightarrow E^{p+q} = D^{0,p+q}_1 \; .
\]
It converges, and is functorial in $(D_1,E_1;i_1,j_1,k_1)$.
\end{Prop}

\begin{Proof} This is \emph{not} \cite{Wei}, Theorem~5.9.7, but rather 
loc.~cit., Exercice~5.9.2 in cohomological notation.
\end{Proof}

In the situation of \ref{13aF}, we refer to $E^{p,q}_1 \Longrightarrow E^{p+q}$
as the \emph{spectral sequence associated to the exact couple
$(D_1,E_1;i_1,j_1,k_1)$}. We will be interested in the functorial
behaviour of the spectral sequence in a very particular case:

\begin{Cor} \label{13aG}
Let
\[
\vcenter{\xymatrix@R-10pt{ 
\left( F^{p,q}_1 \Longrightarrow F^{p+q} \right) \ar[d]^\iota \\
\left( E^{p,q}_1 \Longrightarrow E^{p+q} \right) \\}}
\]
be the morphism associated to a split monomorphism of exact couples.
Then $\iota$ is a split monomorphism of spectral sequences. In particular,
the quotient spectral sequence
\[
\tE^{p,q}_1 \Longrightarrow \tE^{p+q}
\]
exists.
\end{Cor}

\begin{Proof}
Indeed, the derivative of a split monomorphism of exact couples is still
a split monomorphism.
\end{Proof}

\newpage
\subsection{Residue sequences. II} \label{13b}

The aim of this section is to give a proof of Theorem~\ref{9aC}.
In order to do so, we have to apply the abstract material of the
previous subsection to the geometric setting. The following turns out
to be convenient:

\begin{Def} \label{13bA}
(a) Let $X$ be a scheme. A \emph{filtered system of length $\le \vert r \vert$
(with $r \le 0$) of subschemes} is a diagram $\FF$:
\[
\vcenter{\xymatrix@R-10pt{ 
G_0 &
F_0 \ar@{~>}[l] &
F_{-1} \ar@{~>}[l] \ar@{=}[dl] &
G_0 \ar@{~>}[l] \\
G_{-1} &
F_{-1} \ar@{~>}[l] &
F_{-2} \ar@{~>}[l] \ar@{=}[dl] &
G_{-1} \ar@{~>}[l] \\
G_{-2} \ar@{.}[dd] &
F_{-2} \ar@{~>}[l] \ar@{.}[dd] &
F_{-3} \ar@{~>}[l] \ar@{.}[dd] &
G_{-2} \ar@{~>}[l] \ar@{.}[dd] \\
&
&
 \\
G_r &
F_r \ar@{~>}[l] &
F_{r-1} \ar@{~>}[l] &
G_r \ar@{~>}[l] \\}}
\]
where:
\begin{enumerate}
\item[(1)] The $G_i$ and $F_i$ are locally closed subschemes of $X$.
\item[(2)] In each line, the situation
\[
\vcenter{\xymatrix@1{ 
G_{i} &
F_{i} \ar@{~>}[l] &
F_{i-1} \ar@{~>}[l] &
G_{i} \ar@{~>}[l] \\}}
\]
is complementary in the sense that one of the following holds:
\begin{itemize}
\item[(i)] $F_{i-1}$ is an open subscheme of $F_i$, with complement $G_i$.
\item[(ii)] $G_i$ is an open subscheme of $F_{i-1}$, with complement $F_i$.
\item[(iii)] $F_i$ is an open subscheme of $G_i$, with complement $F_{i-1}$.
\end{itemize}
In each line, we shall use a straight arrow 
for the open immersion.
\item[(3)] $G_r = F_r$.
\end{enumerate}
\noindent (b) A filtered system of subschemes as above is called
\emph{real} if in (a)~(2), one is always in the situation (i):
\[
\vcenter{\xymatrix@1{ 
G_{i} &
F_{i} \ar@{~>}[l] &
F_{i-1} \ar[l] &
G_{i} \ar@{~>}[l] \\}}
\]
(c) Let $X$ be smooth over a noetherian
base $S$. A filtered system of subschemes as above 
is called \emph{pure}, if the $G_i$ and $F_i$ are smooth over $S$, and if the
closed immersions in (a)~(2) are of pure codimension.
\end{Def}

\begin{Rem} \label{13bB}
There is a notion of morphism between
filtered systems of subschemes of a fixed scheme $X$. Since it
will be applied only once, and in a very particular 
setting (see Example~\ref{13bC}
below), we do not define it
in a rigid fashion, but rather rely on the reader's ability 
to imagine the scenario: The
above definition is motivated by the desire to obtain (and then, to
organize -- see below) a system of
residue sequences. There are three types of morphisms which induce
(by pullback) the maps in such sequences:
the open immersion (which we think of as being covariant), the closed immersion
(which we think of as being contravariant), and the (virtual) residue morphism
from the closed to the open subscheme (which we think of as being covariant).
Morphisms of filtered systems consist of diagrams, where
each edge is of one of these three types.
\end{Rem}

\begin{Ex} \label{13bC}
Let $k=1$, and consider a real filtered system
\[
\vcenter{\xymatrix@1{ 
G_{0} &
F_{0} \ar@{~>}[l] &
F_{-1} \ar[l] &
G_{0} \ar@{~>}[l] \\}}
\]
(the term $F_{-2} = \emptyset$ is understood).
Then there is a morphism of filtered systems
\[
\vcenter{\xymatrix@R-10pt{ 
G_{0} \ar@{~>}[dd] &
F_{0} \ar@{~>}[l] \ar@{~>}[dd] &
F_{-1} \ar[l] \ar[dd] &
G_{0} \ar@{~>}[l] \ar@{~>}[dd] \\
&&& \\
F_{-1} &
\emptyset \ar@{~>}[l] &
F_{-1} \ar@{~>}[l] &
F_{-1} \ar[l] \\}}
\]
\end{Ex}

\begin{Def} \label{13bCa}
Let $X$ be smooth over a noetherian base $S$, and $\FF$
a pure filtered system of subschemes of $X$. Let $i \le 0$. 
\begin{enumerate}
\item[(a)] The codimension of $G_i$ in $F_i$, denoted by
\[
\codim_{F_i} G_i \; ,
\]
is defined
as follows:
\begin{itemize}
\item[(i)] the usual codimension of the closed subscheme $G_i$ in
$F_i$ in 
case \ref{13bA}~(a)~(2)~(i).
\item[(ii)] minus the usual codimension of the closed subscheme $F_i$ in
$F_{i-1}$ in case \ref{13bA}~(a)~(2)~(ii).
\item[(iii)] zero in case \ref{13bA}~(a)~(2)~(iii).
\end{itemize}
\item[(b)] The codimension of $F_{i-1}$ in $F_i$, denoted by
\[
\codim_{F_i} F_{i-1} \; ,
\]
is defined
as follows:
\begin{itemize}
\item[(i)] zero in case \ref{13bA}~(a)~(2)~(i).
\item[(ii)] minus the usual codimension of the closed subscheme $F_i$ in
$F_{i-1}$ in case \ref{13bA}~(a)~(2)~(ii).
\item[(iii)] the usual codimension of the closed subscheme $F_{i-1}$ in
$G_i$ in case \ref{13bA}~(a)~(2)~(iii).
\end{itemize}
\item[(c)] The codimension of $G_i$ in $F_0$
is defined as 
\[
\codim_{F_0} G_i := \codim_{F_i} G_i + \sum_{i=0}^{i-1} \codim_{F_i} F_{i-1} \; .
\]
\end{enumerate}
\end{Def}

We recall the conventions on our base scheme $B$: if a statement concerns
motivic cohomology, then $B$ is understood to be regular, noetherian and 
connected. In the context of absolute cohomology, $B$ is 
as in Subsection~\ref{3}.

\begin{Prop} \label{13bD}
Let $? \in \{ \Mh , \abs \}$, and $\FF$ a pure filtered system
of length $\le \vert r \vert$
of subschemes of a smooth $B$-scheme $X$. 
There is a natural spectral sequence
\[
E_1^{p,q} = H_?^{-2c_p + p + q}(G_p, -c_p+j) \Longrightarrow 
H_?^{p+q}(F_0,j) \; ,
\]
where $c_p := \codim_{F_0} G_p$. We set $E_1^{p,q} := 0$ for
$p < r$ or $p > 0$.
\end{Prop}

\begin{Proof}
The residue spectral sequences associated to the lines
\[
\vcenter{\xymatrix@1{ 
G_{i} &
F_{i} \ar@{~>}[l] &
F_{i-1} \ar@{~>}[l] &
G_{i} \ar@{~>}[l] \\}}
\]
organize into an exact couple satisfying the conditions of \ref{13aF}. 
\end{Proof}

In the situation of \ref{13bD}, we refer to 
\[
E_1^{p,q} \Longrightarrow 
H_?^{p+q}(F_0,j) 
\]
as the \emph{spectral sequence associated to the
pure filtered system $\FF$}.

\begin{Ex} \label{13bE}
In the case of a real pure filtered system, we
recover Theorem~\ref{1F}.
\end{Ex}

Recall the notation of Subsection~\ref{9a}. 
Let $U \in \Ch_{\pi , P}$, and consider the following real filtered system 
$\FF(U)$ on $\Eh^k$:
\[
F_p \Eh^k := \{ (x_1 , x_2 , y_1 , \ldots , y_{k-2}) \in \Eh^k \tei 
\mbox{at most $2+p$ 
of the $x_i$ lie in} \; U_{\infty} \}
\]
if $p \ge -2$, 
\[
F_p \Eh^k := U^2 \times_B  
\{ (y_1 , \ldots , y_{k-2}) \in \Eh^{k-2} \tei 
\mbox{at most $k+p$ 
of the $y_i$ lie in} \; U_{\infty} \} 
\]
if $p \le -2$. 
Defining $G_p := F_p - F_{p-1}$, we have:
\begin{eqnarray*}
G_0 &=& U_\infty \times_B U_\infty \times_B \Eh^{k-2} \; , \\
G_{-1} &=& U_\infty \times_B U \times_B \Eh^{k-2} \coprod
U \times_B U_\infty \times_B \Eh^{k-2} \; , \\
G_{p} &=& \coprod U^2 \times_B U_\infty^{k+p} \times_B U^{-p-2} 
\end{eqnarray*}
if $p \le -2$, where the disjoint union runs through all subsets
of $\{ 1,\ldots,k-2 \}$ of cardinality $k+p$. As for the codimensions
$c_p$ of $G_p$, we have:
\[
c_p = 2+p \quad \text{if} \quad p \ge -1 \; ,
\]
\[
c_p = k+p \quad \text{if} \quad p \le -2 \; .
\]
We get the version ``before $\widetilde{\quad}$, $\ls{U \in \Ch_{\pi , P , -}}$,
and ${ ^{+, \, \sgn}}\;$''
of Theorem~\ref{9aC}:

\begin{Prop} \label{13bF}
Let $k \ge 2$, and $U \in \Ch_{\pi , -}$. 
\begin{enumerate}
\item[(a)]
There is a 
natural spectral sequence
\[
^k{\stern} \quad \quad \quad \quad 
^kE^{p,q}_1 = {^kE^{p,q}_1(U)} \Longrightarrow 
H^{2k+p+q}_? (\Eh^2 \times_B \Eh^{k-2}, k)_{- , \ldots , -} \; ,
\]
where $^kE^{p,q}_1 = 0$ if $p \le -k-1$. If $-k \le p \le -2$, then
$$^kE^{p,q}_1 = 
\bigoplus 
H^{-p+q}_? (U^2 \times_B U^I , -p)_{-,\ldots,-} \otimes_F 
F [U_\infty(B)]_-^{\otimes (k+p)} \; ,$$
where the sum $\oplus$
runs through all subsets $I$ of $\{ 3, \ldots , k \}$ of cardinality $-p-2$, and
the subscript $-, \ldots ,-$ refers to the action of 
multiplication by $-1$
on all $-p$ components of $U^2 \times_B U^I$.
If $p \ge -1$, then
$$^kE^{p,q}_1 = \bigoplus
F [U_\infty(B)]_-^{\otimes (2+p)} \otimes_F 
H^{2(k-2)-p+q}_? (U^J \times_B \Eh^{k-2} , k-2-p)_{-,\ldots,-} \; ,$$
where the sum $\oplus$ runs through all subsets $J$ 
of $\{ 1 , 2 \}$ of cardinality $-p$, and
the subscript $-, \ldots ,-$ refers to the action of 
multiplication by $-1$
on all components of $U^J \times_B \Eh^{k-2}$.
\item[(b)]
The differential 
$^2\partial^{-2,q}_1$ on $^2E^{-2,q}_1$
is given by the residue (in the two
possible directions) between
\[
H^{2+q}_?(U^2 ,2)_{-,-}
\]
and 
\[
\bigl( F [U_\infty(B)]_- \otimes_F 
H^{1+q}_?(U , 1)_- \bigr)^{\oplus 2} \; .
\]
\item[(c)]
For $k \ge 3$, the differential 
$^k\partial^{-k,q}_1$ on $^kE^{-k,q}_1$
is given by the residue (in the last 
$k-2$ directions) between
\[
H^{k+q}_?(U^2 \times_B U^{k-2},k)_{-, \ldots ,-} 
\]
and 
\[
\bigl( H^{k-1+q}_?(U^2 \times_B U^{k-3},k-1)_{-, \ldots ,-} \otimes_F 
F [U_\infty(B)]_- \bigr)^{\oplus (k-2)} \; .
\]
\item[(d)]
For $k \ge 2$, the edge morphism
\[
^k\varrho^q: H^{k+q}_? (\Eh^2 \times_B \Eh^{k-2}, k)_{- , \ldots , -} \pfeil
H^{k+q}_?(U^2 \times_B U^{k-2},k)_{-, \ldots ,-}
\]
is given by restriction.
\item[(e)] The spectral sequence $^k{\stern}$ is compatible with
the regulators. It is 
contravariantly functorial with respect to restriction from $U$ to
smaller objects $U'$ of $\Ch_{\pi , -}$,
and with respect to change of the base $B$.
\end{enumerate}
\end{Prop}

\begin{Proof}
This is the
spectral sequence \ref{1F} associated to the filtration $\FF(U)$.
\end{Proof} 

For fixed $q$, let us define the quotient complex $^k\tE^{\bullet,q}_1(U)$
of $^kE^{\bullet,q}_1(U)$ by modifying \ref{13bF}~(a) as follows: We replace
$H^i_?(U^r \times_B \Eh^s,j)_{-,\ldots,-}$ by 
\[
\tH^i_?(U^r \times_B \Eh^s,j)_{-,\ldots,-}
\] 
(see Definition~\ref{9aB}),
and $F [U_\infty(B)]_-$ by its image $[U_\infty(B)_F]$ 
under the Abel--Jacobi map $[ \argdot ]$. Thus:\\[0.2cm]
$^k\tE^{p,q}_1(U) = 0$ if $p \le -k-1$. If $-k \le p \le -2$, then
$$^k\tE^{p,q}_1(U) = 
\bigoplus 
\tH^{-p+q}_? (U^2 \times_B U^I , -p)_{-,\ldots,-} \otimes_F 
[U_\infty(B)_F]^{\otimes (k+p)} \; ,$$
where the sum $\oplus$
runs through all subsets $I$ of $\{ 3, \ldots , k \}$ of cardinality $-p-2$.
If $p \ge -1$, then
$$^k\tE^{p,q}_1(U) = \bigoplus
[U_\infty(B)_F]^{\otimes (2+p)} \otimes_F 
\tH^{2(k-2)-p+q}_? (U^J \times_B \Eh^{k-2} , k-2-p)_{-,\ldots,-} \; ,$$
where the sum $\oplus$ runs through all subsets $J$ 
of $\{ 1 , 2 \}$ of cardinality $-p$.\\

The main result of this subsection reads as follows:

\begin{Thm} \label{13bG}
The quotient complexes $^k\tE^{\bullet,q}_1(U)$ of $^kE^{\bullet,q}_1(U)$
organize into a quotient spectral sequence of \ref{13bF} with the same end
terms:
\[
^k\tilde{\stern} \quad \quad \quad \quad 
^k\tE^{p,q}_1 = {^k\tE^{p,q}_1(U)} \Longrightarrow 
H^{2k+p+q}_? (\Eh^2 \times_B \Eh^{k-2}, k)_{- , \ldots , -} \; .
\]
\end{Thm}

\begin{Proofof}{Theorem~\ref{9aC}}
Apply the operations ``$\ls{U \in \Ch_{\pi , P , -}}$''
and ``$\;{ ^{+, \, \sgn}}\;$'' to the above result.
\end{Proofof} 
\quad \\

\begin{Proofof}{Theorem~\ref{13bG}}
Denote by
$^kF^{\bullet,q}_1$ the kernel of the projection
\[
^kE^{\bullet,q}_1 \longonto {^k\tE^{\bullet,q}_1} \; ,
\]
and recall the two criteria \ref{13aA} and \ref{13aG} allowing us to
divide out subcomplexes of complexes of $E_1$-terms. 
Our analysis of the complex $^kF^{\bullet,q}_1$ will 
consist of three steps, corresponding to different types of
coordinate directions:
(i) the first coordinate direction, (ii) the second
coordinate direction, (iii) the remaining $k-2$ coordinate directions.
We shall use \ref{13aG} for (i) and (ii), and \ref{13aA} for (iii).
The point of using unreal filtered systems is that they produce
spectral sequences that converge to $0$.\\

(i) Denote by $\bH^i_?(U^r \times_B \Eh^s,j)$ the quotient
of $H^i_?(U^r \times_B \Eh^s,j)$ by the image of the cup product
\[
H^1_?(U,1) \otimes_F H^{i-1}_?(U^{r-1} \times_B \Eh^s,j-1) \pfeil 
H^i_?(U^r \times_B \Eh^s,j)
\]
in the \emph{first} coordinate direction. Note that because of the
existence of the residue, this cup product is injective. 
For fixed $q$, define the quotient complex $^k\bE^{\bullet,q}_1$ 
as follows: \\[0.2cm]
$^k\bE^{p,q}_1 = 0$ if $p \le -k-1$ or $p \ge 1$. If $-k \le p \le -2$, then
$$^k\bE^{p,q}_1 = 
\bigoplus 
\bH^{-p+q}_? (U^2 \times_B U^I , -p)_{-,\ldots,-} \otimes_F 
F [U_\infty(B)]_-^{\otimes (k+p)} \; ,$$
where the sum $\oplus$
runs through all subsets $I$ of $\{ 3, \ldots , k \}$ of cardinality $-p-2$.
Furthermore,
$$^k\bE^{-1,q}_1 = 
\begin{array}{c}
[U_{\infty} (B)_F] \otimes_F 
               H^{2k-3+q}_? (U \times_B \Eh^{k-2} , k-1)_{-,\ldots,-} \\
\oplus \\
\bH^{2k-3+q}_? (U \times_B \Eh^{k-2} , k-1)_{-,\ldots,-} 
                                       \otimes_F F [U_\infty (B)]_-  
\end{array}
$$
and
$$^k\bE^{0,q}_1 = 
[U_\infty(B)_F] \otimes_F F [U_\infty(B)]_- \otimes_F
H^{2k-4+q}_? (\Eh^{k-2} , k-2)_{-,\ldots,-} \; .$$
We proceed to prove the following claim, which we refer to by the sign $(+)$:
The quotient complexes $^k\bE^{\bullet,q}_1$ of 
$^kE^{\bullet,q}_1$
organize into a quotient spectral sequence of \ref{13bF} with the same end
terms:
\[
^k\bE^{p,q}_1 \Longrightarrow 
H^{2k+p+q}_? (\Eh^2 \times_B \Eh^{k-2}, k)_{- , \ldots , -} \; .
\]
In order to see this, let us study the kernel $\bF^{\bullet,\bullet}_1$
of the projection
\[
^kE^{\bullet,\bullet}_1 \longonto {^k\bE^{\bullet,\bullet}_1} \; :
\]
The terms $p= -1$ and $p=0$ of $\bF^{\bullet,q}_1$ look as follows:
$$\bF^{-1,q}_1 = 
\begin{array}{c}
H^1_?(U,1)_- \otimes_F 
               H^{2k-3+q}_? (U \times_B \Eh^{k-2} , k-1)_{-,\ldots,-} \\
\oplus \\
H^1_?(U,1)_- \otimes_F 
F [U_\infty (B)]_- \otimes_F  H^{2k-4+q}_? (\Eh^{k-2} , k-2)_{-,\ldots,-}                          
\end{array}
$$
and
$$\bF^{0,q}_1 = 
H^1_?(U,1)_- \otimes_F F [U_\infty(B)]_- \otimes_F
H^{2k-4+q}_? (\Eh^{k-2} , k-2)_{-,\ldots,-} \; .$$
Thus, $\bF^{\bullet,q}_1$ contains an obvious acyclic subcomplex 
$G^{\bullet,q}$ concentrated
in degrees $-1$ and $0$, with entries
\[
H^1_?(U,1)_- \otimes_F F [U_\infty(B)]_- \otimes_F
H^{2k-4+q}_? (\Eh^{k-2} , k-2)_{-,\ldots,-} \; .
\]
Let us replace $\bF^{\bullet,\bullet}_1$ by the
quotient of $\bF^{\bullet,\bullet}_1$ by $G^{\bullet,\bullet}\;$; we continue
to refer to this complex as $\bF^{\bullet,\bullet}_1$. The claim $(+)$ will
be a consequence of the following sub-claims:
\begin{enumerate}
\item[$(+a)$] There is a spectral sequence
\[
\bF^{p,q}_1 \Longrightarrow \bF^{p+q} = 0 \; .
\]
\item[$(+b)$] The inclusion of $\bF^{\bullet,\bullet}$ into the quotient
of $^kE^{\bullet,\bullet}_1$ by $G^{\bullet,\bullet}$ extends to a 
morphism 
\[
\vcenter{\xymatrix@R-10pt{ 
\left( \bF^{p,q}_1 \Longrightarrow 0 \right) \ar[d]^\iota \\
\left( {^kE^{p,q}_1} / G^{p,q} \Longrightarrow 
H^{2k+p+q}_? (\Eh^2 \times_B \Eh^{k-2}, k)_{- , \ldots , -} \right) \\}}
\]
of spectral sequences satisfying the hypothesis of \ref{13aG}.
\end{enumerate}
Observe that $\bF^{p,q}_1$ equals $H^1_?(U,1)_- \otimes_F$ the group
\[
A^{-1,q}_1 := H^{2k-3+q}_? (U \times_B \Eh^{k-2} , k-1)_{-,\ldots,-}
\]
if $p=-1$, and $H^1_?(U,1)_- \otimes_F$ the group
$$A^{p,q}_1 := \bigoplus 
H^{-p+q-1}_? (U \times_B U^I , -p-1)_{-,\ldots,-} \otimes_F 
F [U_\infty(B)]_-^{\otimes (k+p)}$$
if $-k \le p \le -2$. 
Now observe that $A^{\bullet,\bullet}$ is the $E_1$-term of the spectral
sequence \ref{13bD} associated to the following pure filtered system $\FF'$
of length $k$ of subschemes of $U \times_B \Eh^{k-2}$:
\[
\vcenter{\xymatrix@R-10pt{ 
G_0'=\emptyset &
F_0'=\emptyset \ar@{~>}[l] &
F_{-1}'=\emptyset \ar@{~>}[l] \ar@{=}[dl] &
G_0'=\emptyset \ar@{~>}[l] \\
G_{-1}'=F_{-2}' &
F_{-1}'=\emptyset \ar@{~>}[l] &
F_{-2}' \ar@{~>}[l] \ar@{=}[dl] &
G_{-1}'=F_{-2}' \ar@{~>}[l] \\
G_{-2}' \ar@{.}[dd] &
F_{-2}' \ar@{~>}[l] \ar@{.}[dd] &
F_{-3}' \ar@{~>}[l] \ar@{.}[dd] &
G_{-2}' \ar@{~>}[l] \ar@{.}[dd] \\
&
&
 \\
G_{-k}' &
F_{-k}' \ar@{~>}[l] &
F_{-k-1}' \ar@{~>}[l] &
G_{-k}' \ar@{~>}[l] \\}}
\]
where
\[
F_p' := U \times_B  
\{ (y_1 , \ldots , y_{k-2}) \in \Eh^{k-2} \tei 
\mbox{at most $k+p$ 
of the $y_i$ lie in} \; U_{\infty} \} 
\]
if $-k \le p \le -2$.
This proves $(+a)$. For $(+b)$, we observe that
the morphism of exact couples required in \ref{13aG}
is given by the cup product in the first coordinate direction. A
splitting can be obtained by \emph{choosing} a left inverse 
(in the category of $F$-vector spaces) of the residue
\[
H^1_?(U,1)_- \longinto F [U_\infty (B)]_- \; .
\]

(ii) We change notation, and 
denote by $\bH^i_?(U^r \times_B \Eh^s,j)$ the quotient
of $H^i_?(U^r \times_B \Eh^s,j)$ by the image of the cup products
\[
H^1_?(U,1) \otimes_F H^{i-1}_?(U^{r-1} \times_B \Eh^s,j-1) \pfeil 
H^i_?(U^r \times_B \Eh^s,j)
\]
in the \emph{first and second} coordinate directions.  
For fixed $q$, define the quotient complex $^k\bE^{\bullet,q}_1$ 
as follows: \\[0.2cm]
$^k\bE^{p,q}_1 = 0$ if $p \le -k-1$ or $p \ge 1$. If $-k \le p \le -2$, then
$$^k\bE^{p,q}_1 = 
\bigoplus 
\bH^{-p+q}_? (U^2 \times_B U^I , -p)_{-,\ldots,-} \otimes_F 
F [U_\infty(B)]_-^{\otimes (k+p)} \; ,$$
where the sum $\oplus$
runs through all subsets $I$ of $\{ 3, \ldots , k \}$ of cardinality $-p-2$.
Furthermore,
$$^k\bE^{-1,q}_1 = 
\begin{array}{c}
[U_{\infty} (B)_F] \otimes_F 
               \bH^{2k-3+q}_? (U \times_B \Eh^{k-2} , k-1)_{-,\ldots,-} \\
\oplus \\
\bH^{2k-3+q}_? (U \times_B \Eh^{k-2} , k-1)_{-,\ldots,-} 
                                       \otimes_F [U_\infty (B)_F] 
\end{array}
$$
and
$$^k\bE^{0,q}_1 = 
[U_\infty(B)_F]^{\otimes 2} \otimes_F
H^{2k-4+q}_? (\Eh^{k-2} , k-2)_{-,\ldots,-} \; .$$
Then the quotient complexes $^k\bE^{\bullet,q}_1$ of 
$^kE^{\bullet,q}_1$
organize into a quotient spectral sequence of \ref{13bF} with the same end
terms:
\[
^k\bE^{p,q}_1 \Longrightarrow 
H^{2k+p+q}_? (\Eh^2 \times_B \Eh^{k-2}, k)_{- , \ldots , -} \; .
\]
The proof of this claim is formally identical as the one of $(+)$
(except that we modify the spectral sequence obtained there instead
of the spectral sequence $^k{\stern}$).
We leave the details to the reader.\\

(iii) Observe that the terms $^k\bE^{p,q}_1$ are already of the right 
shape for $p \ge -1$. Fix $q$, and define the following data:
\begin{enumerate}
\item[(1)] $A:= F [U_\infty(B)]_-$, and $\Fa:= H^1_?(U,1)_-$.
\item[(2)] for any subset $I \subset \{ 3,\ldots,k \}$, let
\[
B^I := \bH^{-p+q}_? (U^2 \times_B U^I , -p)_{-,\ldots,-} \; .
\]
\item[(3)] for any subset $I \subset \{ 3,\ldots,k \}$, and any $m \in I$,
define
\[
\partial_m^I: B^I \pfeil B^{I - \{ m \} }\otimes_F A 
\]
as the residue in the $m$-th coordinate direction,
and 
\[
\epsilon_m^I: B^{I - \{ m \} }\otimes_F {\Fa} \pfeil B^I 
\]
as the cup product.
\end{enumerate}
Now apply \ref{13aB} and \ref{13aA}.
\end{Proofof}

\begin{Rem} \label{13H}
The same method works for other twists. The result 
(which will not be needed in this generality) reads as follows:\\

Let $j \in \BZ$. There is a 
natural spectral sequence
$$E^{p,q}_1 \Longrightarrow 
H^{2k+p+q}_? (\Eh^2 \times_B \Eh^{k-2} , j)_{-,\ldots,-} \; ,$$
where $E^{p,q}_1 = 0$ if $p \le -k-1$.
If $-k \le p \le -2$, then
$$E^{p,q}_1 = 
\ls{U \in \Ch_{\pi , P , -}} \bigoplus 
\tH^{-p+q}_? (U^2 \times_B U^I , j-k-p)_{-,\ldots,-} \otimes_F 
[P_F]^{\otimes (k+p)} \; ,$$
where the sum $\oplus$
runs through all subsets $I$ of $\{ 3, \ldots , k \}$ of cardinality $-p-2$.
If $p \ge -1$, then
$$^k\tE^{p,q}_1 = \bigoplus
[P_F]^{\otimes (2+p)} \otimes_F \ls{U \in \Ch_{\pi , P , -}} 
\tH^{2(k-2)-p+q}_? (U^J \times_B \Eh^{k-2} , j-2-p)_{-,\ldots,-} \; ,$$
where the sum $\oplus$ runs through all subsets $J$ 
of $\{ 1 , 2 \}$ of cardinality $-p$.
\end{Rem}

We conclude the subsection with the following result, which will be needed
in the proof of Main Theorem~\ref{9bI}:

\begin{Prop} \label{13bI}
Let $k \ge 2$. There is a 
natural spectral sequence
$$E^{p,q}_1 \Longrightarrow 
\ls{U \in \Ch_{\pi , P , -}}
        \tH^{2k-2+p+q}_? (U \times_B \Eh^{k-2} , k-1)^{\sgn}_{-,\ldots,-} \; ,$$
where 
$$E^{p,q}_1 = 
\ls{U \in \Ch_{\pi , P , -}} 
\tH^{-p+q}_? (U \times_B U^{-p-1} , -p)^{\sgn}_{-,\ldots,-} \otimes_F 
\bigwedge^{k-1+p} [P_F] \; .$$
Here, the superscript $\sgn$ refers to the action of $\fS_{k-2}$ and
$\fS_{-p-1}$ respectively, while the subscript $-,\ldots,-$
refers to the action on all components.
\end{Prop}

\begin{Proof}
We only sketch the construction of the spectral sequence, because we
have already seen all the necessary techniques. 
\begin{enumerate} 
\item[(1)] Use a suitable filtration $\FF$ to get
$$E^{p,q}_1 \Longrightarrow 
\ls{U \in \Ch_{\pi , P , -}}
        H^{2k-2+p+q}_? (U \times_B \Eh^{k-2} , k-1)^{\sgn}_{-,\ldots,-} \; ,$$
where 
$$E^{p,q}_1 = 
\ls{U \in \Ch_{\pi , P , -}} 
H^{-p+q}_? (U \times_B U^{-p-1} , -p)^{\sgn}_{-,\ldots,-} \otimes_F 
\bigwedge^{k-1+p} F [U_\infty(B)]_- \; .$$ 
\item[(2)] Using \ref{13aB} and \ref{13aA}, modify (1) to get
$$E^{p,q}_1 \Longrightarrow 
\ls{U \in \Ch_{\pi , P , -}}
        H^{2k-2+p+q}_? (U \times_B \Eh^{k-2} , k-1)^{\sgn}_{-,\ldots,-} \; ,$$
where 
$$E^{p,q}_1 = 
\ls{U \in \Ch_{\pi , P , -}} 
\bH^{-p+q}_? (U \times_B U^{-p-1} , -p)^{\sgn}_{-,\ldots,-} \otimes_F 
\bigwedge^{k-1+p} [P_F] \; ,$$ 
and $\bH^i_? (U \times_B U^s , j)$ denotes the quotient of 
$H^i_? (U \times_B U^s , j)$ by the images of the cup products
\[
H^1_?(U,1) \otimes_F H^{i-1}_? (U \times_B U^{s-1} , j-1)
\pfeil H^i_? (U \times_B U^s , j)
\]
in all but the first coordinate direction.
\item[(3)] Imitate the above to obtain 
$$F^{p,q}_1 \Longrightarrow 
\ls{U \in \Ch_{\pi , P , -}} H^1_? (U,1)_- \otimes_F
        H^{2k-3+p+q}_? (\Eh^{k-2} , k-2)^{\sgn}_{-,\ldots,-} \; ,$$
where 
$$F^{p,q}_1 = 
\ls{U \in \Ch_{\pi , P , -}} H^1_? (U,1)_- \otimes_F
\bH^{-p+q-1}_? (U^{-p-1} , -p-1)^{\sgn}_{-,\ldots,-} \otimes_F 
\bigwedge^{k-1+p} [P_F] \; .$$ 
\item[(4)] Use \ref{13aG} to divide; again, the required splitting is obtained
by the choice of a left inverse of the residue
\[
H^1_?(U,1)_- \longinto F [U_\infty (B)]_- \; .
\]
\end{enumerate}
\end{Proof}

\newpage
\subsection{Motives generated by elliptic curves. II} \label{14}

Let $\pi:\Eh \pfeil B$ be an elliptic curve over a regular, noetherian base
$B$.
The aim of this subsection is to give a proof of the following
result, which will be used frequently
in our analysis of the spectral sequence \ref{9aC}:

\begin{Thm} \label{14A}
The exterior cup product
\[
\Eh(B) \otimes_\BZ H^\bullet_? (\Eh^{n-1},\argstern)_{-,\ldots,-}
\pfeil
H^{\bullet +2}_? (\Eh^{n} , \argstern +1)_{-,\ldots,-} 
\]
lands in the $\sgn$-eigenpart of 
$H^{\bullet +2}_? (\Eh^{n} , \argstern +1)_{-,\ldots,-}$ under the action
of $\fS_n$. 
\end{Thm}

In order to prepare the proof of \ref{14A}, 
recall the decomposition
\[
h^\bullet(\Eh) = h^0(\Eh) \oplus h^1(\Eh) \oplus h^2(\Eh) =
\BQ(0) \oplus h^1(\Eh) \oplus \BQ(-1) 
\]
of Subsection~\ref{6}.
The direct summand $h^0(\Eh)$ can be characterized as the image of 
\[
\pi^*: h^\bullet (B) \pfeil h^\bullet (\Eh) \; ,
\]
or, equivalently, of the projector $\pi^* i^*$ of $h^\bullet (\Eh)$.
Similarly, the direct summand $h^0(\Eh) \oplus h^1(\Eh)$ is the kernel of
\[
\pi_*: h^\bullet (\Eh) \pfeil h^{\bullet-2} (B) \otimes \BQ(-1) \; ,
\]
or, equivalently, of the projector $i_* \pi_*$.\\

By \ref{6B}, the decomposition corresponds 
to the decomposition into eigenparts
for the values $1$, $a$, and $a^2$ of the action of
$\tr_{[a]^1}$, for any integer $a$, whose absolute value is at least $2$.
Furthermore, $[-1]$ acts by $-1$ on $h^1(\Eh)$, and as the identity on
$h^0(\Eh)$ and $h^2(\Eh)$ (\ref{6C}). There are analogous
decompositions of $H^\bullet_? (\Eh, \argstern )_-$.

\begin{Thm} \label{14Ba}
Let $s \in \Eh(B)$ be a section.
Then translation by $s$ respects the eigenspace
\[
H^\bullet_? (\Eh, \argstern )_- \; ,
\]
and it acts trivially.
\end{Thm}

\begin{Proof}
As in Subsection~\ref{7}, write
\[
H^\bullet_? (\Eh, \argstern ) = 
                \bigoplus_{r=0}^2 H^\bullet_? (\Eh, \argstern )_{(r)}
\]
for the above decomposition, and abbreviate
\[
H_{(r)} := H^\bullet_? (\Eh, \argstern )_{(r)} \; .
\]
Thus $H_{(1)} = H^\bullet_? (\Eh, \argstern )_-$.
Denote the translation by $+_s$. The direct summand
\[
H_{(0)} \oplus H_{(1)} \subset H^\bullet_? (\Eh, \argstern )
\]
can be characterized as the kernel of the morphism
\[
\pi_*: H^\bullet_? (\Eh, \argstern ) \pfeil
H^{\bullet -2}_? (B, \argstern -1 ) \; .
\]
Since $\pi \circ +_s = \pi$, we see that $H_{(0)} \oplus H_{(1)}$
is respected by $(+_s)_*$. Next, we note that the
direct summand
\[
H_{(0)} \subset H^\bullet_? (\Eh, \argstern )
\]
can be characterized as the image of the morphism
\[
\pi^*: H^\bullet_? (B, \argstern ) \pfeil 
H^\bullet_? (\Eh, \argstern )  \; .
\]
Since $\pi \circ +_{-s} = \pi$, this image is respected by $+_{-s}^* = (+_s)_*$.
Furthermore, the action of $+_s$ on $H_{(0)}$ is trivial.
Let us write
\[
(+_s)_* = ((+_s)_*^0 , (+_s)_*^1) : H_{(1)} \pfeil H_{(0)} \oplus H_{(1)} \;.
\]
By comparing coefficients in 
\[
(+_s)_* \circ (+_s)_* ([2]_* \alpha) = (+_{2s})_* ([2]_* \alpha) = 
[2]_* (+_s)_* (\alpha) \; ,
\]
for $\alpha \in H_{(1)}$, we obtain the identities
\[
(+_s)_*^1 \circ (+_s)_*^1 = (+_s)_*^1 
\]
and
\[
(+_s)_*^0 + 2 (+_s)_*^0 \circ (+_s)_*^1 = 0 \; .
\]
But $(+_s)_*^1$ is invertible (the inverse being $(+_{-s})_*^1$), hence
\[
(+_s)_*^1 = \id \quad \mbox{and} \quad (+_s)_*^0 = 0 \; .
\]
\end{Proof}

\begin{Rem} \label{14C}
(a) We expect a corresponding statement on the level of
motives: The direct summand $h^1(\Eh) \subset h^\bullet(\Eh)$ should
be respected by translations by sections, and they should act trivially.
If the base is a field, this is easy to show (write down the corresponding
equation in the Chow ring of $\Eh^2$). In the relative case, note that
translations do not define morphisms in $\Mh^{(\Eh)}$.\\[0.2cm]
(b) If $s$ is non-torsion, then translation by $s$ does \emph{not} act
trivially on motivic cohomology
$H^2_\Mh (\Eh,1)$.\\[0.2cm]
(c) If $s$ is torsion, then translation by $s$ \emph{does} act
trivially on $H^\bullet_? (\Eh,\argstern)$: see \cite{B5}, Lemma~3.2.1, whose proof
works also for families of elliptic curves.
\end{Rem}

\begin{Cor} \label{14D}
The exterior cup product
\[
\Eh(B) \otimes_\BZ H^\bullet_? (\Eh,\argstern)_{-}
\pfeil
H^{\bullet +2}_? (\Eh^{2} , \argstern +1)_{-,-} 
\]
lands in 
$H^{\bullet +2}_? (\Eh^{2} , \argstern +1)_{-,-}^\sgn$. 
\end{Cor}

\begin{Proof}
By Proposition~\ref{9aZ}, the projection
\[
H^{\bullet +2}_? (\Eh^{2} , \argstern +1)_{-,-} \pfeil
H^{\bullet +2}_? (\Eh^{2} , \argstern +1)_{-,-}^+
\]
has the same kernel as the map
\[
H^{\bullet +2}_? (\Eh^{2} , \argstern +1)_{-,-} \pfeil
H^{\bullet}_? (B , \argstern)
\]
given by the composition of $\Sigma_*$ and $i^*$. It thus suffices to show that
the composition $\gamma$
of the exterior cup product and $\Sigma_*$ is the trivial map
\[
\Eh(B) \otimes_\BZ H^\bullet_? (\Eh,\argstern)_{-} \pfeil
H^\bullet_? (\Eh,\argstern) \; .
\]
Let $s \in \Eh(B)$, and $x \in H^\bullet_? (\Eh,\argstern)_{-}$. Then
we have
\[
\gamma(s \otimes x) = (+_s)_* (x) - x \; .
\]
Our claim follows from \ref{14Ba}.
\end{Proof}
\quad \\

\begin{Proofof}{Theorem~\ref{14A}}
Apply \ref{14D} in all $n-1$ coordinate directions. Then use the fact that
$\fS_n$ is generated by the transpositions
\[
(1,2) \; , \; (1,3) \; , \ldots , \; (1,n) \; .
\]
\end{Proofof}

\newpage
\subsection{The case $k=3$} \label{11}

The purpose of this subsection is to give proofs of 
Theorem~\ref{9bA} and Main Theorem~\ref{9bI} for $k=3$.
As in Subsection~\ref{10}, $P \subset \Eh(B)$ 
satisfies $(DP)$ (see \ref{9aA}), and
$? \in \{ \Mh , \abs \}$.\\

Let $U \in \Ch_{\pi , P}$. The filtration $\FF(U)$ on $\Eh^3$ introduced
in Subsection~\ref{13b} has the following shape:
\[
F_p \Eh^3 := \{ (x_1 , x_2 , y) \in \Eh^3 \tei \mbox{at most $2+p$ 
of the $x_i$ lie in} \; U_{\infty} \} 
\quad \text{if} \quad p \ge -2 \; ,
\]
\[
F_{-3} \Eh^3 := U^3 \; ,
\]
and $F_{-4} := \emptyset$. 
We put
\[
G_p := F_p - F_{p-1} \; .
\]
Thus we have:
\begin{eqnarray*}
G_0 &=& U_\infty \times_B U_\infty \times_B \Eh \; , \\
G_{-1} &=& U_\infty \times_B U \times_B \Eh \coprod
U \times_B U_\infty \times_B \Eh \; , \\
G_{-2} &=& U^2 \times_B U_\infty \; , \\
G_{-3} &=& U^3 \; .
\end{eqnarray*}
For the codimensions $c_p$ of $G_p$, we have:
\[
c_0 = 2 \; , \; c_{-1} = 1 \; , \; 
c_{-2} = 1 \; , \; c_{-3} = 0 \; .
\]
From Theorem~\ref{13bG}, we get:

\begin{Prop} \label{11C}
There is a 
natural spectral sequence
$$^3\tilde{\stern} \quad \quad \quad \quad ^3\tE^{p,q}_1 \Longrightarrow 
H^{6+p+q}_? (\Eh^3 , 3)_{-,-,-} \; ,$$
where $^3\tE^{p,q}_1 = 0$ if $p \le -4$,
$$^3\tE^{p,q}_1 = 
\ls{U \in \Ch_{\pi , P , -}} 
\tH^{-p+q}_? (U^2 \times_B U^{-p-2} , -p)_{-,\ldots,-} \otimes_F 
[P_F]^{\otimes (3+p)}$$
if $-3 \le p \le -2$, and 
$$^3\tE^{p,q}_1 = \bigoplus
[P_F]^{\otimes (2+p)} \otimes_F 
\ls{U \in \Ch_{\pi , P , -}} 
\tH^{2-p+q}_? (U^{-p} \times_B \Eh , 1-p)_{-,\ldots,-} $$
if $p \ge -1$,
where the sum $\oplus$ runs through all subsets of $\{ 1 , 2 \}$ of cardinality
$-p$.
\end{Prop}

Recall that due to our definition of the filtrations $\FF(U)$,
the full symmetric group $\fS_3$ does \emph{not} act on our complexes.
The situation is however equivariant under $\fS_2$,
which acts on the \emph{first} two of the three coordinates. We obtain
the spectral sequence of
Theorem~\ref{9aC} for $k=3$ by passing
to the $\fS_2$-invariants of the spectral sequence
$^3\tilde{\stern}$, and applying the isomorphism
\[
\jmath: H^{6+p+q}_? (\Eh^2 \times_B \Eh , 3)_{-,-,-}^+ \isoto 
                          H^{4+p+q}_? (\Eh^{(1)} , 2)^\sgn
\]
of \ref{9aZ}.
We introduce the following abbreviations: 
\begin{eqnarray*}
\tH^{q} (U^3) &:=& \tH^{q}_? (U^3 , 3)_{-,-,-} \; , \\
\tH^{q} (U^2) &:=& \tH^{q}_? (U^2 , 2)_{-,-} \; , \\
\tH^{q} (U \times_B \Eh) &:=& \tH^{q}_? (U \times_B \Eh , 2)_{-,-} \; .
\end{eqnarray*}
Let us study the lines $q=1$ and $q=0$ of $^3\tilde{\stern}$. 
The columns $-3,-2,-1$ have the following shape:
\[
\vcenter{\xymatrix@R-10pt{ 
\ld \tH^4 (U^3)
         \ar[r]^-{\partial_1} \ar@{-->}[drr]^{\partial_2} & 
\ld \tH^3 (U^2) \otimes_F [P_F] \ar[r] &
\quad\quad\quad\quad\quad\quad\quad\quad\quad\quad\quad\quad \\
\quad\quad\quad\quad\quad \ar[r] &
\ld \tH^2 (U^2) \otimes_F [P_F] \ar[r] &
[P_F] \otimes_F \ld \tH^3 (U \times_B \Eh)^{\oplus 2} \\}}
\]
The differential $\partial_2$ is only
defined on the kernel of $\partial_1$.\\

Let us consider the situation after passing to the $\fS_2$-invariants,
and together with the Bloch groups:
\[
\vcenter{\xymatrix@R-10pt{ 
\Bl_{3,P,?} \ar[r]^-{d_3} \ar@{_{(}->}[d] & 
\Bl_{2,P,?} \otimes_F [P_F] \ar@{_{(}->}[d] & \\
\ld \tH^4 (U^3)^+
         \ar[r]^-{\partial_1} \ar@{-->}[drr]^{\partial_2} & 
\ld \tH^3 (U^2)^+ \otimes_F [P_F] \ar[r] &
\quad\quad\quad\quad\quad\quad\quad\quad\quad\quad\quad \\
\quad\quad\quad\quad\quad\;\;\; \ar[r] &
0 \ar[r] &
[P_F] \otimes_F \ld \tH^3 (U \times_B \Eh) \\}}
\]
Here, we have used \ref{10K} to see that $^3\tE^{-2,0,+}_1$ is trivial.
By definition of the differential $d_3$, the diagram is commutative.

\begin{Prop} \label{11E}
Restriction to the $U^3$ induces a canonical isomorphism
\[
H^2_? (\Eh^{(1)} , 2)^\sgn
 \stackrel{\ref{9aZ}}{=}
H^4_? (\Eh^2 \times_B \Eh , 3)_{-,-,-}^+ \isoto \kerr(\partial_2) \; .
\]
\end{Prop}

The proof of \ref{11E} will require the following result:

\begin{Lem} \label{11Ea}
For $i \le 2$, restriction induces an isomorphism
\[
H^i_? (\Eh \times_B \Eh,2)_{-,-} \isoto H^i_? (U \times_B \Eh,2)_{-,-}
          \isoto \tH^i_? (U \times_B \Eh,2)_{-,-} \; .
\]
\end{Lem}

\begin{Proof}
Consider the localization sequence
\begin{eqnarray*}
\ldots \! & \! \stackrel{\res}{\pfeil} \! & \! 
                        F [U_\infty]_- \otimes_F H^{i-2}_?(\Eh,1)_- \\
\pfeil H^i_? (\Eh \times_B \Eh,2)_{-,-} \pfeil H^i_? (U \times_B \Eh,2)_{-,-} 
\! & \! \stackrel{\res}{\pfeil} \! & \! 
                          F [U_\infty]_- \otimes_F H^{i-1}_?(\Eh,1)_- \; .
\end{eqnarray*}
\end{Proof}

\begin{Proofof}{Proposition~\ref{11E}}
The target of the differential $\partial_3$ on $\kerr(\partial_2)$ is
a quotient of
\[
\Sym^2 [P_F] \otimes_F H^1_?(\Eh , 1)_- \; ,
\]
which is trivial. Therefore, we have
\[
\kerr(\partial_2) = {^3\tE^{-3,1,+}_\infty} \; .
\]
We know already that $^3\tE^{-2,0,+}_\infty = 0$.
Since $H^0_?(\Eh , 1)_- = 0$, the term $^3\tE^{0,-2,+}_\infty$
is trivial. We therefore have an exact sequence
\[
0 \pfeil {^3\tE^{-1,-1,+}_\infty} \pfeil 
H^4_? (\Eh^2 \times_B \Eh , 3)_{-,-,-}^+ \pfeil 
\kerr(\partial_2) \pfeil 0 \; .
\]
Let us study the term $^3\tE^{-1,-1,+}_\infty$. It is a subquotient
of
\[
^3\tE^{-1,-1,+}_1 = [P_F] \otimes_F \ld \tH^2_? (U \times_B \Eh , 2)_{-,-} 
\stackrel{\ref{11Ea}}{=}
[P_F] \otimes_F H^2_? (\Eh \times_B \Eh,2)_{-,-} \; .
\]
The edge morphism
\[
^3\tE^{-1,-1,+}_1 \pfeil 
H^4_? (\Eh^2 \times_B \Eh , 3)_{-,-,-}^+
\]
equals the composition of the exterior cup product
\[
[P_F] \otimes_F H^2_? (\Eh \times_B \Eh,2)_{-,-}
\stackrel{\cup}{\pfeil} 
H^4_? (\Eh^2 \times_B \Eh , 3)_{-,-,-} \; ,
\]
and the projection onto the $+$-eigenspace.
By Theorem~\ref{14A},
this composition is 
trivial. Therefore, the term $^3\tE^{-1,-1,+}_\infty$ is trivial.
\end{Proofof}

In order to prove the analogue of Proposition~\ref{10F} for $k=3$
(see \ref{11G}), we need the following result:

\begin{Lem} \label{11F}
Let $U \in \Ch_{\pi , P , -}$. 
Restriction to $U \times_B U$ induces an injection
\[
\tH^3_? (U \times_B \Eh , 2)_{-,-} \pfeil
\tH^3_? (U \times_B U , 2)_{-,-} \; .
\]
\end{Lem}

\begin{Proof}
Consider the localization sequence
\begin{eqnarray*}
\ldots \! & \! \stackrel{\res}{\pfeil} \! & \! 
                        H^1_?(U,1)_- \otimes_F F [U_\infty]_- \\
\pfeil H^3_? (U \times_B \Eh,2)_{-,-} \pfeil H^3_? (U \times_B U,2)_{-,-} 
\! & \! \stackrel{\res}{\pfeil} \! & \! \ldots 
\end{eqnarray*}
Together with surjectivity of
\[
H^2_? (\Eh,1)_- \pfeil H^2_? (U,1)_- \; , 
\]
it shows that restriction induces an
injection
\[
\tH^3_? (U \times_B \Eh,2)_{-,-} \longinto
\bar{H}^3_? (U \times_B U , 2)_{-,-} \; ,
\]
where $\bar{H}^3_? (U \times_B U , 2)_{-,-}$ 
is defined as the quotient of $H^3_? (U \times_B U , 2)_{-,-}$ by the
image of the cup product
\[
\cup: H^1_? (U , 1)_- \otimes_F H^2_? (U , 1)_- \pfeil
H^3_? (U \times_B U , 2)_{-,-} 
\]
(in this single direction). 
Now consider the map
\[
\bar{H}^3_? (U \times_B U , 2)_{-,-} \pfeil
\left( \tH^3_? (U \times_B U , 2)_{-,-} \oplus 
(\tH^2_? (U , 1)_- \otimes_F F [U_\infty (B)]_-) \right)
\]
given by the projection, and the residue $\res_2$
in the second coordinate direction.
It is injective because of
the definition of $\tH^3_? (U \times_B U , 2)_{-,-}$, and because
the residue 
\[
H^1_? (U , 1)_- \pfeil F [U_\infty (B)]_-
\]
is injective. But $\res_2$ is trivial on the image
of $\tH^3_? (U \times_B \Eh,2)_{-,-}$.
\end{Proof}

\begin{Prop} \label{11G}
The sequence
\[
0 \pfeil H^2_? (\Eh^{(1)} , 2)^\sgn \stackrel{\varrho_3}{\pfeil} \Bl_{3,P,?}
\stackrel{d_3}{\pfeil} \Bl_{2,P,?} \otimes_F [P_F]
\]
is exact.
\end{Prop}

\begin{Proof}
By definition of $d_3$, we have 
$\kerr(d_3) = \Bl_{3,P,?} \cap \kerr(\partial_1)$.
In view of \ref{11E},
the only thing that remains to be shown is that the kernel of $d_3$
is contained in $\kerr(\partial_2)$.\\ 

In order to do so, let us study the map
\[ 
\partial_2: {^3\tE^{-3,1,+}_2} \pfeil {^3\tE^{-1,0,+}_2} \; .
\]
Since $^3\tE^{-2,0,+}_1 = 0$, we have
\[
^3\tE^{-1,0,+}_2 \subset {^3\tE^{-1,0,+}_1} = 
[P_F] \otimes_F \ld \tH^3_? (U \times_B \Eh , 2)_{-,-} \; .
\]
By \ref{11F}, the latter group is contained in 
\[
[P_F] \otimes_F \ld \tH^3_? (U \times_B U , 2)_{-,-} \; .
\]
With these identifications, the differential
$\partial_2$ is induced by
the residue map
\[
\Res: \ld \tH^4_? (U^2 \times_B U , 3)_{-,-,-}^+ \pfeil
[P_F] \otimes_F \ld \tH^3_? (U \times_B U , 2)_{-,-} 
\]
(in any of the first two coordinate directions). Next, recall the definition
of $\Bl_{3,P,?}$ (see \ref{9aD}): 
it is the image of the map $p$ from the kernel of
\[
\ld H^4_?(\Eh \times_B U^2 , 3)^\sgn_{-,-,-}
\pfeil
\ld H^4_?(U^3 , 3)^\Sgn_{-,-,-}
\]
($\sgn$ referring to the action of $\fS_2$ on the \emph{last} two
coordinates, and $\Sgn$ to the action of $\fS_3$)
to 
\[
\ld \tH^4_? (U^2 \times_B U , 3)_{-,-,-}^+ \; .
\]
The map $p$ is defined as the following composition:
\[
\vcenter{\xymatrix@R-10pt{ 
\kerr \left(
\ld H^4_?(\Eh \times_B U^2 , 3)^\sgn_{-,-,-}
\pfeil
\ld H^4_?(U^3 , 3)^\Sgn_{-,-,-} 
\right) 
\ar[d] \\
\kerr \left(
\ld \tH^4_?(U \times_B U^2 , 3)^\sgn_{-,-,-}
\pfeil
\ld \tH^4_?(U^3 , 3)^\Sgn_{-,-,-} 
\right) 
\ar@{_{(}->}[d]^{\symm_{1,2}} \\
\ld \tH^4_? (U^2 \times_B U , 3)_{-,-,-}^+ \\}}
\]
Here, $\symm_{1,2}$ denotes the 
symmetrization with respect to the first two coordinates.
It will be vital to observe that it is injective on 
\[
\kerr \left(
\ld \tH^4_?(U \times_B U^2 , 3)^\sgn_{-,-,-}
\pfeil
\ld \tH^4_?(U^3 , 3)^\Sgn_{-,-,-} 
\right) \; .
\]
We have to show that $\partial_2$ is trivial on the intersection of
$\kerr(\partial_1)$ and the image of $p$.\\

The composition $\Res \circ p$ equals 
the symmetrization with respect to the first two coordinates
$\halb (\res_1 + \res_2)$ of the residue map from 
\[
\kerr \left(
\ld H^4_?(\Eh \times_B U^2 , 3)^\sgn_{-,-,-}
\pfeil
\ld H^4_?(U^3 , 3)^\Sgn_{-,-,-} 
\right) 
\]
to $[P_F] \otimes_F \ld \tH^3_? (U \times_B U , 2)_{-,-}$. Obviously, $\res_1$
is trivial. For $\res_2$, observe that 
\[
\res_2 = \res_3 \quad \text{on} \quad
\ld \tH^4_?(U \times_B U^2 , 3)^\sgn_{-,-,-} \; .
\]
Therefore, $\res_2$ vanishes on 
\[
\kerr \left( \res_3 \tei
\tH^4_?(U \times_B U^2 , 3)^\sgn_{-,-,-}
\right) \; .
\]
Because of the definition of $\partial_1$,
it remains to observe that
\[
\res_3 \circ \symm_{1,2} (x) = 0 \Longleftrightarrow
\res_3 (x) = 0
\]
for all $x \in \kerr (\ld \tH^4_?(U \times_B U^2 , 3)^\sgn_{-,-,-}
\to \ld \tH^4_?(U^3 , 3)^\Sgn_{-,-,-})$.
\end{Proof}

This proves Theorem~\ref{9bA} and Main Theorem~\ref{9bI} for $k=3$. \\

\begin{Rem} \label{11H}
The exact sequence of \ref{11G} should be compared to 
\cite{GL}, Theorem 1.5, which concerns the case when the base $B$ is
the spectrum of a field $K$. The result is stronger than \ref{11G} in that
\begin{enumerate}
\item[(a)] it contains a statement about the cohomology
of the complex $B(\Eh,3)^\bullet$ of \cite{GL}, (9),
which is a continuation of (the analogue of)
\[
\Bl_{3,\Mh}
\stackrel{d_3}{\pfeil} \Bl_{2,\Mh} \otimes_\BQ \, \Eh(K) 
\]
to the right.
\item[(b)] the third Bloch group of loc.~cit.\ is described in terms of
generators and relations. 
\end{enumerate}
\end{Rem}

We conclude the subsection with a summary of the situation
for $k=3$, when the base $B$ is of arithmetic nature. When $B$ is the
spectrum of a number field $K$, the result is implied by \cite{GL},
Theorem~1.5:

\begin{Thm} \label{11I}
Let $B$ be a smooth, separated and connected
scheme of finite type over a number field $K$, or a ring of integers $O_K$, 
$\Eh$ an elliptic curve
over $B$, and $P \subset \Eh(B)$ a subset satisfying the disjointness
property $(DP)$. 
In the $\Q$-vector space
\[
\Q [P - \{0\}] \; ,
\]
consider the subspace $B_{3,P}^{\sharp}$ of those divisors 
$\sum_{\alpha} \lambda_{\alpha} (s_{\alpha})$ satisfying the 
following conditions:
\begin{enumerate}
\item[(i)] For any homomorphism $X : P_\BQ \to \Q$, one has
$$\sum_{\alpha} \lambda_{\alpha} X (s_{\alpha})^3 = 0 
\quad \mbox{in} \; \Q \; .$$
\item[(ii)] For any closed point $b$ of the generic fibre
of $B$ with residue field $K(b)$ and fibre
$\Eh_b$ of $\Eh$, and any homomorphism $X : P_\BQ \to \Q$, one has
$$\sum_{\alpha} \lambda_{\alpha} X (s_{\alpha}) h_v (s_{\alpha,b}) = 0$$ 
for any (finite or infinite) place $v$ of $K(b)$. Here, $h_v$ denotes the 
local N\'eron height function.
\end{enumerate}
(a) The map $(s) \mapsto \{ s \}_{3}$ induces a morphism
\[
B_{3,P}^{\sharp} \pfeil \kerr (d_{3,\Mh}) \isoto 
                        H^2_\Mh (\Eh^{(1)} , 2)^\sgn \; .
\]
(b) For any $\BC$-valued point $b$ of $B$, and any
\[
S = \sum_{\alpha} \lambda_{\alpha} (s_{\alpha}) \in B_{3,P}^{\sharp} \; ,
\]
the value of the regulator to absolute Hodge cohomology
\[
H^2_{\abs} (\Eh^{(1)}_b , \R (2))^{\sgn} = 
H^1 (\Eh_b (\C) , 2\pi i \R)
\]
(see Remark~\ref{2J}~(a)) on $S$ is given by
\[
3 \cdot 
\sum_{\alpha} \lambda_{\alpha} G_{\Eh_b,3} (s_{\alpha,b}) \; ,
\]
for the functions $G_{\Eh_b,3}$ of \ref{2H}.
\end{Thm}

\begin{Proof}
Observe that
(i) is a reformulation of the condition
$$\sum_{\alpha} \lambda_{\alpha} X (s_{\alpha}) \cdot \{ s_{\alpha} \}_2
\in \kerr (d_{2,\Mh}) \; ,$$
while (ii) is equivalent to
$$\Eis^0_{\Mh} 
\left( \sum_{\alpha} \lambda_{\alpha} X (s_{\alpha}) \cdot 
   \{ s_{\alpha} \}_2 \right) 
    = 0 \quad \mbox{in} \; \Oh^*(B) \otimes_{\Z} \Q$$
for any $X: \Eh(B) \to \Q\,$; note that a regular function on
$B$ is equal to the constant function $1$ if and only if it takes the 
value $1$ at any closed point of the generic fibre of $B$. The statement 
about the regulator is \ref{9bE}.
\end{Proof}

\begin{Rem} \label{11K} 
In order to
check whether or not a given function on $B$ is constant, it suffices to
consider a dense subset of the generic fibre of $B$. Therefore, 
one is reduced to checking condition (ii) of \ref{11I} for a set of
closed points of the generic fibre, which lies dense inside $B$. 
\end{Rem} 

\newpage
\subsection{The torsion case. II} \label{12}

The modest aim of this short subsection is to show that the
two constructions of the Eisenstein symbol on torsion
given in Subsection~\ref{8} and in \ref{9bC} yield the same result.
Readers mainly interested in the proofs of the main results
may therefore decide to skip this subsection. \\

Recall the situation: $s \in \tEh (B)$ is torsion in $\Eh(B)$, and $k \ge 2$. 
We define
\[
U_s := \tEh -  s(B) \; ,
\]
which is an element of $\Ch_{\pi, (-)}$.
Set $P := \{ i,s \}$.
First, we have the following variant of \ref{8A} and \ref{8B}:

\begin{Prop} \label{12A}
There is an isomorphism,
canonical up to the choice of generator of
$F [U_{s,\infty} (B) ]_{(-)}$,
\[
H^{k+1}_{?} (U_s^{k} ,k)^{\Sgn}_{(-, \ldots , -)} \isoto
\bigoplus_{r=0}^{k} 
H^{r+1}_{?} (\Eh^{r} ,r)^{\Sgn}_{-, \ldots , -} \; .
\]
It is compatible with the residue.
\end{Prop}

Denote by $\alpha$ the natural map
\[
H^{k+1}_{?} (\Eh \times_B U_s^{k-1} ,k)^{\sgn}_{-, (-, \ldots , -)} \pfeil
H^{k+1}_{?} (U_s^{k} ,k)^{\Sgn}_{(-, \ldots , -)} \; .
\]

\begin{Cor} \label{12B}
There is an isomorphism,
canonical up to the choice of generator of
$F [U_{s,\infty} (B) ]_{(-)}$,
\[
\kerr (\alpha) \isoto 
\bigoplus_{r=2}^{k} 
H^{r+1}_{?} (\Eh^2 \times_B \Eh^{r-2} ,r)^{+,\sgn}_{-, \ldots , -} \; .
\]
It is compatible with the residue.
\end{Cor}

\begin{Proof}
This follows from \ref{6Eb}.
\end{Proof}

Denote by $\can$ the canonical map from $\kerr (\alpha)$
to the Bloch group $\Bl_{k,P,?}$.
The compatibility of the constructions is a consequence of the
following:

\begin{Prop} \label{12C}
The diagram
\[
\vcenter{\xymatrix@R-10pt{ 
\kerr (\alpha)
\ar[r]^-{\cong}_-{\ref{12B}} \ar[d]_{\can} &
\bigoplus_{r=2}^{k} 
H^{r+1}_{?} (\Eh^2 \times_B \Eh^{r-2} ,r)^{+,\sgn}_{-, \ldots , -}
\ar@{>>}[d] \\
\Bl_{k,P,?} &
H^{k-1}_{?} (\Eh^{(k-2)} , k-1)^{\sgn}
\ar[l]^-{\ref{9bC}}_-{\cong}  
\\}}    
\]
is commutative.
\end{Prop}

\begin{Proof}
Assume that $[N] s = 0$.
All the maps in the diagram are compatible with the action of $[N+1]$.
This shows already that $\can$ is trivial on elements of $\kerr (\alpha)$
with trivial component in
\[
H^{k+1}_{?} (\Eh^2 \times_B \Eh^{k-2} ,k)^{+,\sgn}_{-, \ldots , -} 
\]
under the isomorphism of \ref{12B}. On the part of $\kerr (\alpha)$
corresponding to the direct summand
\[
H^{k+1}_{?} (\Eh^2 \times_B \Eh^{k-2} ,k)^{+,\sgn}_{-, \ldots , -} \; , 
\]
$\can$ gives the right map. This is not quite obvious, because the
projection maps
\[
\Sigma \quad \mbox{and} \quad q: \Eh^k \pfeil \Eh^{k-1}
\]
used in Section~\ref{II} and Subsection~\ref{8} respectively, are not the same. 
In fact, we have $\Sigma = q \circ t$, where
\[
t: \Eh^k \pfeil \Eh^k 
\]
is given by
\[
(x_1,x_2,x_3,x_4,\ldots,x_k) \longmapsto 
            (x_1,x_2,x_1+x_3,x_1+x_4,\ldots,x_1+x_k) \; .
\]
By \ref{14Ba}, $t$ acts trivially on
\[
H^\bullet_{?} (\Eh^k , \argstern)_{-, \ldots , -} \; .
\]
\end{Proof}

\newpage
\subsection{Proof of the Main Theorem} \label{15}

This subsection will provide proofs of Theorem~\ref{9bA}, as well as of 
Main Theorem~\ref{9bI}. We will show them in parallel, which
is possible because the analogues of the vanishing hypotheses
made in \ref{9bI} hold in absolute cohomology:

\begin{Prop} \label{15A}
$H^i_{\abs} (\Eh^n , j)_{- , \ldots , -} = 0$ if $i < n$, or if 
$i = n$ and $2j \neq n$.
\end{Prop}

\begin{Proof}
First, observe that
\[
H^{\bullet}_{\abs} (\Eh^n , \argstern) = 
H^{\bullet}_{\abs} (B , \Rh_B (\Eh^n , F (\argstern)_{\Eh^n})) \; ,
\]
from which we conclude that
\[
H^{\bullet}_{\abs} (\Eh^n , \argstern)_{- , \ldots , -} = 
H^{\bullet - n- d (B)}_{\abs} (B , \fH^{\otimes n} (\argstern - n))
\]
since $[-1]$ acts trivially on 
\[
F (1) = \Hh^{-1}_B (\Eh , F (1)) \quad \mbox{and} \quad 
F (0) = \Hh^1_B (\Eh , F(1)) \; .
\]
We thus have
\[
H^i_{\abs} (\Eh^n , j)_{- , \ldots , -} = 
\Ext^{i-n}_{\Sh B} (F (0) , \fH^{\otimes n} (j-n)) \; .
\]
The sheaf $\fH^{\otimes n} (j-n)$ is pure of weight $n-2j$, so
there are no nontrivial morphisms 
\[
F(0) \pfeil \fH^{\otimes n} (j-n) 
\] 
if $2j \ne n$.
\end{Proof}

\quad \\

\begin{Proofof}{Theorem~\ref{9bA} and Main Theorem~\ref{9bI}}
Our main technical tools will be the spectral sequences $^m\tilde{\stern}$
for $m \le k$ of Theorem~\ref{9aC}, or more generally, of
Theorem~\ref{13bG}, and of their direct limits over 
$U \in \Ch_{\pi , P , -}$. By \ref{15A}, and by assumption
respectively, we have:
\[
H^i_?(\Eh^m , m)_{- , \ldots , -} = 0
\]
for $2 \le m \le k-2$ and $-k+2m+2 \le i \le m$:
\begin{eqnarray*}
H^i_? (\Eh^2 , 2)_{-,-} & = & 0 \; , \quad -k + 6 \le i \le 2 \; , \\
H^i_? (\Eh^3 , 3)_{-,-,-} & = & 0 \; , \quad -k + 8 \le i \le 3 \; ,\\
& \vdots & \\
H^{k-2}_? (\Eh^{k-2} , k-2)_{- , \ldots , -} & = & 0 \; .
\end{eqnarray*}

Let us consider the following two claims:
\begin{enumerate}
\item[(1)] Let $U \in \Ch_{\pi,-}$. Then
\[
\tH^i_?(U^m , m)_{- , \ldots , -} = 0
\]
for $2 \le m \le k-2$ and $-k+2m+2 \le i \le m$:
\begin{eqnarray*}
\tH^i_? (U^2 , 2)_{-,-} & = & 0 \; , \quad -k + 6 \le i \le 2 \; , \\
\tH^i_? (U^3 , 3)_{-,-,-} & = & 0 \; , \quad -k + 8 \le i \le 3 \; ,\\
& \vdots & \\
\tH^{k-2}_? (U^{k-2} , k-2)_{- , \ldots , -} & = & 0 \; .
\end{eqnarray*}
\item[(2)] Let $U \in \Ch_{\pi,-}$. Then the restriction
\[
H^{2m-k+1}_?(\Eh^2 \times_B \Eh^{m-2} , m)_{- , \ldots , -}^{+,\sgn} \pfeil
\tH^{2m-k+1}_?(U^2 \times_B U^{m-2} , m)_{- , \ldots , -}^{+,\sgn}
\]
is surjective for $2 \le m \le k-1$.
\end{enumerate}
Before proving these claims, let us show that they imply the theorems.
We have really seen
all the arguments in the proof of the case $k=3$:
\begin{enumerate}
\item[(i)] Show that $\ker(d_k)$ is already contained in the term
$^k\tE^{-k,1}_\infty$.
\item[(ii)] Show that the terms $^k\tE^{p,q}_\infty$, for $p+q = -(k-1)$ and
$p \ne -k$ are trivial.
\end{enumerate}

As for (i), observe that because of (1), \emph{almost} all the
targets of the differentials on the $^k\tE^{-k,1}_r$, for $r \ge 2$, 
are trivial, except one: the differential
\[
^k\partial_{k-1}: {^k\tE^{-k,1}_{k-1}} \pfeil {^k\tE^{-1,3-k}_{k-1}} \; .
\]
Recall that 
\[
{^k\tE^{-1,3-k}_1} = 
[P_F] \otimes_F 
          \ld \tH^k_? (U \times_B \Eh^{k-2} , k-1)_{-,\ldots,-}^\sgn \; .
\]
Thus, restriction gives a map
\[
\alpha: {^k\tE^{-1,3-k}_1} \pfeil
[P_F] \otimes_F 
           \ld \tH^k_? (U \times_B U^{k-2} , k-1)_{-,\ldots,-}^\sgn \; ,
\]
which induces a map from ${^k\tE^{-1,3-k}_{k-1}}$ to a certain subquotient of
\[
[P_F] \otimes_F 
           \ld \tH^k_? (U \times_B U^{k-2} , k-1)_{-,\ldots,-}^\sgn \; .
\]
The composition $\beta$ of this map and the
differential $^k\partial_{k-1}$ is induced by the residue map
\[
\ld \tH^{k+1}_? (U^2 \times_B U^{k-2} , k)_{-,\ldots,-}^{+,\sgn} \longto
[P_F] \otimes_F 
           \ld \tH^k_? (U \times_B U^{k-2} , k-1)_{-,\ldots,-}^\sgn 
\]
(in any of the first two coordinate directions). The map $\beta$ is trivial on
$\kerr(d_k)$ for the same reasons as those which appeared in the proof
of \ref{11G}. Step (i) will thus be completed once we have shown that $\alpha$
is injective. For this, we observe that $\alpha$ is just the edge morphism
in the spectral sequence constructed in \ref{13bI}. Injectivity of $\alpha$ is
thus a consequence of the vanishing of its $E^{p,q}_\infty$-terms for
\[
p+q = -k+2 \; \quad p \ne -(k-1) \; ,
\]
which in turn follows from (1).\\

Because of (2), the
terms $^k\tE^{p,q}_\infty$, for $p+q = -(k-1)$ and
$-(k-1) \le p \le -2$, are subquotients of terms of the form
\[
[P_F]^{\otimes {k-m}} \otimes
H^{2m-k+1}_?(\Eh^m , m)_{- , \ldots , -} 
\]
for some $m \le k-1$. (Actually, because of (1), these terms are
trivial except for $p = -(k-1)$). The same holds for 
the term $^k\tE^{0,-(k-1)}_\infty$ 
(for trivial reasons), and for $^k\tE^{-1,-(k-2)}_\infty$ since the
restriction
\[
H^{k-1}_? (\Eh \times \Eh^{k-2}, k-1)_{- , \ldots , -} \pfeil
H^{k-1}_? (U \times \Eh^{k-2}, k-1)_{- , \ldots , -}
\]
is surjective. Indeed, the residue sequence shows that the cokernel of
this map injects into 
\[
F [U_\infty (B)]_- \otimes_F H^{k-2}_? (\Eh^{k-2}, k-2)_{- , \ldots , -} \; ,
\]
which by assumption is trivial.
The edge morphisms 
\[
^k\tE^{p,q}_\infty \pfeil 
H^{k+1}_? (\Eh^2 \times_B \Eh^{k-2} , k)_{- , \ldots , -}^{+,\sgn}
\]
are induced by the exterior cup products, followed by the projections
onto the $({+,\sgn})$-eigenspaces. By Theorem~\ref{14A}, these maps 
are trivial. This shows (ii).\\

Now for the proofs of (1) and (2). First, one shows as above 
that the restriction
\[
H^{2m-k}_? (\Eh \times \Eh^{m-2}, m-1)_{- , \ldots , -}^\sgn \pfeil
H^{2m-k}_? (U \times \Eh^{m-2}, m-1)_{- , \ldots , -}
\]
is surjective. By assumption, these terms are thus trivial.
The spectral sequence $^m\tilde{\stern}$ then shows that (1) implies
(2).\\

Finally, (1) will be proved by induction on $k$, the claim being trivial
if $k \le 3$. For fixed $k \ge 4$, we know already that  
\[
\tH^i_?(U^m , m)_{- , \ldots , -} = 0
\]
for $2 \le m \le k-3$ and $-k+2m+3 \le i \le m$:
\begin{eqnarray*}
\tH^i_? (U^2 , 2)_{-,-} & = & 0 \; , \quad -k + 7 \le i \le 2 \; , \\
\tH^i_? (U^3 , 3)_{-,-,-} & = & 0 \; , \quad -k + 9 \le i \le 3 \; ,\\
& \vdots & \\
\tH^{k-3}_? (U^{k-3} , k-3)_{- , \ldots , -} & = & 0 \; .
\end{eqnarray*}
For all $m$ between $2$ and $k-2$, we have to show in addition that
\[
\tH^{-k+2m+2}_? (U^m , m)_{- , \ldots , -} = 0 \; ,
\]
and this will in turn be achieved by induction on $m$. Assume the vanishing
has been established for all $m = 1,2,\ldots,M-1$. The
spectral sequence $^M\tilde{\stern}$ of Theorem~\ref{13bG}, together with
our assumptions on the vanishing of the
\[
H^i_?(\Eh^m , m)_{- , \ldots , -}
\]
shows then that 
\[
\tH^{-k+2M+2}_? (U^M , M)_{- , \ldots , -} = 0 \; ,
\]
thereby proving (1), and hence completing the proof.
\end{Proofof}

\newpage
\subsection{Bloch groups and Eisenstein symbol \\
in absolute cohomology} \label{13c}

By construction, the Eisenstein symbol is compatible with the regulators. 
The purpose of this subsection is to connect our {\em geometric} construction 
of $\Eis^{k-2}$ to the {\em sheaf theoretical} one of 
\cite{W5} sketched in Subsection~\ref{2} (Theorem~\ref{13cG}). 
In particular, we get 
the desired relation to Kronecker double series (Theorem~\ref{9bE}).\\

As always, we fix $P \subset \Eh (B)$ satisfying $(DP)$.
Our first aim is a sheaf theoretical interpretation of the groups
occurring in the definition of
$\Bl_{k,P,\abs}$ for $k \ge 2$ (Proposition~\ref{13cDa}). 
Recall the definition and basic properties
of the sheaves $\Gh^{(n)}_U \in \Sh^{s,W} \!\! B$, for
$U \in \Ch_{\pi , -}$ and $n \ge 0$ (\ref{4A}--\ref{4Dc}).

\begin{Def} \label{13cA}
For $U \in \Ch_{\pi , -}$, define
\[
\CF_U := \Hh^0_B (U^2 , F (2))^+_{- , -} = 
\Hh^{2+d (B)}_B (U^2 , F (2)_{U^2})^+_{- , -} \; .
\]
\end{Def}

\begin{Prop} \label{13cB}
(a) There is a canonical isomorphism
\[
\CF_U \isoto \bigwedge^2 \Gh^{(1)}_U \; .
\]
\noindent (b) $\CF_U$ has a weight filtration, and
\[
W_{-2} \CF_U = \bigwedge^2 \fH = F(1) \; .
\]
The identification of $\bigwedge^2 \fH$ and $F(1)$ is 
induced by the Poincar\'e pairing
\[
\langle \argdot , \argdot \rangle : (\fH(-1)) \otimes_F \fH \pfeil F(0)
\]
(see Subsection~\ref{2}).
\end{Prop}

\begin{Proof}
(a) follows from the K\"unneth formula and the definitions. 
(b) is a consequence of \ref{4Db}.
\end{Proof}

\begin{Cor} \label{13cC}
Let $k \ge 2$. Then the sheaf $\CF_U \otimes_F \Gh^{(k-2)}_U$ has a 
weight filtration, and
\[
W_{-k} \left( \CF_U \otimes_F \Gh^{(k-2)}_U \right) = 
\Sym^{k-2} \fH (1) \; .
\]
\end{Cor}

\begin{Proof}
This is \ref{4Dc}.
\end{Proof}

\begin{Prop} \label{13cDa}
(a) Let $k \ge 2$. The identity
\[
\Gh^{(1)}_U \otimes_F \bigl( \Gh^{(1)}_U)^{\otimes (k-1)} \bigr) 
\stackrel{=}{\pfeil}
\bigl( \Gh^{(1)}_U)^{\otimes 2} \bigr) \otimes_F
                       \bigl( \Gh^{(1)}_U)^{\otimes (k-2)} \bigr) \; ,
\]
together with the canonical isomorphisms of \ref{13cB}, \ref{4B}, and \ref{4Db},
induces a natural morphism
\[
p: \fH \otimes_F \Gh^{(k-1)}_U \pfeil
\CF_U \otimes_F \Gh^{(k-2)}_U \; .
\]
(b) The map $p$ on 
$\Ext^1_{\SB} (F(0) , \argdot)$-niveau induced by the morphism $p$ 
equals the composition of the following two maps: first, the
restriction 
\[
H^{k+1}_\abs (\Eh \times_B U^{k-1} , k)^\sgn_{-,-,\ldots,-} \pfeil
H^{k+1}_\abs(U \times_B U^{k-1},k)^\sgn_{-,-,\ldots,-}
\]
where $\sgn$ refers to the action of $\fS_{k-1}$; 
second, the symmetrization with respect to the first two coordinates
on $H^{k+1}_\abs(U \times_B U^{k-1},k)^\sgn_{-,-,\ldots,-}$:
\[
H^{k+1}_\abs(U \times_B U^{k-1},k)^\sgn_{-,-,\ldots,-} \pfeil
H^{k+1}_\abs(U^2 \times_B U^{k-2},k)^{+,\sgn}_{-,\ldots,-} \; .
\]
$p$ is therefore compatible with the map denoted $p$ in Definition~\ref{9aD}.
\end{Prop}

\begin{Proof}
This follows from \ref{4Da} and \ref{4K}~(a).
\end{Proof}

We now use the machinery developed in \cite{W5}, Sections~2 and 3. 
Think of an element
\[
\E \in H^{k+1}_{\abs} (U^2 \times_B U^{k-2} , k)^{+,\sgn}_{- , \ldots , -} = 
\Ext^1_{\SB} (F(0) , \CF_U \otimes_F \Gh^{(k-2)}_U)
\]
as a sheaf, together with fixed morphisms

\begin{eqnarray*}
\E & \longonto & F(0) \quad \mbox{and} \\
\CF_U \otimes_F \Gh^{(k-2)}_U & \longinto & \E \; .
\end{eqnarray*}
Passing to weight graded objects, 
and using \ref{13cC}, we see that the second of these morphisms defines
\[
y : \Gr^W_{\bullet} \E  \longonto  
\Gr^W_{-k} \left( \CF_U \otimes_F \Gh^{(k-2)}_U \right)
= \Sym^{k-2} \fH (1) \; .
\]
We would like to define 
\[
x : F(0) \longinto  \Gr^W_{\bullet} \E
\]
in a \emph{canonical} way from $\E$.
Since weights $0$ occur in $\CF_U \otimes_F \Gh^{(k-2)}_U$, this
is not a priori possible unless $\E$ is of a special shape. 
So assume that $\E$ is actually contained in the image of the morphism
$p$ defined in \ref{13cDa}.
Since $\fH \otimes_F \Gh^{(k-1)}_U$ is of strictly negative weights,
there is a canonical choice of
\[
x : F(0) \longinto  \Gr^W_{\bullet} \E \; .
\]
In the terminology of \cite{W5}, 2.3, we have defined a ``coefficient''
\[
cl (\E)_k \in \Gamma (B , \Lie^{\vee}_B \otimes_F \Sym^{k-2} \fH (1)) \; ,
\]
where as in loc.\ cit., 3.3,
we denote by $\Lie_B$ the Lie algebra of the 
pro-unipotent part of the Tannakian dual of $\Sh^{s,W} \!\! B$.\\

We thus get a morphism
\[
cl_{k,U} : \imm (p) 
\longrightarrow \Gamma (B ,  \Lie^{\vee}_B \otimes_F \Sym^{k-2} \fH (1)) \; ,
\]
where $p$ denotes the morphism
\[
H^{k+1}_\abs (\Eh \times_B U^{k-1} , k)^\sgn_{-,-,\ldots,-} \pfeil
H^{k+1}_{\abs} (U^2 \times_B U^{k-2} , k)^{+,\sgn}_{- , \ldots , -}
\]
of Proposition~\ref{13cDa} and Definition~\ref{9aD}.

\begin{Prop} \label{13cE}
For $k \ge 2$, the morphisms $cl_{k,U}$ induce a morphism
\[
cl_k : \Bl_{k,P,\abs} \longrightarrow 
\Gamma (B , \Lie^{\vee}_B \otimes_F \Sym^{k-2} \fH (1)) \; .
\]
\end{Prop}

\begin{Proof}
We need to show that every element $\E$ in the image of the cup product
\[
\cup_j: H^1_\abs (U , 1)_- \otimes H^k_\abs (U^{k-1} , k-1)_{- \ldots,-} \pfeil
H^{k+1}_{\abs} (U^2 \times_B U^{k-2} , k)^{+,\sgn}_{- , \ldots , -} 
\]
in the $j$-th coordinate direction has a trivial coefficient
($j = 1, \ldots, k$). For this, a sheaf theoretical interpretation
of this cup product is required. Recall the exact sequence
\[
0 \pfeil H^1_\abs (U , 1)_- \stackrel{\res}{\pfeil}
F [U_\infty (B)] _- \stackrel{\alpha}{\pfeil} 
H^2_\abs (\Eh , 1)_- = \Ext^1_{\SB} (F(0),\fH) \; .
\]
Given $x \in F [U_\infty (B)] _-$, the associated extension $\alpha (x)$
is obtained by pulling back the exact sequence
\[
0 \pfeil \fH \pfeil \Gh^{(1)}_U \pfeil
F [U_\infty (B)] _- (0) \pfeil 0 
\]
of \ref{4Db} via $x$, which we consider as a morphism
\[
F(0) \pfeil F [U_\infty (B)] _- (0) \; .
\]
For $x \in H^1_\abs (U , 1)_-$, this extension is trivial, and because
of weight reasons, there is then a unique splitting
\[
\beta_x: F(0) \pfeil \Gh^{(1)}_U \; .
\]
This in turn induces a morphism
\[
\beta'_x : \bigl( \Gh^{(1)}_U \bigr)^{\otimes (k-1)} \pfeil 
\bigl( \Gh^{(1)}_U \bigr)^{\otimes k} 
\longonto \CF_U \otimes_F \Gh^{(k-2)}_U \; ;
\]
for the first of the two arrows, 
we have used $\beta_x$ in the $j$-th coordinate. 
We then have the following compatibility: the image of 
\[
x \otimes y \in 
H^1_\abs (U , 1)_- \otimes H^k_\abs (U^{k-1} , k-1)_{- \ldots,-} 
\]
under $\cup_j$ equals the extension in
\[
H^{k+1}_{\abs} (U^2 \times_B U^{k-2} , k)^{+,\sgn}_{- , \ldots , -} = 
\Ext^1_{\SB} (F(0) , \CF_U \otimes_F \Gh^{(k-2)}_U)
\]
obtained by pushing out 
\[
y \in H^k_\abs (U^{k-1} , k-1)_{- \ldots,-} = 
\Ext^1_{\SB} (F(0) , (\Gh^{(1)}_U)^{\otimes (k-1)} )
\]
via $\beta'_x$. Since the minimal weight contained in $y$ is $-(k-1)$, this
extension has a trivial coefficient $cl(x \otimes y)_k$.
\end{Proof}

There is a canonical isomorphism
\[
cl_1: H^2_\abs(\Eh,1)_- = \Ext^1_{\SB} (F(0),\fH) \isoto
\Gamma (B , \Lie^{\vee}_B \otimes_F \fH) 
\]
(see \cite{W5}, Remark~3.3~(a)). Proposition~\ref{13cE} thus extends to the
case $k=1$. Now recall the elements 
\[
\{ s \}_k \in \Gamma (B , \Lie^{\vee}_B \otimes_F \Sym^{k-2} \fH (1)) 
\] 
defined in \cite{W5}, 3.3 and recalled in \ref{2J}~(c). 

\begin{Prop} \label{13cF} 
Let $s \in P \cap \tEh (B)$.
\begin{enumerate}
\item[(a)] The morphism $cl_1$ maps 
$\{ s \}_{1,\abs} \in \Bl_{1,P,\abs}$ to the element 
\[
\{ s \}_1 \in \Gamma (B , \Lie^{\vee}_B \otimes_F \fH) \; .
\]
\item[(b)] For $k \ge 2$, the morphism $cl_k$ maps 
$\{ s \}_{k,\abs} \in \Bl_{k,P,\abs}$ to $\frac{k!}{k-1}$ times the element
\[
\{ s \}_k \in \Gamma (B , \Lie^{\vee}_B \otimes_F \Sym^{k-2} \fH (1)) \; .
\]
\end{enumerate}
\end{Prop}

\begin{Proof} 
The factor 
\[
\frac{k!}{k-1} = (k-1)! \cdot \frac{k}{k-1}
\]
in (b) is explained as follows:
\begin{enumerate}
\item[(1)] The 
different normalizations
of the epimorphism
\[
\fH \otimes_F \Sym^{k-1} \fH \longonto 
\Sym^{k-2} \fH (1) 
\]
account for the factor $k/(k-1)$
(see \ref{2J}~(c)). 
\item[(2)] The factor $(k-1)!$ comes from \ref{5F} -- observe that
the constructions of 
\cite{W5} were performed
in the ``coordinates'' $\kappa$, while our definition of
$\{ s \}_{k,\abs}$ uses $\eta$.
\end{enumerate}
(a) follows from 
\ref{9aG}~(a) and \cite{W5}, Remark 3.3~(a).
\end{Proof}

Our next task is to compare the differential
$d_{k,\abs}$ on $\Bl_{k,P,\abs}$ 
to the differential
\[
d \otimes \id : \Gamma (B , \Lie^{\vee}_B \otimes_F \argdot) 
\longrightarrow \Gamma (B , \bigwedge^2 \Lie^{\vee}_B \otimes_F \argdot)
\]
of \cite{W5}, Section 2.
Recall (Definition~\ref{9aE}) that for $k \ge 3$,
\[
d_k: \Bl_{k,P,\abs} \pfeil
\Bl_{k-1,P,\abs} \otimes_F [P_F] \subset 
                       \Bl_{k-1,P,\abs} \otimes_F \Bl_{1,P,\abs}
\]
is induced by the
limit over $U \in \Ch_{\pi , P , -}$ of the residue maps
\[
H^{k+1}_{\abs} (U^2 \times_B U^{k-2} , k)^{+,\sgn}_{- , \ldots , -} \!\pfeil\!
H^k_{\abs} (U^2 \times_B U^{k-3} , k-1)^{+,\sgn}_{- , \ldots , -} 
                  \otimes_F F [U_{\infty} (B)]_- 
\]
For $k=2$, 
\[ 
d_2: \Bl_{2,P,\abs} \pfeil 
\Sym^2 \Bl_{1,P,\abs} = \Sym^2 H^2_\abs(\Eh, 1)_-
\]
is induced by the
limit over $U \in \Ch_{\pi , P , -}$ of the symmetrization of 
the composition of the residue maps
\[
H^{3}_\abs(\Eh \times_B U,2)_{-, -}
\pfeil H^2_\abs(\Eh, 1)_- \otimes_F F[U_\infty (B)]_- \; ,
\]
and of the Abel--Jacobi map
\[
H^2_\abs(\Eh, 1)_- \otimes_F F[U_\infty (B)]_-
\stackrel{[\argdot]}{\pfeil} H^2_\abs(\Eh, 1)_-^{\otimes 2} \; .
\]
The sheaf theoretical interpretation of the residue 
looks as follows:

\begin{Prop} \label{13cD}
(a) Let $k \ge 3$. 
The residue map
\[
\Gh^{(k-2)}_U \longrightarrow 
\Gh^{(k-3)}_U \otimes_F F [U_{\infty} (B)]_- (0)
\]
of Subsection~\ref{4} induces a morphism
\[
\CF_U \otimes_F \Gh^{(k-2)}_U \longrightarrow 
\CF_U \otimes_F \Gh^{(k-3)}_U \otimes_F F [U_{\infty} (B)]_- (0) \; .
\]
\noindent (b) The map on 
$\Ext^1_{\SB} (F(0) , \argdot)$-niveau induced by the morphism in (a)
equals the residue
\[
H^{k+1}_{\abs} (U^2 \times_B U^{k-2} , k)^{+,\sgn}_{- , \ldots , -} \!\pfeil\!
H^k_{\abs} (U^2 \times_B U^{k-3} , k-1)^{+,\sgn}_{- , \ldots , -} 
                  \otimes_F F [U_{\infty} (B)]_- 
\]
\noindent (c) The residue map
\[
\Gh^{(1)}_U \longrightarrow F [U_{\infty} (B)]_- (0)
\]
of Subsection~\ref{4} induces a morphism
\[
\fH \otimes_F \Gh^{(1)}_U \longrightarrow 
\fH \otimes_F F [U_{\infty} (B)]_- (0) \; . 
\]
\noindent (d) The map on
$\Ext^1_{\SB} (F(0) , \argdot)$-niveau induced by the morphism in (c)
equals the residue
\[
H^{3}_\abs(\Eh \times_B U,2)_{-, -}
\pfeil H^2_\abs(\Eh, 1)_- \otimes_F F[U_\infty (B)]_- \; .
\]
\end{Prop}

\begin{Proof}
This follows from \ref{4Da} and \ref{4K}~(a).
\end{Proof}

Denote by
\[
\pr_k : \Lie^{\vee}_B \otimes_F \Sym^{k-3} \fH (1) \otimes_F
\Lie^{\vee}_B \otimes_F \fH 
\longrightarrow \bigwedge^2 
\Lie^{\vee}_B \otimes_F \Sym^{k-2} \fH (1)
\]
the morphism $\bigwedge \otimes$ mult. 

\begin{Thm} \label{13cG}
Write $S^l$ for the sheaf $\Sym^l \fH$. 
($S^{-1} = \fH^{\vee}$.) Let $k \ge 2$. The following diagram commutes:
\[
\vcenter{\xymatrix@R-10pt{ 
\Bl_{k,P,\abs} \ar[r]^-{cl_k} \ar[dd]_{d_{k,\abs}} &
\Gamma (B , \Lie^{\vee}_B \otimes_F S^{k-2} (1))
\ar[d]^{d \otimes \id} \\
&
\Gamma (B , \bigwedge^2 \Lie^{\vee}_B \otimes_F S^{k-2} (1)) \\
\Bl_{k-1 , P, \abs} \otimes_{F} \Bl_{1 , P, \abs} 
                 \ar[r]^-{cl_{k-1} \otimes cl_1} &
\; \Gamma (B , \Lie^{\vee}_B \otimes_F S^{k-3} (1)) \otimes_F 
\Gamma (B , \Lie^{\vee}_B \otimes_F \fH) \ar[u]_{(k-1) \cdot \pr_k} \\}}
\]
\end{Thm}

\begin{Rem} \label{13cH}
(a) For the elements $\{ s \}_{k,\abs}$ of $\Bl_{k,P,\abs}$, 
Theorem~\ref{13cG} follows
directly from \ref{13cF}, \ref{9aI}, and \cite{W5}, Thm.~3.4. \\[0.2cm]
(b) Unfortunately, \cite{W5}, Thm.~3.4 is
stated incorrectly. For $k \ge 3$, the factor
\[
\frac{k-1}{k}
\]
of loc.~cit.\ should be replaced by 
\[
\frac{(k-1)^2}{k(k-2)} \; . 
\]
In fact, the proof of loc.~cit.\ is correct from its second line onwards.
Since the epimorphisms
\[
\fH \otimes_F \Sym^{k-1} \fH \longonto 
\Sym^{k-2} \fH (1) 
\]
and 
\[
\fH \otimes_F \Sym^{k-2} \fH \longonto 
\Sym^{k-3} \fH (1) 
\]
used in loc.~cit.\ \emph{both} involve factors -- namely 
\[
(k-1)/k \quad \mbox{and} \quad (k-2)/(k-1)
\] 
respectively -- their quotient occurs in the correct 
version. The mistake committed implicitly in the first line of the
proof of loc.~cit.\ occurred because we forgot to take 
into account the second factor. \\[0.2cm]
(c) The same remark applies to \cite{W5}, Lemma~3.5, and to the definition
of the differential $d_k^{\sharp\sharp}$
on page 393 of \cite{W5}, which again has to be modified by
a factor $(k-1)/(k-2)$ when $k \ge 3$. 
All this has no effect on the main results of
loc.~cit.
\end{Rem}

\begin{Proofof}{Theorem~\ref{13cG}}
The proof consists of a faithful imitation of the
one of the equation
\[
d \{ s \}_k^{\sim} = 
\pr_k^{\sim} ( \{ s \}_{k-1}^{\sim}  \otimes \{ s \}_1 )
\]
on page 390 of \cite{W5}. The modification concerns the extension
denoted $[\Delta]$ in loc.~cit., which has to be replaced by $\Gh^{(1)}_U$.
While $\Gr_0^W[\Delta]$ has rank one, $\Gr_0^W \Gh^{(1)}_U$ has a
rank $r$ which will usually be greater than one. The analogue of the Lie
algebra on page 391 of \cite{W5} will thus by generated by elements
$e_1^j$, $e_2^j$, $d_1^j$, and $d_2^j$, for $1 \le j \le r$, where nontrivial
commutator relations occur only between elements having the same superscript
$j$. The nontrivial relations are of the same shape as in loc.~cit.,
page 392. Note that the calculations of \cite{W5} were done using the
``coordinates'' $\kappa$. 
This explains the factor $k-1$, which again
comes from \ref{5F}.
\end{Proofof}

Theorem \ref{13cG} applies to all the cohomology theories of Subsection~\ref{1}. 
One could of course imagine applications to other theories, e.g., a $p$-adic 
one. \ref{13cG} connects the construction of the earlier paragraphs (which 
were possible because of the existence of Grothendieck's functors and of 
a ``formalism of weights'') to the Tannakian construction of 
\cite{W5}, Section~3. The latter is possible 
as soon as the axioms of loc.\ cit., 
3.1 are fulfilled.\\

In the Hodge theoretic setting, we obtain:\\

\begin{Proofof}{Theorem~\ref{9bE}}
This is Theorem~\ref{2I}
and Proposition~\ref{13cF}.
\end{Proofof}

\newpage

%
%

\end{document}